С.В. Курапов
М.В. Давидовский

# АЛГОРИТМИЧЕСКИЕ МЕТОДЫ КОНЕЧНЫХ ДИСКРЕТНЫХ СТРУКТУР

# ТЕОРЕМА О ЧЕТЫРЕХ КРАСКАХ. ТЕОРИЯ, МЕТОДЫ, АЛГОРИТМЫ

(монография)

Проблема четырех красок тесно связана с другими разделами математики и практическими задачами. Известно более 20 ее переформулировок, которые связывают эту проблему с задачами алгебры, статистической механики и задачами планирования. И это тоже характерно для математики: решение задачи, изучаемой из чистого любопытства, оказывается полезным в описании реальных и совершенно различных по своей природе объектов и процессов. Несмотря на опубликованные машинные методы комбинаторного подтверждения гипотезы четырех красок, до сих пор отсутствует чёткое описание механизма раскраски плоского графа четырьмя красками, его естественной сути и его связи с явлением планарности графа. Нужно не только доказать (желательно дедуктивными методами), что любой плоский граф раскрашиваем четырьмя цветами, но нужно и показать, как его раскрасить. В работе рассмотрен подход, основанный на возможности сведения максимально плоского графа к регулярному плоскому кубическому графу с дальнейшей его раскраской. На основании теоремы Тейта-Волынского, вершины максимально плоского графа можно раскрасить четырьмя цветами, если ребра его двойственного кубического графа можно раскрасить тремя цветами. Рассматривая свойства раскрашенного кубического графа, можно показать, что сложение цветов подчиняется законам преобразования группы Клейна четвёртого порядка. Используя это свойство, удается создать линейные алгоритмы раскраски плоского графа.

Для научных работников, преподавателей, студентов и аспирантов высших учебных заведений специализирующихся в области прикладной математики и информатики.

# Содержание

# Введение

История задачи о раскраске вершин плоского графа в четыре цвета довольна драматична. По-видимому, впервые проблема четырёх красок была поставлена немецким математиком А. Мёбиусом (A. Mobius); имеются сведения об устном сообщении на лекциях в 1840г. Однако первое упоминание о гипотезе относят к 1852 г., когда студент Лондонского университетского колледжа Френсис Гутри (Francis Guthrie) изложил эту проблему де Моргану (DeMorgan), а последний описал ее в письме к Гамильтону:

«Мой студент попросил меня сегодня объяснить одну задачу, которая мне не была ранее известна и пока не понятна до конца. Он утверждает, что если любую фигуру разделить любым способом на части и раскрасить их различными красками так, чтобы фигуры, имеющие общий отрезок граничной линии, были окрашены в различные цвета, то всего потребуются четыре краски, но не больше. Мне известен случай, когда требуется четыре краски. Вопрос: нельзя ли придумать случай, когда необходимы, пять или более красок?..»

Гамильтону не удалось придумать карту, для раскраски которой потребовалось бы пять цветов, но он не сумел и доказать, что такой карты не существует. Весть о проблеме четырех красок быстро распространилась по Европе, но, несмотря на все усилия, проблема упорно не поддавалась решению, хотя казалась простой.

Френсис Гутри, раскрашивая от нечего делать карту графств Британии наткнулся на головоломку, которая показалась ему простой, хотя решить ее он так и не сумел. Гутри просто хотел узнать, каково минимальное число красок необходимо взять для раскраски любой мыслимой карты при условии, чтобы никакие две смежные области (имеющие общую границу) не оказались окрашенными в один и тот же цвет.

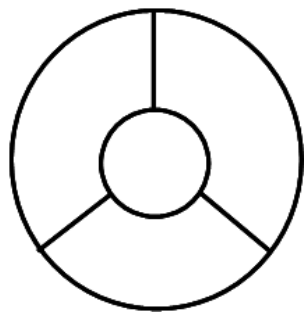

Рис. 1. Карта

Например, для раскрашивания карты, изображенной на рис.1, трех красок недостаточно. Следовательно, для раскрашивания некоторых карт необходимы, по крайней мере, четыре краски. Гутри хотел узнать, окажется ли четырех красок достаточно для раскрашивания всех карт, или для некоторых карт могут потребоваться пять, шесть или больше красок.

Фрэнсис Гутри понял, что карту графств Британии он мог бы раскрасить всего лишь в четыре цвета, причем так, что никакие два соседних графства не оказались бы раскрашенными в один и тот же цвет. Затем он стал размышлять над тем, хватит ли четырех цветов для аналогичной раскраски любой другой карты [14].

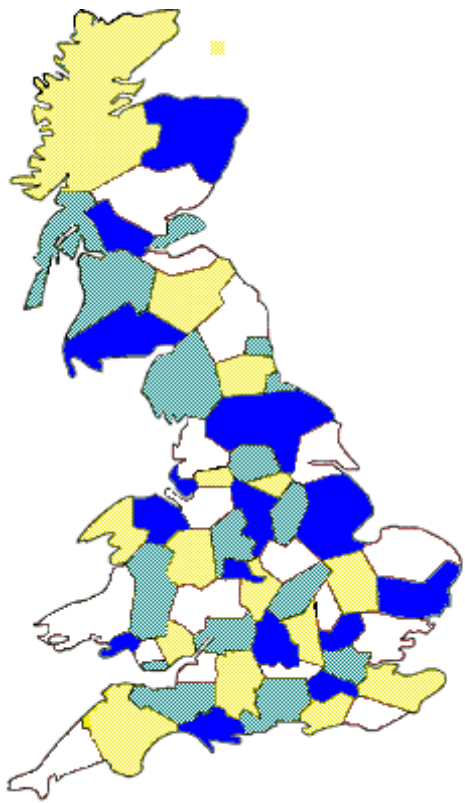

Рис. 2. Карта графств Британии

Любопытно, что Герман Минковский однажды на лекции заявил, что проблема четырех красок не была решена потому, что найти решение пытались только третьеразрядные математики. Но и его собственные усилия в течение нескольких недель не увенчались успехом. «Небеса разгневались на меня за мое высокомерие, — вынужден был признать Минковский. — Мое доказательство также оказалось с изъяном»[14].

Проблему можно сформулировать в терминах графов, если заменить географическую карту графом, расположенным на плоскости (или, что равносильно, на сфере), и тогда вопрос заключается в такой раскраске граней графа наименьшим числом цветов, при которой любые две грани (страны), имеющие общее ребро, не раскрашиваются одинаково. Если построить двойственный граф, вершины которого соответствуют граням данного графа, а ребра соединяют вершины, соответствующие граням, имеющим общее ребро, то задача сводится к раскраске вершин плоского графа так, чтобы любая пара его смежных вершин раскрашивалась разными цветами, и в этом случае раскраска вершин графа называется правильной.

Притягательная сила данной проблемы заключается в доступности ее формулировки. За ее решение до сих пор берутся различные любители, которые приводят свое, очень простое «доказательство» проблемы. Действительно, если воспользоваться элементарными понятиями из теории графов, то формулировка проблемы четырех красок звучит просто.

Необходимо доказать, что все вершины произвольного плоского графа могут быть раскрашены четырьмя красками таким образом, что любые две смежные вершины получат различные цвета.

В 1878 г. знаменитый английский математик А. Кэли (Cayley) выступил перед членами Лондонского математического общества с сообщением о задаче четырех красок, сказав при этом, что сам он вот уже несколько дней (!) как не может ее решить.

В 1879 г. в Шотландии Тэйтом было опубликовано первое решение изложенной выше проблемы. Он свел раскраску вершин исходного графа к раскраске ребер двойственного графа и предположил, что эта задача всегда имеет решение. Тэйт высказал предположение, что всякий плоский кубический граф без перешейков имеет гамильтонов цикл. Если бы ему это удалось доказать, то тем самым была бы решена проблема четырех красок, хотя сам Тэйт и не обнаружил этой зависимости. Однако Татт [17] показал, что это предположение неверно, и привел пример кубического графа с 46 вершинами, не являющегося гамильтоновым (см. рис. 1.3). Так как ребра графа однозначно соответствуют ребрам двойственного графа **H**, то теорему Тэйта можно выразить в эквивалентной форме относительно максимальных плоских графов.

Проблема четырех красок оставалась нерешенной четверть века. Надежда на успех появилась в 1879 году, когда британский математик Альфред Брей Кемпе (Kempe) опубликовал в «American Journal of Mathematics» статью, в которой, по его утверждению, содержалось решение головоломки Гутри. Казалось, Кемпе удалось доказать, что для раскраски любой карты требуется самое большее четыре краски, и тщательное изучение доказательства вроде бы подтверждало его правильность. Кемпе был тотчас же избран членом Королевского общества, а позднее возведен в рыцарское звание за вклад в развитие математики. Это доказательство не подвергалось никаким сомнениям в течение десяти лет. Как остроумно замечает Г.Рингель в своей книге «Проблема раскраски графов»: «Десятилетняя жизнь доказательства Кемпе без какого-либо его подтверждения может служить свидетельством того, что математики того времени были столь же мало склонны читать работы своих коллег, как и в наши дни»[12].

Но в 1890 году лектор Дурхэмского университета Перси Джон Хивуд (Heawood) опубликовал работу, потрясшую математический мир. Через десять лет после того, как Кемпе, казалось бы, решил проблему четырех красок, Хивуд не оставил от его решения камня на камне, показав, где в решении Кемпе была допущена принципиальная ошибка. Единственной хорошей новостью было то, что Хивуду удалось получить оценку для максимального числа красок: оно могло быть равно четырем или пяти, но не более. Проблеме четырех красок Хивуд посвятил всю свою долгую жизнь (он прожил почти 90 лет) и написал много статей, где он обобщал эту проблему на более сложные, чем плоскость, поверхности, но саму ее так и не решил.

Хотя Кемпе, Хивуду и другим так и не удалось решить проблему четырех красок, их попытки внесли большой вклад в новый раздел математики — топологию. В отличие от геометрии, которая занимается изучением точной формы и размеров объекта, топологию интересуют только самые фундаментальные свойства объекта, составляющие его суть.

Математики надеялись, что рассматривая карты через упрощающие линзы топологии, они сумеют постичь самую суть проблемы четырех красок [14].

В 1913 г. Биркгоф ввел понятие неприводимого графа и доказал ряд теорем о свойствах таких графов. Пользуясь этими результатами, американский математик из Филадельфии Франклин в 1920 г. доказал, что гипотеза четырех красок верна для всех плоских графов с числом вершин, не превышающим 25. Это направление получило дальнейшее развитие в работах Эрреры, Винна. В 1926 г. Рейнольдс довел число вершин графов, для которых справедлива гипотеза четырех красок, до 27. В 1936 г. снова Франклин увеличил это число до 31. Затем Винн в 1938 г. довел это число до 35. И хотя в 1970 году Оре и Стемпл увеличили число областей до 39, стало ясно, что дальнейшее его увеличение не приведет к решению проблемы, необходим коренной поворот в исследованиях, который бы дал возможность решить проблему целиком.

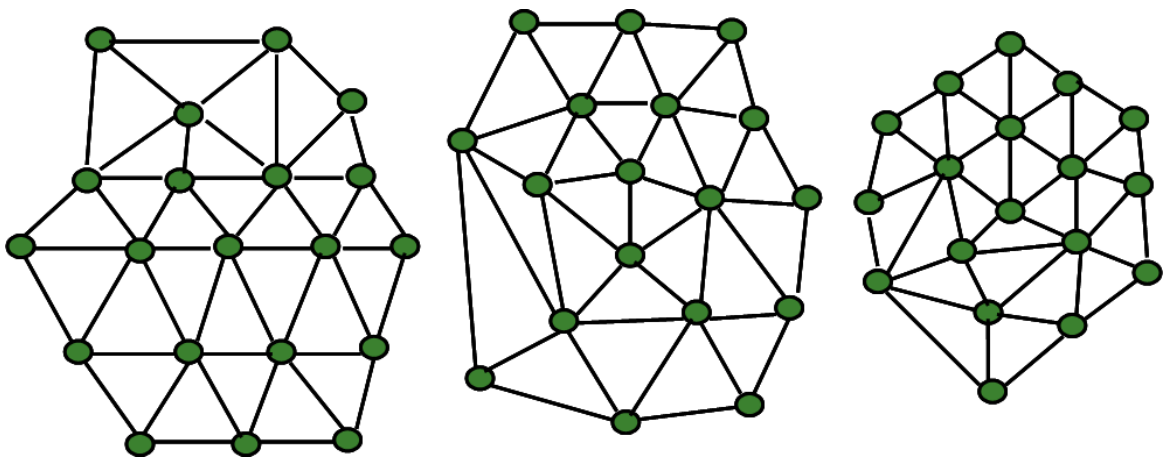

Рис. 3. Неприводимые плоские конфигурации.

Таким поворотом можно считать появление в 1969 г. работы Хееша, в которой он свел вопрос о справедливости гипотезы четырех красок к исследованию достаточно большого так на-зываемого неустранимого множества конфигурации. При этом основной метод, используемый Хеешем, — изобретенный им метод «нейтрализации». Хееш предполагал, что этот метод, если его в достаточной степени развить, приведет к решению проблемы четырех красок.

Для многих математиков эта работа осталась почти незамеченной, остальные были настроены пессимистически относительно перспектив данного направления. Семь лет понадобилось последователям Хееша Аппелю и Хейкену, чтобы этот путь привел их к желаемому результату.

В 1976 г. проблема четырех красок была положительно решена американскими математиками Аппелем. и Хейкеном, работающими в Иллинойском университете. Этому событию предшествовала работа ученых (и не только математиков) различных стран на протяжении более ста лет.

Аппель и Хейкен пришли к успеху после семилетнего кропотливого труда, при этом огромную часть вычислений провели с помощью современных ЭВМ. Но их доказательство настолько сложно и объемно, что за много лет, прошедших с тех пор, не появилось, ни одного сообщения о том, что оно кем-то полностью проверено. Поэтому не прекращаются попытки решить проблему четырех красок с использованием других подходов.

Вольфганг Хейкен и Кеннет Аппель, предложили новый метод, перевернувший традиционные представления о математическом доказательстве. Значительную часть рутинных проверок выполнил компьютер, и это революционное нововведение в сложившуюся практику дедуктивных рассуждений в чистой математике служит основанием для некоторого естественного скептицизма по отношению к данному доказательству и по сей день

В июне 1976 года, затратив 1200 часов машинного времени, Хейкен и Аппель заявили во всеуслышание[18,19], что им удалось проанализировать все 1936 карты и для раскрашивания ни одной из них не требуется более четырех красок. Проблема четырех красок была, наконец решена. Следует особенно подчеркнуть, что решение проблемы четырех красок стало первым математическим доказательством, в котором роль компьютера не сводилась к ускорению вычислений, — компьютер привнес в решение проблемы нечто гораздо большее: его роль была столь значительной, что без компьютера получить доказательство было бы невозможно. Решение проблемы четырех красок с помощью компьютера было выдающимся достижением, но в то же время, оно вызвало у математического сообщества чувство тревоги, так как проверка доказательства в традиционном смысле не представлялась возможной.

Доказательство Аппеля и Хейкена в целом, хотя и принятое математическим сообществом, вызывает до сих пор определенный скептицизм. Для читателя, поверхностно знакомого с математикой, предыдущая фраза должна вызвать изумление: а как же обязательный для математики принцип исключенного третьего, в соответствии с которым утверждение либо справедливо, ли-бо нет? Дело обстоит не так просто. Вот что пишут сами авторы доказательства: «Читатель должен разобраться в 50 страницах текста и диаграмм, 85 страницах с почти 2500 дополнительными диаграммами, 400 страницами микрофишей, содержащими еще диаграммы, а также тысячи отдельных проверок утверждений, сделанных в 24 леммах основного текста». Вдобавок читатель узнает, что проверка некоторых фактов потребовала 1200 часов компьютерного времени, а при проверке вручную потребовалось бы гораздо больше. Статьи устрашающи по стилю и длине, и немногие математики прочли их сколько-нибудь подробно[14].

Прежде, чем опубликовать решение Хейкена и Аппеля на страницах «Illinois Journal of Mathematics», редакторам было необходимо подвергнуть его тщательному рецензированию в каком-то не известном ранее смысле. Традиционное рецензирование было невозможно, поэтому было решено ввести программу Хейкена и Аппеля в независимый компьютер с тем, чтобы убедиться, что результат останется тем же.

Такое нестандартное рецензирование привело в ярость некоторых математиков, утверждавших, будто компьютерная поверка неадекватна, так как не дает гарантии от внезапного отказа в недрах компьютера, который может стать причиной сбоя в логике. Х.П.Ф.Суиннертон-Дайер высказал следующее замечание по поводу компьютерных доказательств: «Когда теорема доказана с помощью компьютера, невозможно изложить доказательство в соответствии с традиционным критерием — так, чтобы достаточно терпеливый читатель смог шаг за шагом повторить доказательство и убедиться в том, что оно верно. Даже если бы кто-нибудь взял на себя труд распечатать все программы и все данные, использованные в доказательстве, нельзя быть уверенным в абсолютно правильной работе компьютера. Кроме того, у любого современного компьютера по каким-то неясным причинам могут быть слабые места, как в программном обеспечении, так и в электронном оборудовании, которые могут приводить к сбоям так редко, что остаются необнаруженными на протяжении нескольких лет, и поэтому в работе каждого компьютера могут быть незамеченные ошибки».

Те математики, которые отказались признать работу Хейкена и Аппеля, не могли отрицать, что все математики соглашались принимать традиционные доказательства, даже если они сами не проверяли их. В случае доказательства Великой теоремы Ферма, представленного Уайлсом, менее 10% специалистов по теории чисел полностью понимали его рассуждения, но все 100% сочли, что доказательство правильное. Те, кто не смог до конца понять все тонкости доказательства, приняли его потому, что доказательство признали другие — те, кто все понял, шаг за шагом проследил весь ход доказательства и проверил каждую деталь.

Еще более ярким примером может служить так называемое доказательство классификации конечных простых групп, состоящее из 500 отдельных работ, написанных более чем сотней ма-тематиков. Говорят, что полностью разобрался в этом доказательстве (общим объемом в 15000 страниц) один-единственный человек на свете — скончавшийся в 1992 году Дэниэл Горенстейн. Тем не менее, математическое сообщество в целом могло быть спокойным: каждый фрагмент доказательства был изучен группой специалистов, и каждая строка из 15000 страниц была десятки раз проверена и перепроверена. Что же касается проблемы четырех красок, то с ней дело обстояло иначе: она никем не была и не будет полностью проверена [14].

За несколько лет, прошедших с тех пор, как Хакен и Аппель сообщили о доказательстве теоремы о четырех красках, компьютеры неоднократно использовались для решения других, менее известных, но столь, же важных проблем. В математике — области, не ведавшей ранее вмешательства столь современной технологии, как компьютеры, — все больше и больше специалистов неохотно осваивали использование кремниевой логики и разделяли мнение У. Хейкена: «Всякий, в любом месте доказательства, может полностью вникнуть в детали и проверить их. То, что компьютер может за несколько часов «просмотреть» столько деталей, сколько человек не сможет просмотреть за всю свою жизнь, не меняет в принципе представление о математическом доказательстве. Меняется не теория, а практика математического доказательства».

Лишь совсем недавно математики наделили компьютеры еще большей властью, используя так называемые генетические алгоритмы. Это компьютерные программы, общая структура которых составлена математиком, но тонкие детали определяются самим компьютером. Некоторые направления, или «линии», в программе обладают способностью мутировать и эволюционировать наподобие индивидуальных генов в органической ДНК. Отправляясь от исходной материнской программы, компьютер может порождать сотни дочерних программ, слегка отличающихся из-за введенных компьютером случайных мутаций. Дочерние программы используются в попытках решения проблемы. Большинство программ бесславно не срабатывают, а та, которой удается дальше других продвинуться к желанному результату, используется в качестве материнской программы, порождающей новые поколение дочерних программ. Выживание наиболее приспособленного интерпретируется как выделение той из дочерних программ, которая позволяет особенно близко подойти к решению проблемы. Математики надеются, что, повторяя этот процесс, программа без вмешательства извне приблизится к решению проблемы. В некоторых случаях такой подход оказался весьма успешным.

Специалист в области «computer science» Эдвард Френкин даже заявил, что когда-нибудь компьютер найдет решение какой-нибудь важной проблемы без помощи математиков. Десять лет назад Френкин учредил премию Лейбница размером в 100000 долларов. Премия будет присуждена первой компьютерной программе, способной сформулировать и доказать теорему, которая окажет «глубокое влияние на развитие математики». Будет ли когда-нибудь присуждена премия Лейбница — вопрос спорный, но одно можно сказать со всей определенностью: компьютерной программе всегда будет недоставать прозрачности традиционных доказательств, и в сравнении с ними она будет проигрывать, уступая им в глубине. Математическое доказательство должно не только давать ответ на поставленный вопрос, но и способствовать пониманию, почему ответ именно таков, каков он есть, и в чем именно состоит его суть. Задавая вопрос на входе в черный ящик и получая ответ на выходе из него, мы увеличиваем знание, но не понимание.

Математик Рональд Грэхем описывает недостаточную глубину компьютерных доказательств на примере одной из великих не доказанных по сей день гипотез — гипотезы Римана: «Я был бы весьма и весьма разочарован, если бы можно было подключиться к компьютеру, спросить у него, верна ли гипотеза Римана, и получить в ответ: "Да, верна, но Вы не сможете понять доказательство"». Математик Филип Дэвис, похожим образом отреагировал на решение проблемы четырех красок: «Моей первой реакцией было: "Потрясающе! Как им удалось решить эту проблему?". Я ожидал какой-то блестящей новой идеи, красота которой перевернула бы всю мою жизнь. Но когда я услышал в ответ: "Они решили проблему, перебрав тысячи случаев и пропустив все варианты один за другим через компьютер", — меня охватило глубочайшее уныние. Я подумал: "Значит, все сводилось к простому перебору, и проблема четырех красок вовсе не заслуживала названия хорошей проблемы"» [13].

Говоря прямо, компьютерную часть доказательства К. Аппеля и У. Хейкена (Appel, Haken) невозможно проверить вручную, а традиционная часть доказательства длинна и сложна настолько, что ее никто целиком и не проверял. Между тем доказательство, не поддающееся проверке, есть нонсенс. Согласиться с подобным доказательством означает примерно то же, что просто поверить авторам [14]. Ведь даже проверка распечаток всех программ и всех входных данных не может гарантировать от случайных сбоев или даже от скрытых пороков электроники (вспомним, что ошибки при выполнении деления у первой версии процессора Pentium были случайно обнаружены спустя полгода после начала его коммерческих продаж (кстати, математиком, специалистом в теории чисел).

Все это казалось говорит о том, что недавно появились зловещие признаки того, что доказательство Уайлса Великой теоремы Ферма, возможно, стало одним из последних примеров классического доказательства, и будущие доказательства столь сложных проблем будут полагаться не столько на изящные рассуждения, сколько на грубую силу.

Рассмотрим вопрос построение алгоритма на основе работ [18,19]. В работах показано, что можно построив полиномиальный алгоритм раскраски плоского графа четырьмя цветами со сложностью $O(n^2)$. Однако простота этого алгоритма обманчива. Константа, стоящая при $n^2$, имеет порядок $10^5$— $10^6$, а требуемая память порядка $10^6$. Таким образом, несмотря на свою полиномиальность, алгоритм четырех раскрасок, построенный на основе доказательства теоремы о четырех красках [17,18], может оказаться практически не очень эффективным даже для сравнительно небольших n (n - количество вершин) [1,2].

В октябре 2002 г. в Великобритании отмечалось 150-летие проблемы четырех красок и 25-летие ее решения, данного К. Аппелем (K. Appel), У. Хейкеном (W. Haken) и Дж. Кохом (J. Koch). В октябре того же года эти события отмечались в течение специальной праздничной недели [13].

Проблема четырех красок кажется на первый взгляд изолированной задачей, мало связанной с другими разделами математики и практическими задачами. На самом деле это не так. Известно более 20 ее переформулировок, которые связывают эту проблему с задачами алгебры, статистической механики и задачами планирования. И это тоже характерно для математики: решение задачи, изучаемой из чистого любопытства, оказывается полезным в описании реальных и совершенно различных по своей природе объектов и процессов.

Методы доказательства приведенные К. Аппелем и У. Хейкеном можно отнести к комбинаторному подтверждению гипотезы четырех красок. Несмотря на опубликованное, до сих пор отсутствует чёткое описание механизма раскраски плоского графа четырьмя красками, его естественной сути и его связи с явлением планарности графа. Ведь мало сказать, что любой плоский граф раскрашиваем, нужно и показать, как его раскрасить. При таком подходе становится непонятным и возникновение уравнений Хивуда, и их роль в раскраске плоских графов. С точки зрения прикладной математики нужно получить такое доказательство, которое описывает механизм раскраски плоских графов и на основе которого можно строить простые и эффективные алгоритмы раскраски плоских графов.

Поэтому в данной работе рассмотрен другой подход, основанный на возможности сведения максимально плоского графа к регулярному плоскому кубическому графу с дальнейшей его раскраской. На основании теоремы Тейта-Волынского [1,5,6] вершины максимально плоского графа можно раскрасить четырьмя цветами, если ребра его двойственного кубического графа можно раскрасить тремя цветами. Рассматривая свойства раскрашенного кубического графа, можно показать, что сложение цветов подчиняется законам преобразования группы Клейна четвёртого порядка.

В работе рассмотрен рекурсивный процесс построения последующего кубического графа путём введения нового ребра в предыдущий плоский кубический граф. При таком рекурсивном рассмотрении можно показать, что теорема о четырёх красках есть следствие теоремы о существовании цветного диска проходящего по двум сцепленным рёбрам в плоском кубическом графе без мостов. Приведено дедуктивное доказательство данной теоремы. Следствием этой теоремы является раскрашиваемость вновь введённого ребра в цвет цветного диска для сцепленных ребер. Полученные следствия из данной теоремы позволяют утверждать о существовании правильной раскраски произвольного плоского кубического графа без мостов, а также устанавливают, что согласно теореме Петерсена существует цветной 2-фактор с цветными дисками четной длины.

Для поиска цветного диска проходящего по сцепленным рёбрам введена операция ротации цветных дисков. Приведены примеры для поиска такого цветного диска.

Рассмотрен метод раскраски кубического графа способом рекурсивного построения плоского кубического графа без мостов. Представлен способ раскраски методом удаления рёбер в плоском кубическом графе с чётной длиной циклов.

Рассмотрен переход от плоского кубического графа H к плоскому биквадратному графу $B_2$. На основании свойства эйлеровых циклов в биквадратном графе доказана теорема о существовании эйлерова маршрута с дисками, длина которых кратна трём. Полученный таким образом эйлеров цикл устанавливает раскраску рёбер плоского биквадратного графа в три цвета. Показано, что существование такого эйлерова маршрута порождает цветное вращение базовых треугольных циклов.

Цветное вращение базовых треугольных граней можно описать уравнениями Хивуда. Рассмотрены методы решения уравнений Хивуда с помощью цветного вращения треугольных граней в максимально плоском графе. Установлено соответствие между раскраской максимально плоского графа $G'$, плоского кубического графа H и плоского биквадратного графа $B_2$.

Рекурсивный подход к построению плоских кубических графов и введённая операция ротации цветных дисков позволяют строить простые и эффективные алгоритмы для раскраски произвольных максимально плоских графов.

Множество вершин и рёбер графа нумеруются на рисунках числами, с целью наглядности изображения. При этом предполагается, что вершины принадлежат множеству $v_j \in V, j = 1,2,...,n$. В свою очередь, рёбра принадлежат множеству $e_i \in E, i = 1,2,...,m$.

В основном работа написана для прикладных математиков. Представлен алгоритмический подход к раскраске плоских графов. В работе приведено большое количество примеров для пояснения рассматриваемых явлений и механизма раскраски. Данная работа написана с целью дальнейшего развития математических методов раскраски плоских графов. Описывается алгоритмический подход к раскраске плоских графов для построения вычислительного процесса с полиномиальной вычислительной сложностью $o(m^2)$. Параллельно рассматриваются методы описания и построения плоских рисунков графа, так как современные методы описывают графы с точностью до изоморфизма, игнорируя понятие рисунка плоского графа и вывода его на экран компьютера.

# Глава 1. Кубические графы и группа Клейна четвертого порядка

## 1.1. Максимально плоский граф

Рассмотрим произвольный связный плоский граф G, у которого степень каждой вершины больше двух (см. рис. 1.1). Рассмотрим в графе G каждую нетреугольную грань и, путем добавления ребер, разобьем ее на несколько треугольных граней (см. рис. 1.2). В результате получим *максимально плоский граф* G′, все грани которого будут треугольными. И хотя данный процесс не однозначен, однако если граф G′ раскрашен, то удаление вновь введенных ребер из графа **G′** не нарушит раскраску графа G′[7,17].

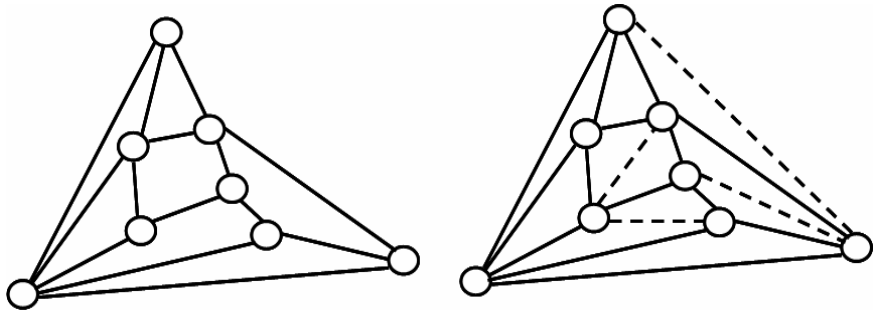

Рис. 1.1. Плоский граф G**.**     Рис. 1.2. Максимально плоский граф G′**.**

Выберем в максимально плоском графе G′ произвольную грань и обозначим её буквой *s* . Поставим в соответствие грани *s* вершину $x^0$ какого-то нового графа. Тогда каждому произвольному ребру *e* графа G′ поставим в соответствие ребро $e^0$, соединяющее те вершины $x^0_i$ и $x^0_j$, которые соответствуют граням $s_i$ и $s_j$ по обе стороны ребра *e* максимально плоского графа. Полученный таким образом граф $G^0$ является плоским, связным и называется двойственным к G′ или однородным графом H степени три, или кубическим графом H.

## 1.2. Плоский кубический граф

Построение дуального плоского кубического графа H представлено на рис. 1.3, где вершины закрашены более темным цветом, а ребра пунктирными линиями.

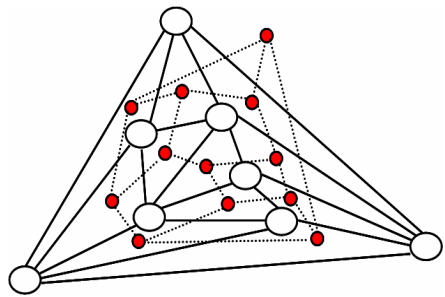

Рис. 1.3. Построение плоского кубического графа H**.**

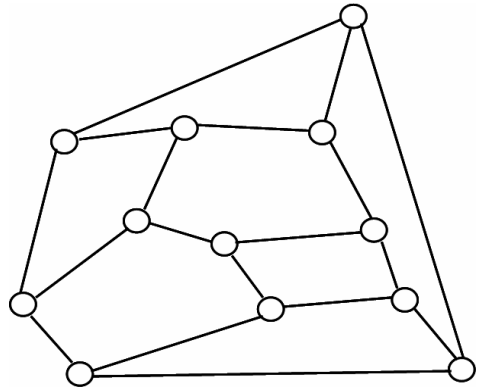

Рис. 1.4. Плоский кубический граф H.

Вершины в плоском кубическом графе соответствуют граням максимально плоского графа, а грани плоского кубического графа - соответствуют вершинам максимально плоского графа. Количество ребер в максимально плоском графе равно количеству ребер двойственного плоского кубического графа. Количество ребер в кубическом графе H с n вершинами определяется по формуле:

$$m = 3n/2. \qquad (1.1)$$

И, следовательно, количество ребер в таком графе всегда кратно трем:

$$m \equiv 0 \pmod 3. \qquad (1.2)$$

Так как в кубическом графе количество ребер *m* - целое число, следовательно, количество вершин *n* в однородном кубическом графе должно быть четным числом:

$$n \equiv 0 \pmod 2. \qquad (1.3)$$

### 1.3. Раскраска кубического графа

С точки зрения раскраски, между максимальным плоским графом G′ и двойственным к нему плоским кубическим графом H существует связь, устанавливаемая следующей теоремой.

**Теорема 1.1.** (Тэйт) [5]. *Пусть H — плоский однородный кубический граф без перешейков; необходимое и достаточное условие возможности такого раскрашивания граней графа четырьмя цветами, при котором никакие две смежные грани не окрашиваются в одинаковый цвет, состоит в том, чтобы хроматический класс графа H был равен трем.*

Теорема 1.1 устанавливает связь между раскраской вершин максимально плоского графа и раскраской ребер двойственному к нему кубического графа H. Однако, она не утверждает, что любой кубический граф без мостов имеет хроматический класс равный трем. Тэйт высказал предположение, что всякий плоский кубический граф без перешейков имеет гамильтонов цикл. Если бы ему это удалось доказать, то, тем самым, была бы решена проблема четырех красок, хотя сам Тэйт и не обнаружил этой зависимости. Однако Татт [1,5,6] показал, что это предпо-ложение неверно и привел пример кубического графа с 46 вершинами, не являющегося гамильтоновым (рис. 1.5).

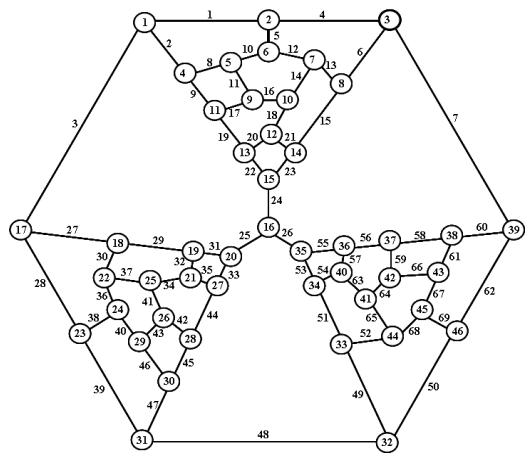

Рис. 1.5. Граф Татта.

Отсутствие в некоторых плоских кубических графах гамильтонова цикла (см. рис. 1.5 и рис. 1.6) породило сомнение у многих математиков относительно возможности построения механизма раскраски вершин плоского графа, используя раскраску ребер кубического графа. И поэтому, дальнейшее развитие методов раскраски плоских графов пошло по пути исследования достаточно большого так называемого неустранимого множества конфигураций [1,18,19].

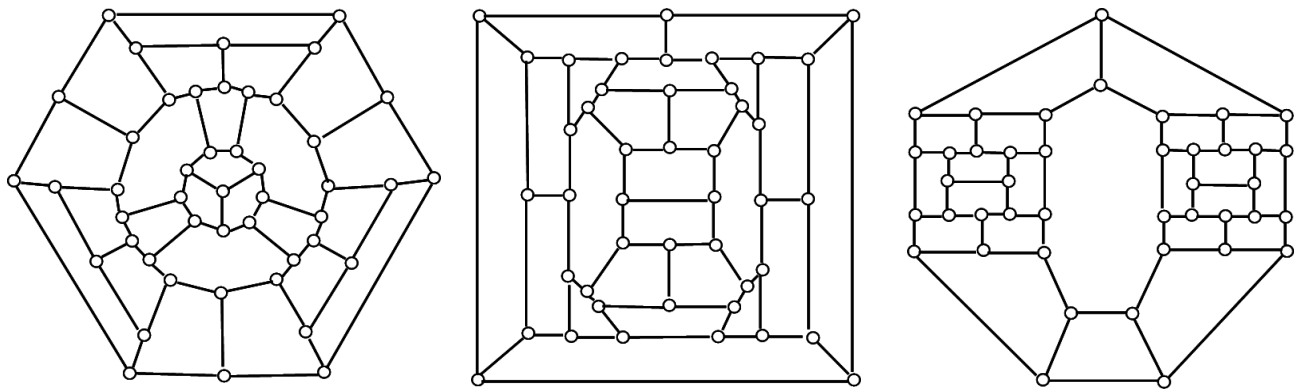

Рис. 1.6. Негамильтоновы плоские кубические графы.

Рассмотрение механизма такой раскраски ребер плоского кубического графа, в общем случае, начнем с описания группы симметрий тетраэдра.

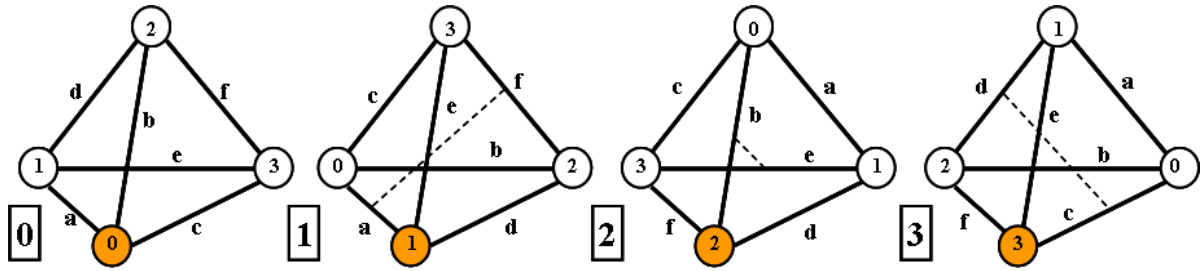

Рис. 1.7. Различные самосовмещения тетраэдра.

Известно, что группа симметрий тетраэдра содержит подгруппу $C_4$, называемую группой Клейна четвертого порядка [3,5] (см. рис.1.7). Если обозначить положения симметрии тетраэдра через 0,1,2,3, то таблица преобразования для этой группы имеет вид:

| + | 0 | 1 | 2 | 3 |
|---|---|---|---|---|
| 0 | 0 | 1 | 2 | 3 |
| 1 | 1 | 0 | 3 | 2 |
| 2 | 2 | 3 | 0 | 1 |
| 3 | 3 | 2 | 1 | 0 |

(1.4)

Будем раскрашивать ребра тетраэдра так, чтобы примкнувшие к каждой вершине ребра были окрашены тремя разными цветами.

Здесь, красный цвет будем обозначать буквой R или цифрой 1 для преобразования Клейна, а на рисунках графа ребра данного цвета будем представлять сплошной линией.

Синий цвет будем обозначать буквой B или цифрой 2, а на рисунках графа ребра данного цвета будем представлять точечной линией.

Зеленый цвет будем обозначать буквой G или цифрой 3, а на рисунках графа ребра данного цвета будем представлять пунктирной линией.

Белый цвет будем обозначать буквой W или цифрой 0, а на рисунках графа ребра данного цвета будем представлять штрихпунктирной линией.

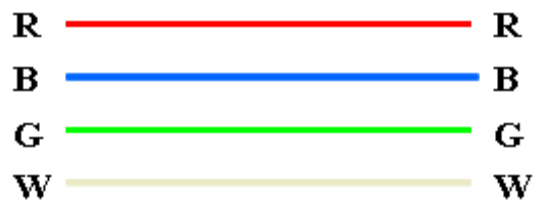

Рис. 1.8. Обозначение цветных линий.

Для удобства представления, вершины и ребра графа иногда будем обозначать цифрами или буквами, каждый раз оговаривая их применение или употребляя по смыслу.

Раскрасив ребра тетраэдра, перейдем к раскрашенным плоским кубическим графам. Самым минимальным планарным плоским кубическим графом К$_4$ без кратных ребер (см. рис. 1.10) может являться отображение тетраэдра на плоскость (см. рис. 1.9).

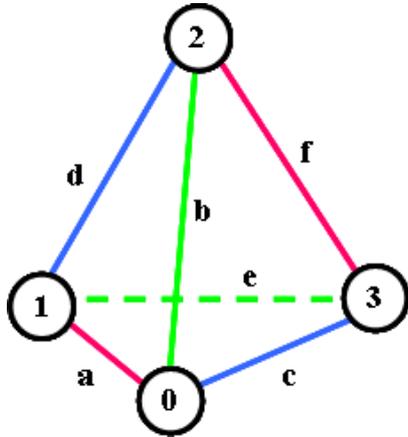 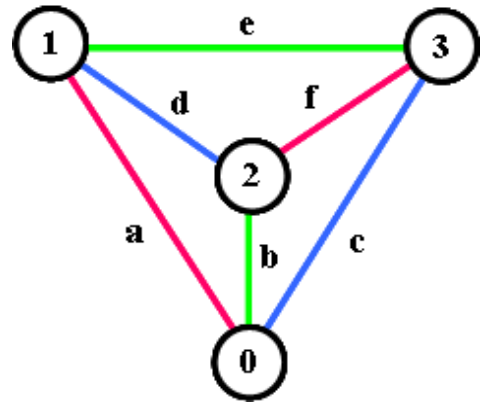

Рис. 1.9. Тетраэдр с раскрашенными ребрами.    Рис. 1.10. Плоский кубический граф К$_4$

Если кубический граф Н имеет хроматический класс равный трем, то его ребра должны быть раскрашены в три цвета. Раскраска предполагает наличие трех разноцветных ребер для каждой вершины графа Н. Кубический граф, имеющий хроматический класс равный трем будем называть раскрашенным плоским кубическим графом Н. Например, кубический граф, представленный на рис.1.11, является раскрашенным кубическим графом Н.

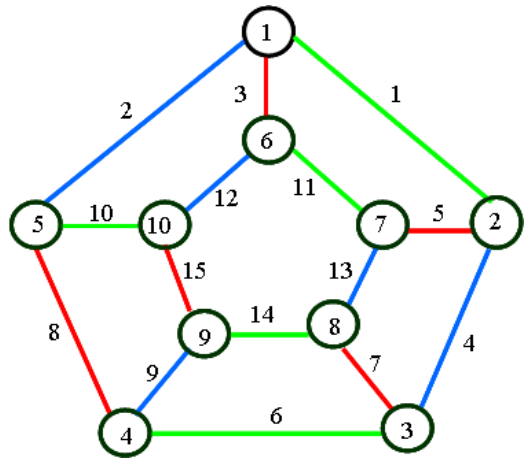

Рис. 1.11. Раскраска ребер в кубическом графе Н.

Раскраска ребер, в свою очередь, индуцирует раскраску граней плоского кубического графа Н. На рис.1.12 представлены 4 способа раскраски граней плоского кубического графа для раскрашенного кубического графа на рис 1.11.

*Квазициклом* будем называть цикл с четной валентностью вершин[5,6]. Гамильтоновым квазициклом будем называть квазицикл с валентностью вершин равной двум и проходящий по всем вершинам графа.

Пусть G = (X,E;P) - граф с пронумерованным множеством рёбер E = {$e_1, e_2,...,e_m$}, а $L_G$ - множество всех суграфов этого графа. Относительно операции сложения (будем называть её кольцевой суммой):

$$(X,E_1;P) \oplus (X,E_2;P) = (X,(E_1 \cup E_2)\setminus(E_1\cap E_2);P) \qquad (1.5)$$

множество $L_G$ образует абелеву 2-группу [15].

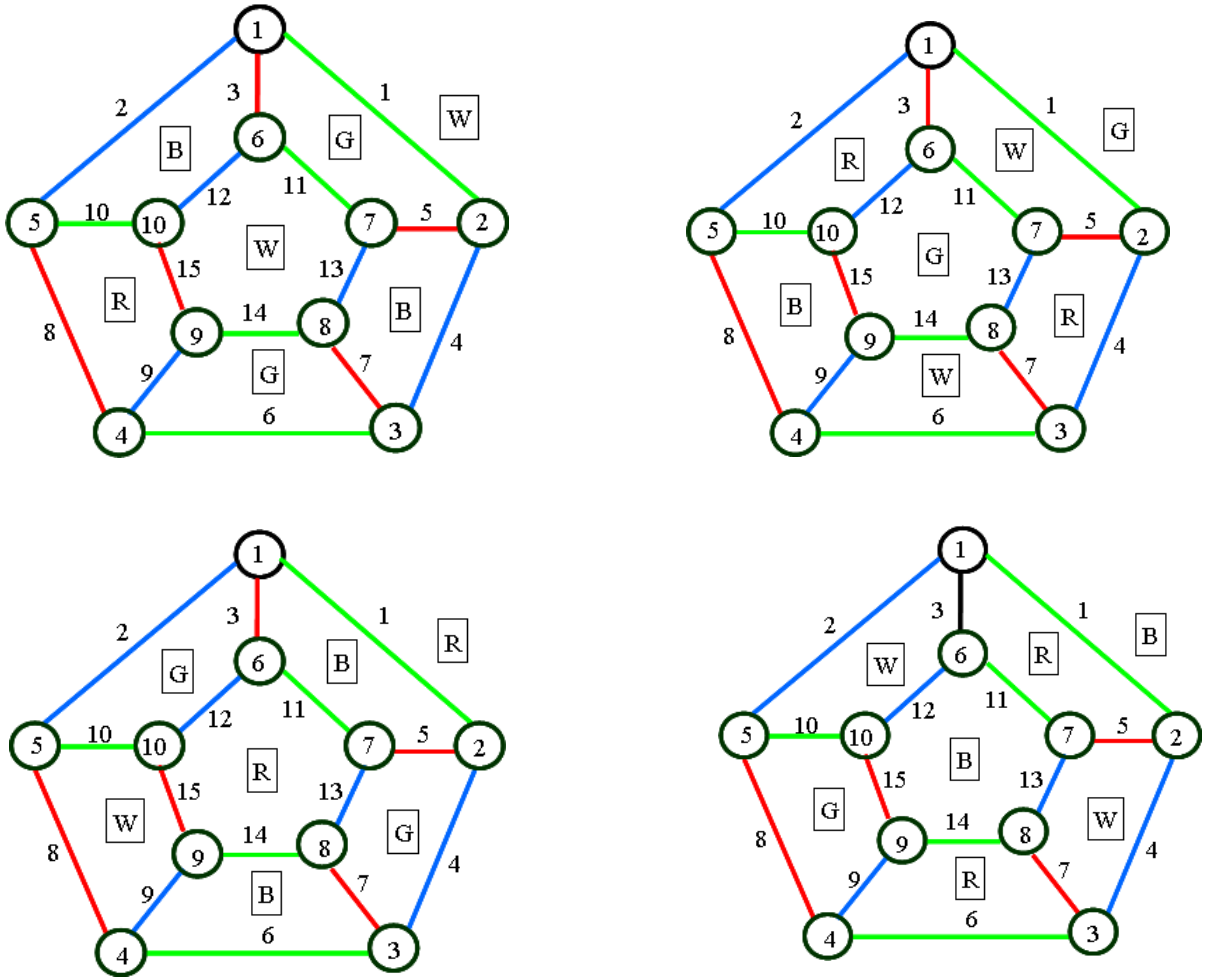

Рис. 1.12. Раскраска граней плоского раскрашенного кубического графа.

Действительно, $L_G$ заведомо является группоидом. Относя каждому суграфу G = (X,E;P) строку чисел ($\alpha_1,\alpha_2,...,\alpha_i,...,\alpha_m$), в которой i = (1,2,...,m), и определяя сложение строк как покомпонентное по модулю 2, мы получим изоморфный $L_G$ группоид, элементами которого служат всевозможные строки длины *m* из нулей и единиц и который представляет собой абелеву 2-группу.

В дальнейшем группу $L_G$ удобно рассматривать как линейное пространство над полем коэффициентов GF(2) = {0,1}, называемое пространством суграфов данного графа G. Размерность этого пространства dim $L_G$ = *m*, ибо множество элементов (1,0,....0), (0.1.....0),....(0.0....,1), представляющих однорёберные суграфы, образует базис.

В дальнейшем нам понадобится понятие изометрического цикла в графе.

Расстоянием $\rho_G(x,y)$ в графе G между вершинами x и y графа G = (X,E;P) называется длина кратчайшего из маршрутов (и, значит, кратчайшей из простых цепей), соединяющих эти вершины. Если x и y отделены в G, то $\rho_G(x,y) = +\infty$. Функция $\rho = \rho(x,y) = \rho_G(x,y)$, определенная на множестве всех пар вершин графа G и принимающая целые неотрицательные значения (к числу которых мы относим и бесконечное), заслуживает названия метрики графа, поскольку она удовлетворяет трем аксиомам Фреше:

$$\forall x, y \in X [\rho(x, y) = 0 \Leftrightarrow x = y], \qquad (1.6)$$

$$\forall x, y \in X [\rho(x, y) = \rho(y, x_1)], \qquad (1.7)$$

$$\forall x, y, z \in X [\rho(x, y) + \rho(y, z) \geq \rho(x, z)]. \qquad (1.8)$$

Введем следующее понятия, связанные с метрикой графа

**Определение 1.1**[20]. *Изометрический подграф* – подграф G′ графа G, у которого все расстояния внутри G′ те же самые, что и в G.

**Определение 1.2**[20]. *Изометрическим циклом* в графе называется простой цикл, для которого кратчайший путь между любыми двумя его несмежными вершинами состоит из ребер этого цикла. Изометрический цикл представляет собой частный случай изометрического подграфа.

Будем также считать, что для плоских графов граница грани, кроме обода определяется изометрическим циклом.

С другой стороны, каждый гамильтонов квазицикл есть 2-фактор графа, так как n – фактор – это регулярный суграф степени n. Тогда имеет место следующая теорема для кубических графов.

**Теорема 1.2**. (Петерсен) [17]. *Любой кубический граф, не содержащий мостов, можно представить в виде суммы 1-фактора и 2-фактора.*

Характерный непланарный кубический граф Петерсена представлен на рис. 1.13.

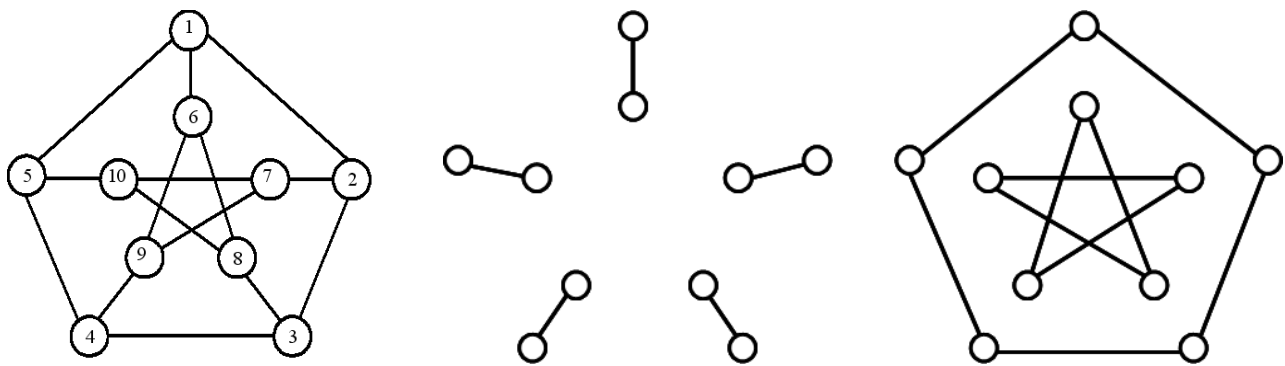

б) 1-фактор графа Петерсена    в) 2-фактор графа Петерсена

Рис. 1.13. Кубический граф Петерсена.

В свою очередь, гамильтонов квазицикл (2-фактор) может состоять из нескольких простых циклов. Каждый такой простой цикл четной длины будем называть – *цветным диском*, так как его ребра можно последовательно раскрасить двумя цветами. На рис. 1.14,г представлены два диска одного цвета, один проходит по ребрам $\{e_2,e_3,e_8,e_9,e_{12},e_{15}\}$, а другой по ребрам $\{e_4,e_5,e_7,e_{13}\}$. Гамильтонов квазицикл состоящий из цветных дисков, будем называть цветным гамильтоновым квазициклом или цветным 2-фактором.

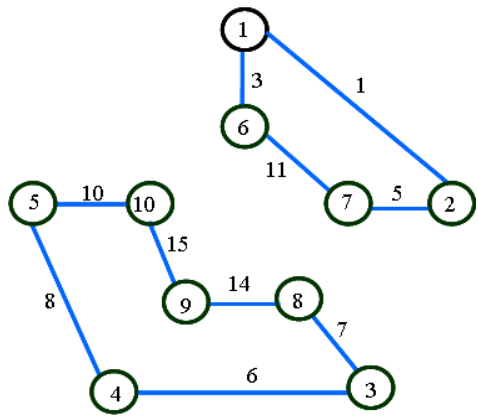

а) синий гамильтонов квазицикл (2-фактор)

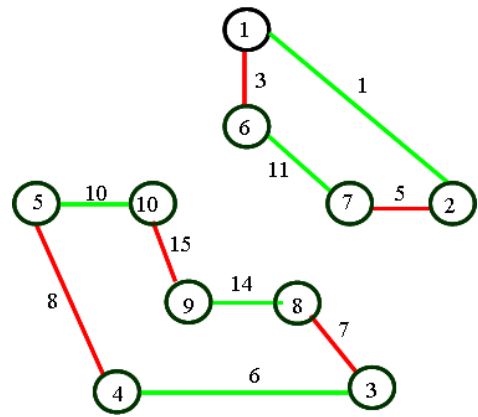

б) тот же квазицикл состоящий из красных и зеленых ребер

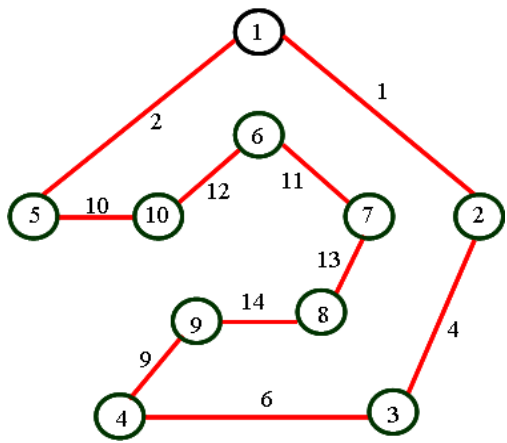

в) красный гамильтонов квазицикл

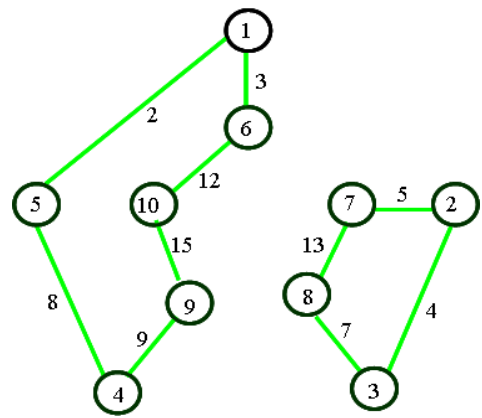

г) зеленый гамильтонов квазицикл

Рис. 1.14. Цветные индуцированные гамильтоновы квазициклы (2-факторы).

В раскрашенном кубическом графе обязательно существует цветная конфигурация, состоящая из цветных 2-факторов с дисками только четной длины и цветных 1-факторов.

Любой цветной гамильтонов квазицикл состоит из дисков. Если цветной диск один, то такой гамильтонов квазицикл - гамильтонов цикл.

### 1.4. Свойства раскрашенных плоских кубических графов

Цветные гамильтоновы квазициклы (цветные 2-факторы) обладают следующими свойствами:

**Свойство 1.** *В любом цветном диске число вершин четно.*

Действительно, для того чтобы диск был цветным, необходимо чтобы его ребра были раскрашены в два цвета. Это, в свою очередь, возможно только в случае, если количество вершин в цикле четно.

**Свойство 2.** *Любой цветной гамильтонов квазицикл (цветной 2-фактор) порождает два других цветных гамильтоновых квазицикла (цветных 2-фактора) и производит раскраску графа H.*

Если задан цветной гамильтонов квазицикл, то ребра графа этого 2-фактора можно раскрасить в два цвета отличных от цвета 2-фактора, а нераскрашенные ребра графа - в цвет самого гамильтонова квазицикла (цветного 2-фактора). Таким образом, кубический граф будет раскрашен, а раскраска ребер, в свою очередь, порождает три цветных гамильтоновых квазицикла (цветных 2-фактора).

**Свойство 3.** *Сумма по модулю 2 трех цветных гамильтоновых квазициклов (цветных 2-факторов), порожденных раскраской кубического графа H, есть пустое множество.*

Так как по каждому ребру проходит два цветных гамильтоновых квазицикла и каждое ребро в сумме появляется дважды, а это, в свою очередь, при сложении по модулю 2 тождественно $\varnothing$. Поэтому сумма цветных гамильтоновых квазициклов, порожденных правильной раскраской кубического графа H, есть пустое множество.

На базисе подпространства циклов графа G можно определить функционал Маклейна принимающий нулевое значение, если граф планарен [10]:

$$F(C) = \sum_{i=1}^{m}(a_i - 1)(a_i - 2) = \sum_{i=1}^{m}a_i^2 - 3\sum_{i=1}^{m}a_i + 2m \ \square . \tag{1.9}$$

Здесь $a_i$ – количество изометрических циклов проходящих по ребру *i*.

Кроме того, в раскрашенном кубическом графе любой цветной 2-фактор образуется как конечная сумма изометрических циклов:

$R_c$ - сумма изометрических циклов образующих красный 2-фактор;

$B_c$ - сумма изометрических циклов образующих синий 2-фактор;

$G_c$ - сумма изометрических циклов образующих зеленый 2-фактор;

$W_c$ - сумма изометрических циклов не вошедших ни в красный, ни в синий и не в зеленый 2-факторы.

С учетом сказанного и введенной операции кольцевого сложения циклов, цветные гамильтоновы квазициклы и цветные 1-факторы обладают следующими свойствами:

$R_c \oplus G_c \oplus B_c = \varnothing;$ (1.10)
$R_f \oplus G_f \oplus B_f = H;$ (1.11)

$$R_f \oplus R_c = H; \quad (1.12)$$
$$G_f \oplus G_c = H; \quad (1.13)$$
$$B_f \oplus B_c = H. \quad (1.14)$$

Здесь индекс *c* обозначает множество ребер, принадлежащих цветному 2-фактору (гамильтонову квазициклу), индекс *f* обозначает множество ребер, принадлежащих 1-фактору. H - множество ребер исходного плоского кубического графа.

Проявление свойства двойственности построения цветных 2-факторов описывается следующими выражениями:

$$R_c = \overline{R}_c, \quad B_c = \overline{B}_c, \quad G_c = \overline{G}_c, \quad (1.15)$$

где: $\overline{R}_c$, $\overline{B}_c$, $\overline{G}_c$ являются соответственными дополнениями изометрических циклов $R_c$, $B_c$, $G_c$ до полного базиса плоского графа и обода графа H. То есть:

$$R_c \oplus \overline{R}_c = B_c \oplus \overline{B}_c = G_c \oplus \overline{G}_c = \varnothing \quad (1.16)$$

**Пример 1.1.** Сказанное поясним на примере плоского кубического графа H, представленного на рис. 1.11.

Для данного кубического графа H базисные изометрические циклы (границы граней плоского графа) и обод:

$c_1 = \{e_1,e_3,e_5,e_{11}\}$, $c_2 = \{e_4,e_5,e_7,e_{13}\}$, $c_3 = \{e_6,e_7,e_9,e_{14}\}$, $c_4 = \{e_8,e_9,e_{10},e_{15}\}$, $c_5 = \{e_2,e_3,e_{10},e_{12}\}$,
$c_6 = \{e_{11},e_{12},e_{13},e_{14},e_{15}\}$, $c_0 = \{e_1,e_2,e_4,e_6,e_8\}$.

Для раскраски, представленной на рис.1.11, имеем:

$R_c = c_1 \oplus c_2 \oplus c_3 \oplus c_5 = \{e_1,e_3,e_5,e_{11}\} \oplus \{e_4,e_5,e_7,e_{13}\} \oplus \{e_6,e_7,e_9,e_{14}\} \oplus \{e_2,e_3,e_{10},e_{12}\} =$
$= \{e_1,e_2,e_4,e_6,e_9,e_{10},e_{11},e_{12},e_{13},e_{14}\}$;
$\overline{R}_c = c_4 \oplus c_6 \oplus c_0 = \{e_8,e_9,e_{10},e_{15}\} \oplus \{e_{11},e_{12},e_{13},e_{14},e_{15}\} \oplus \{e_1,e_2,e_4,e_6,e_8\} =$
$= \{e_1,e_2,e_4,e_6,e_9,e_{10},e_{11},e_{12},e_{13},e_{14}\}$;
$B_c = c_1 \oplus c_3 \oplus c_4 = \{e_1,e_3,e_5,e_{11}\} \oplus \{e_6,e_7,e_9,e_{14}\} \oplus \{e_8,e_9,e_{10},e_{15}\} =$
$= \{e_1,e_3,e_5,e_6,e_7,e_8,e_{10},e_{11},e_{14},e_{15}\}$;
$\overline{B}_c = c_2 \oplus c_5 \oplus c_6 \oplus c_0 = \{e_4,e_5,e_7,e_{13}\} \oplus \{e_2,e_3,e_{10},e_{12}\} \oplus \{e_{11},e_{12},e_{13},e_{14},e_{15}\} \oplus$
$\oplus \{e_1,e_2,e_4,e_6,e_8\} = \{e_1,e_3,e_5,e_6,e_7,e_8,e_{10},e_{11},e_{14},e_{15}\}$;
$G_c = c_2 \oplus c_4 \oplus c_5 = \{e_4,e_5,e_7,e_{13}\} \oplus \{e_8,e_9,e_{10},e_{15}\} \oplus \{e_2,e_3,e_{10},e_{12}\} =$
$= \{e_2,e_3,e_4,e_5,e_7,e_8,e_9,e_{12},e_{13},e_{15}\}$;
$\overline{G}_c = c_1 \oplus c_3 \oplus c_6 \oplus c_0 = \{e_1,e_3,e_5,e_{11}\} \oplus \{e_6,e_7,e_9,e_{14}\} \oplus \{e_{11},e_{12},e_{13},e_{14},e_{15}\} \oplus$
$\oplus \{e_1,e_2,e_4,e_6,e_8\} = \{e_2,e_3,e_4,e_5,e_7,e_8,e_9,e_{12},e_{13},e_{15}\}$;
$W_c = c_6 \oplus c_0 = \{e_{11},e_{12},e_{13},e_{14},e_{15}\} \oplus \{e_1,e_2,e_4,e_6,e_8\} = \{e_1,e_2,e_4,e_6,e_8,e_{11},e_{12},e_{13},e_{14},e_{15}\}$;
$R_c \oplus G_c \oplus B_c = \{e_1,e_2,e_4,e_6,e_9,e_{10},e_{11},e_{12},e_{13},e_{14}\} \oplus \{e_1,e_3,e_5,e_6,e_7,e_8,e_{10},e_{11},e_{14},e_{15}\} \oplus$
$\oplus \{e_2,e_3,e_4,e_5,e_7,e_8,e_9,e_{12},e_{13},e_{15}\} = \varnothing$.
$R_f = \{e_3,e_5,e_7,e_8,e_{15}\}$; $B_f = \{e_2,e_4,e_9,e_{12},e_{13}\}$; $G_f = \{e_1,e_6,e_{10},e_{11},e_{14}\}$;
$R_f \oplus G_f \oplus B_f = H$;
$R_f \oplus R_c = \{e_1,e_3,e_5,e_6,e_7,e_8,e_{10},e_{11},e_{13},e_{14}\} \oplus \{e_3,e_5,e_7,e_8,e_{15}\} = H$;
$B_f \oplus B_c = \{e_1,e_3,e_5,e_6,e_7,e_8,e_{10},e_{11},e_{14},e_{15}\} \oplus \{e_2,e_4,e_9,e_{12},e_{13}\} = H$;
$G_f \oplus G_c = \{e_2,e_3,e_4,e_5,e_7,e_8,e_9,e_{12},e_{13},e_{15}\} \oplus \{e_1,e_6,e_{10},e_{11},e_{14}\} = H$.

## 1.5. Раскраска граней плоского кубического графа H

Будем называть гранью - 2 клетку ограниченную изометрическим циклом или ободом. Три цветных 2-фактора, в свою очередь, порождают (индуцируют) раскраску граней плоского кубического графа H. В свою очередь, раскраска граней отождествляет раскраску изометрических циклов в цвет грани (так как границей грани служит изометрический цикл). Для раскрашенных плоских кубических графов кольцевая сумма базисных изометрических циклов и обода согласно теореме Маклейна [10] равна:

$$\sum_{i=1}^{m-n+1} c_i \oplus c_0 = \varnothing. \tag{1.17}$$

С другой стороны, кольцевая сумма трех цветных гамильтоновых квазициклов (цветных 2-факторов) согласно выражению (1.10):

$$R_c \oplus B_c \oplus G_c = \varnothing. \tag{1.18}$$

Следующее равенство связывает условие планарности Маклейна записанное в виде суммы изометрических циклов и обода с кольцевой суммой трех цветных 2-фактора для плоского кубического графа H:

$$\sum_{i=1}^{m-n+1} c_i \oplus c_0 = R_c \oplus G_c \oplus B_c = \overline{R}_c \oplus \overline{B}_c \oplus \overline{G}_c = \overline{R}_c \oplus B_c \oplus G_c = R_c \oplus \overline{B}_c \oplus \overline{G}_c =$$
$$= R_c \oplus \overline{B}_c \oplus G_c = \overline{R}_c \oplus B_c \oplus \overline{G}_c = R_c \oplus B_c \oplus \overline{G}_c = \overline{R}_c \oplus \overline{B}_c \oplus G_c = \varnothing. \tag{1.19}$$

Будем рассматривать раскрашенный плоский кубический граф. Обозначим кольцевую сумму изометрических циклов, описывающих границы всех граней красного цвета, через $\sum g_R$. Кольцевую сумму изометрических циклов описывающих границы граней синего цвета, обозначить как $\sum g_B$. Кольцевую сумму изометрических циклов описывающих границы граней зеленого цвета обозначим как $\sum g_G$. Кольцевую сумму всех изометрических циклов описывающих границы граней белого цвета как $\sum g_W$. Тогда можно записать следующие выражения:

$$G_c = \sum g_R \oplus \sum g_B = \sum g_G \oplus \sum g_W ; \tag{1.20}$$

$$B_c = \sum g_R \oplus \sum g_G = \sum g_B \oplus \sum g_W ; \tag{1.21}$$

$$R_c = \sum g_G \oplus \sum g_B = \sum g_R \oplus \sum g_W . \tag{1.22}$$

В свою очередь, три цветных 2-фактора индуцируют раскраску П граней, где границы граней суть базисные изометрические циклы и обод плоского графа H, и раскраска граней описывается соответствием Г = {П, С, К}. Здесь: П - множество пар кортежей <$c_i$, $a_{i,1} \times 1 + a_{i,2} \times 2 + a_{i,3}$

× 3 + $a_{i,4}$ × 0>, С – множество базисных изометрических циклов и обода плоского кубического графа Н, К – множество четырех цветов.

Множество П характеризует раскраску граней плоского кубического графа Н**.** Данное множество также представляет собой произведение элементов множества С и элементов множества К, а его элементы – это кортежи записанные в виде:

{<$c_1$, $a_{1,1}$ × 1 + $a_{1,2}$ × 2 + $a_{1,3}$ × 3 + $a_{1,4}$ × 0>,
<$c_2$, $a_{2,1}$ × 1 + $a_{2,2}$ × 2 + $a_{2,3}$ × 3 + $a_{2,4}$ × 0>,
……………………………………
<$c_{m-n+1}$, $a_{m-n+1,1}$ ×1 + $a_{m-n+1,2}$ × 2 + $a_{m-n+1,3}$ × 3 + $a_{m-n+1,4}$ × 0>,
<$c_0$, $a_{m-n+2,1}$ × 1 + $a_{m-n+2,2}$ × 2 + $a_{m-n+2,3}$ × 3 + $a_{m-n+2,4}$ × 0>}.  (1.23)

Здесь: $a_{i,j} \in \{0,1\}$; i = (1,2,…,m-n+2); j = (1,2,3,4). Коэффициенты $a_{i,j}$ при изометрических циклах определяют принадлежность изометрического цикла, описывающего границу грани соответствующему цветному 2-фактору. Сложение производится по законам операции преобразрвания группы Клейна, а умножение по обычным правилам арифметики.

Множество изометрических циклов и обода обозначим С = {$c_1$,$c_2$, ,$c_{m-n+1}$,$c_0$}, множество цветов можно записать как: К = {1,2,3,0}. В следующей записи пары трех цветных 2-факторов индуцируют одинаковую раскраску граней:

($R_c \oplus B_c \oplus G_c$) ∨ ($\overline{R}_c \oplus \overline{B}_c \oplus \overline{G}_c$) индуцирует раскраску $Г_1$ = {$П_1$, С, К},  (1.24)

($\overline{R}_c \oplus B_c \oplus G_c$) ∨ ($R_c \oplus \overline{B}_c \oplus \overline{G}_c$) индуцирует раскраску $Г_2$ = {$П_2$, С, К},  (1.25)

($R_c \oplus \overline{B}_c \oplus G_c$) ∨ ($\overline{R}_c \oplus B_c \oplus \overline{G}_c$) индуцирует раскраску $Г_3$ = {$П_3$, С, К},  (1.26)

($R_c \oplus B_c \oplus \overline{G}_c$) ∨ ($\overline{R}_c \oplus \overline{B}_c \oplus G_c$) индуцирует раскраску $Г_4$ = {$П_4$, С, К}.  (1.27)

Сказанное рассмотрим для раскраски графа представленной на рисунке 1.11.

**Пример 1.1.** Пусть раскраска граней, где границы граней суть изометрические циклы, индуцируется следующим набором трех цветных 2-факторов $R_c \oplus B_c \oplus G_c$ = ($c_1 \oplus c_2 \oplus c_3 \oplus c_5$) ⊕ ⊕ ($c_1 \oplus c_3 \oplus c_4$) ⊕ ($c_2 \oplus c_4 \oplus c_5$) = ∅. Здесь цикл $c_1$ принадлежит и $R_c$ и $B_c$ и т.д.

Тогда коэффициенты $a_{i,j}$ можно записать в виде следующей матрицы:

|       | 1(R) | 2(B) | 3(G) | 0(W) |
|-------|------|------|------|------|
| $c_1$ | 1    | 1    | 0    | 0    |
| $c_2$ | 1    | 0    | 1    | 0    |
| $c_3$ | 1    | 1    | 0    | 0    |
| $c_4$ | 0    | 1    | 1    | 0    |
| $c_5$ | 1    | 0    | 1    | 0    |
| $c_6$ | 0    | 0    | 0    | 0    |
| $c_0$ | 0    | 0    | 0    | 0    |

Раскраска граней $П_1$ индуцированное $R_c \oplus B_c \oplus G_c$ может быть записано в виде:

$$\Pi_1 = \{<c_1, 1 \times 1 + 1 \times 2 + 0 \times 3 + 0 \times 0>,$$
$$<c_2, 1 \times 1 + 0 \times 2 + 1 \times 3 + 0 \times 0>,$$
$$<c_3, 1 \times 1 + 1 \times 2 + 0 \times 3 + 0 \times 0>,$$
$$<c_4, 0 \times 1 + 1 \times 2 + 1 \times 3 + 0 \times 0>,$$
$$<c_5, 1 \times 1 + 0 \times 2 + 1 \times 3 + 0 \times 0>,$$
$$<c_6, 0 \times 1 + 0 \times 2 + 0 \times 3 + 0 \times 0>,$$
$$<c_0, 0 \times 1 + 0 \times 2 + 0 \times 3 + 0 \times 0>\}.$$

После проведения операций сложения и умножения получим:

$\Pi_1 = \{<c_1, 3>,<c_2, 2>,<c_3, 3>,<c_4, 1>,<c_5, 2>,<c_6, 0>, <c_0, 0>\}$, или в виде

$\Pi_1 = \{<c_1, G>,<c_2, B>,<c_3, G>,<c_4, R>,<c_5, B>,<c_6, W>, <c_0, W>\}$.

Пусть раскраска граней индуцируется следующими другими тремя цветными 2-факторами $\overline{R}_c \oplus \overline{B}_c \oplus \overline{G}_c = (c_4 \oplus c_6 \oplus c_0) \oplus (c_2 \oplus c_5 \oplus c_6 \oplus c_0) \oplus (c_1 \oplus c_3 \oplus c_6 \oplus c_0) = \varnothing$.

Тогда коэффициенты $a_{i,j}$ можно записать в виде следующей матрицы:

|        | 1(R) | 2(B) | 3(G) | 0(W) |
|--------|------|------|------|------|
| $c_1$  | 0    | 0    | 1    | 0    |
| $c_2$  | 0    | 1    | 0    | 0    |
| $c_3$  | 0    | 0    | 1    | 0    |
| $c_4$  | 1    | 0    | 0    | 0    |
| $c_5$  | 0    | 1    | 0    | 0    |
| $c_6$  | 1    | 1    | 1    | 0    |
| $c_0$  | 1    | 1    | 1    | 0    |

Раскраска граней $\Pi_1$ индуцированное $\overline{R}_c \oplus \overline{B}_c \oplus \overline{G}_c$ может быть записано в виде:

$$\Pi_1 = \{<c_1, 0 \times 1 + 0 \times 2 + 1 \times 3 + 0 \times 0>,$$
$$<c_2, 0 \times 1 + 1 \times 2 + 0 \times 3 + 0 \times 0>,$$
$$<c_3, 0 \times 1 + 0 \times 2 + 1 \times 3 + 0 \times 0>,$$
$$<c_4, 1 \times 1 + 0 \times 2 + 0 \times 3 + 0 \times 0>,$$
$$<c_5, 0 \times 1 + 1 \times 2 + 0 \times 3 + 0 \times 0>,$$
$$<c_6, 1 \times 1 + 1 \times 2 + 1 \times 3 + 0 \times 0>,$$
$$<c_0, 1 \times 1 + 1 \times 2 + 1 \times 3 + 0 \times 0>\}.$$

И хотя коэффициенты $a_{i,j}$ отличны от предыдущего случая, после проведения операций сложения и умножения получим то же множество цветов для граней плоского кубического графа:

$\Pi_1 = \{<c_1, 3>,<c_2, 2>,<c_3, 3>,<c_4, 1>,<c_5, 2>,<c_6, 0>, <c_0, 0>\}$, или в виде:

$\Pi_1 = \{<c_1, G>,<c_2, B>,<c_3, G>,<c_4, R>,<c_5, B>,<c_6, W>, <c_0, W>\}$.

Если раскраска граней индуцируется следующими тремя цветными 2-факторами $\overline{R}_c \oplus B_c \oplus \oplus G_c = (c_4 \oplus c_6 \oplus c_0) \oplus (c_1 \oplus c_3 \oplus c_4) \oplus (c_2 \oplus c_4 \oplus c_5) = \varnothing$.

Тогда коэффициенты $a_{i,j}$ можно записать в виде следующей матрицы:

| 1(R) | 2(B) | 3(G) | 0(W) |

| | | | | |
|---|---|---|---|---|
| $c_1$ | 0 | 1 | 0 | 0 |
| $c_2$ | 0 | 0 | 1 | 0 |
| $c_3$ | 0 | 1 | 0 | 0 |
| $c_4$ | 1 | 1 | 1 | 0 |
| $c_5$ | 0 | 0 | 1 | 0 |
| $c_6$ | 1 | 0 | 0 | 0 |
| $c_0$ | 1 | 0 | 0 | 0 |

Множество раскраски изометрических циклов $П_2$ индуцированное $\overline{R}_c \oplus B_c \oplus G_c$, может быть записано в виде:

$П_2 = \{<c_1, 0 \times 1 + 1 \times 2 + 0 \times 3 + 0 \times 0>,$
$<c_2, 0 \times 1 + 0 \times 2 + 1 \times 3 + 0 \times 0>,$
$<c_3, 0 \times 1 + 1 \times 2 + 0 \times 3 + 0 \times 0>,$
$<c_4, 1 \times 1 + 1 \times 2 + 1 \times 3 + 0 \times 0>,$
$<c_5, 0 \times 1 + 0 \times 2 + 1 \times 3 + 0 \times 0>,$
$<c_6, 1 \times 1 + 0 \times 2 + 0 \times 3 + 0 \times 0>,$
$<c_0, 1 \times 1 + 0 \times 2 + 0 \times 3 + 0 \times 0>\}.$

После проведения операций сложения и умножения получим то же множество:

$П_2 = \{<c_1, 2>, <c_2, 3>, <c_3, 2>, <c_4, 0>, <c_5, 3>, <c_6, 1>, <c_0, 1>\}$,

или в виде: $П_2 = \{<c_1, B>, <c_2, G>, <c_3, B>, <c_4, W>, <c_5, G>, <c_6, R>, <c_0, R>\}$.

Раскраска граней индуцированное тремя цветными 2-факторами $R_c \oplus \overline{B}_c \oplus \overline{G}_c$ равна множеству $П_2$.

Проделав аналогичные вычисления получим следующее множество раскрасок граней $П_3$ индуцированного тремя цветными 2-факторами $R_c \oplus \overline{B}_c \oplus G_c$ или $\overline{R}_c \oplus B_c \oplus \overline{G}_c$:

$П_3 = \{<c_1, R>, <c_2, W>, <c_3, R>, <c_4, G>, <c_5, W>, <c_6, B>, <c_0, B>\}$.

Для множества $П_4$, индуцированного тремя цветными 2-факторами $R_c \oplus B_c \oplus \overline{G}_c$ или $\overline{R}_c \oplus \overline{B}_c \oplus G_c$:

$П_4 = \{<c_1, W>, <c_2, R>, <c_3, W>, <c_4, B>, <c_5, R>, <c_6, G>, <c_0, G>\}$.

Переход от раскраски граней индуцированной тремя цветными 2-факторами $(R_c \oplus B_c \oplus G_c) \vee (\overline{R}_c \oplus \overline{B}_c \oplus \overline{G}_c)$ к другой раскраске граней индуцированной тремя цветными 2-факторами $(\overline{R}_c \oplus B_c \oplus G_c) \vee (R_c \oplus \overline{B}_c \oplus \overline{G}_c)$ осуществляется путем привнесения цвета R ко всем цветам граней плоского кубического графа по законам преобразования группы Клейна. Например, для нашего примера, добавление цвета R к цветам граней $П_1$ приведет к следующей перекраски цветов граней плоского кубического графа $П_2$:

$П_2 = П_1 + R = \{<c_1, G+R>, <c_2, B+R>, <c_3, G+R>, <c_4, R+R>, <c_5, B+R>, <c_6, W+R>,$
$<c_0, W+R>\} = \{<c_1, B>, <c_2, G>, <c_3, B>, <c_4, W>, <c_5, G>, <c_6, R>, <c_0, R>\}$.

На следующей диаграмме представлен процесс замены трех цветных 2-факторов индуцирующих раскраску граней плоского кубического графа. Нижние три цветные 2-факторы индуцируют раскраску граней плоского кубического графа путем добавления цвета R к исходной раскраске граней плоского кубического графа индуцируемой тремя цветными 2-факторами расположенными сверху:

$$(R_c \oplus B_c \oplus G_c) \vee (\overline{R}_c \oplus \overline{B}_c \oplus \overline{G}_c)$$
$$(\overline{R}_c \oplus B_c \oplus G_c) \vee (R_c \oplus \overline{B}_c \oplus \overline{G}_c)$$

Здесь $R_c$ и $\overline{R}_c$ меняются местами.

Аналогично строится диаграмма и при добавлении цвета B к исходной раскраски граней плоского кубического графа:

$$(R_c \oplus B_c \oplus G_c) \vee (\overline{R}_c \oplus \overline{B}_c \oplus \overline{G}_c)$$
$$(R_c \oplus \overline{B}_c \oplus G_c) \vee (\overline{R}_c \oplus B_c \oplus \overline{G}_c)$$

И при добавлении цвета G к исходной раскраски граней плоского кубического графа:

$$(R_c \oplus B_c \oplus G_c) \vee (\overline{R}_c \oplus \overline{B}_c \oplus \overline{G}_c)$$
$$(R_c \oplus B_c \oplus \overline{G}_c) \vee (\overline{R}_c \oplus \overline{B}_c \oplus G_c)$$

Сказанное можно представить в виде следующей таблицы:

Таблица раскраски изометрических циклов.

|  | $c_1$ | $c_2$ | $c_3$ | $c_4$ | $c_5$ | $c_6$ | $c_0$ | Индуцирующие цветные 2-факторы |
|---|---|---|---|---|---|---|---|---|
| Исходная раскраска | G | B | G | R | B | W | W | $(R_c \oplus B_c \oplus G_c) \vee (\overline{R}_c \oplus \overline{B}_c \oplus \overline{G}_c)$ |
| + R | B | G | B | W | G | R | R | $(\overline{R}_c \oplus B_c \oplus G_c) \vee (R_c \oplus \overline{B}_c \oplus \overline{G}_c)$ |
| + B | R | W | R | G | W | B | B | $(R_c \oplus \overline{B}_c \oplus G_c) \vee (\overline{R}_c \oplus B_c \oplus \overline{G}_c)$ |
| + G | W | R | W | B | R | G | G | $(R_c \oplus B_c \oplus \overline{G}_c) \vee (\overline{R}_c \oplus \overline{B}_c \oplus G_c)$ |

С другой стороны, можно сказать, что раскраска граней плоского кубического графа в четыре цвета порождает (индуцирует) раскраску ребер в три цвета для плоского кубического графа H согласно операции преобразования группы Клейна. Кроме того будем отождествлять раскраску граней плоского графа и раскраску изометрических базисных циклов. На рис. 1.15 представлен раскрашенный кубический граф его раскраска ребер, раскраска граней и раскраска изометрических циклов.

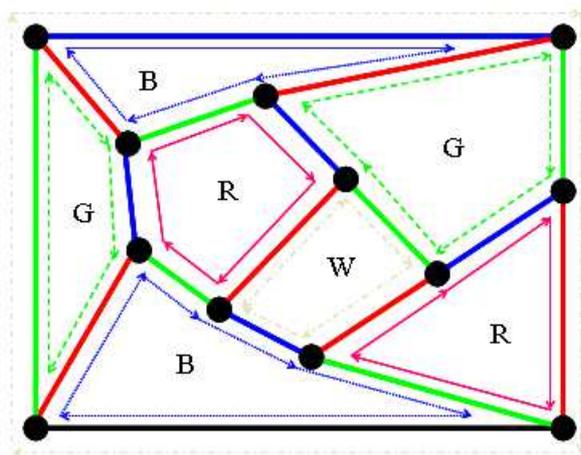

Рис. 1.15. Цветные изометрические циклы и цветные грани.

## 1.6. Ротация цветных дисков

Под **ротацией** цветного гамильтонова диска будем подразумевать изменение последовательности раскраски рёбер данного диска. Ротация диска приводит к изменению других цветных гамильтоновых квазициклов [9]. На рис. 1.16,а представлен синий диск до ротации, а на рис. 1.16,б после ротации.

Операция ротации дисков может быть описана как перестановка рёбер в других цветных 2-факторах. Например, для графа, представленного на рис. 1.11, цветные 2-факторы до ротации синего диска имеют вид:

синий 2-фактор - $\{e_1, e_3, e_5, e_6, e_7, e_8, e_{10}, e_{11}, e_{14}, e_{15}\}$,

красный 2-фактор - $\{e_1, e_2, e_4, e_6, e_9, e_{10}, e_{11}, e_{12}, e_{13}, e_{14}\}$,

зелёный 2-фактор - $\{e_2, e_3, e_4, e_5, e_7, e_8, e_9, e_{12}, e_{13}, e_{15}\}$.

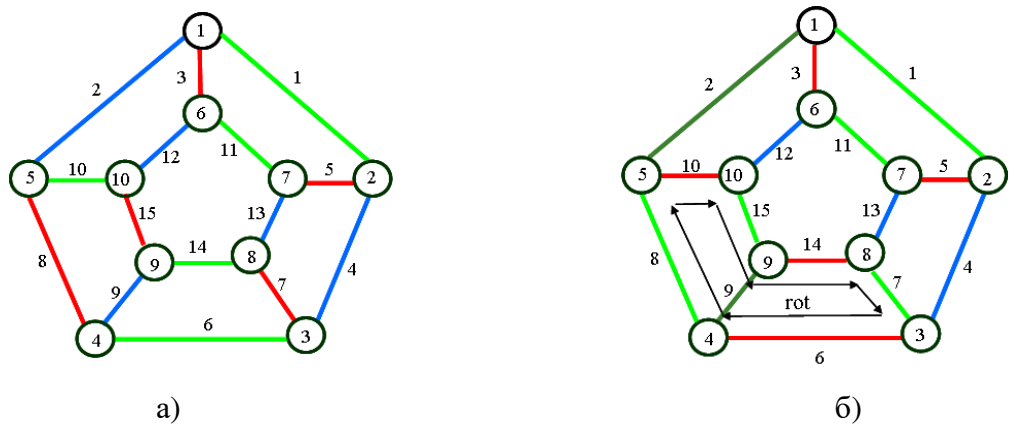

а) б)

Рис. 1.16. Раскраска графа H до ротации синего диска а), и после ротации б).

После ротации синего диска:

синий 2-фактор - $\{e_1, e_3, e_5, e_6, e_7, e_8, e_{10}, e_{11}, e_{14}, e_{15}\}$,

красный 2-фактор - $\{e_1, e_2, e_4, e_7, e_8, e_9, e_{11}, e_{12}, e_{13}, e_{15}\}$,

зеленый 2-фактор - $\{e_2, e_3, e_4, e_5, e_6, e_9, e_{10}, e_{12}, e_{13}, e_{14}\}$.

В красном и зеленом 2-факторах произошла перестановка следующих пар ребер ($e_6 \leftrightarrow e_7$), ($e_{10} \leftrightarrow e_8$), ($e_{14} \leftrightarrow e_{15}$).

### 1.7. Рекурсивный принцип построения плоских кубических графов

Что касается построения кубического графа без мостов на плоскости, то любой кубический граф может быть построен рекурсивным способом путем введения новых ребер в предыдущий граф для получения последующего кубического графа следующим образом [11]:

***Способ 1.*** Расположением двух новых вершин на двух ребрах принадлежащих одному изометрическому циклу исходного кубического графа без мостов и проведения нового ребра соединяющего эти вершины (см. рис. 1.17).

***Способ 2.*** Расположением двух новых вершин на одном ребре исходного кубического графа без мостов и проведения нового ребра, соединяющего эти вершины (см. рис. 1.18).

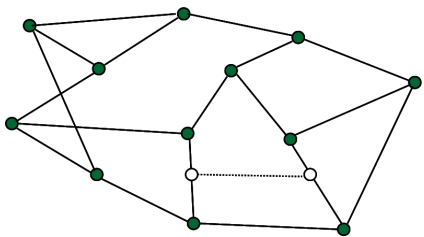   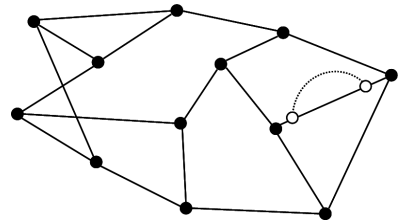

Рис. 1.17. Первый способ построения кубического графа H**.**     Рис. 1.18. Второй способ построения кубического графа H**.**

Минимальным графом для первоначального построения является кубический граф с двумя вершинами и тремя ребрами (см. рис. 1.20).

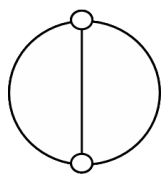

Рис. 1.19. Минимальный кубический граф.

Любой кубический граф может быть построен из предыдущего кубического графа без мостов путем введения нового ребра одним из двух указанных выше способов. Нетрудно заметить, что цвет вновь вводимого ребра соответствует цвету гамильтонового диска (диска 2-фактора), на котором находятся концы вновь вводимого ребра. И действительно, вновь введенное ребро увеличивает количество вершин диска ровно на 2, что естественно не влияет на раскраску диска. Это относится также и к удаленным ребрам, так как удаление ребра уменьшает количество вершин в диске ровно на 2. Если ребро введено или удалено с учетом последнего факта, то в результате получим правильную раскраску нового кубического графа [9].

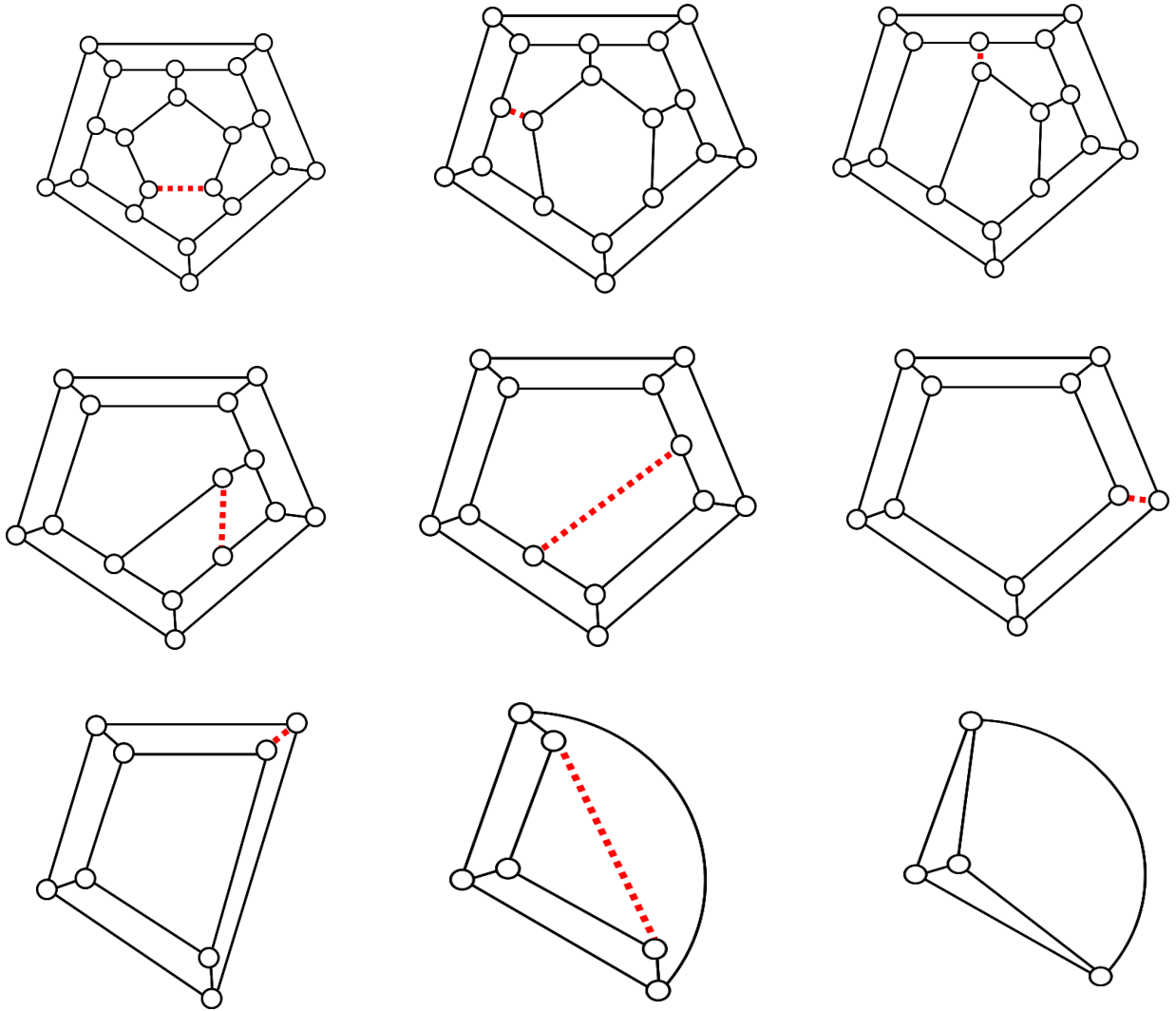

Рис. 1.20. Последовательный процесс удаления ребер из исходного кубического графа.

Плоский кубический граф без мостов до введения нового ребра будем называть *предыдущим*. Плоский кубический граф без мостов, полученный после введения нового ребра, будем называть *последующим*. Ребра, на которых находятся вершины вновь вводимого ребра для предыдущего графа, впредь будем называть *сцепленными ребрами*. Для получения последующего плоского кубического графа сцепленные ребра обязательно должны принадлежать изометрическому циклу.

С учетом описанного, раскраску любого плоского кубического графа можно свести к последовательному процессу:

- удаление ребер из исходного плоского кубического графа без мостов с целью получения плоского кубического графа, раскраска которого известна;
- раскраска полученного плоского кубического графа без мостов;
- определение сцепленных ребер для введения нового ребра;

- поиск цветного диска проходящего по сцепленным ребрам;
- раскраски вновь введенного ребра в цвет диска, принадлежащего сцепленным ребрам;
- перекраски ребер диска, принадлежащего сцепленным ребрам после введения нового ребра.

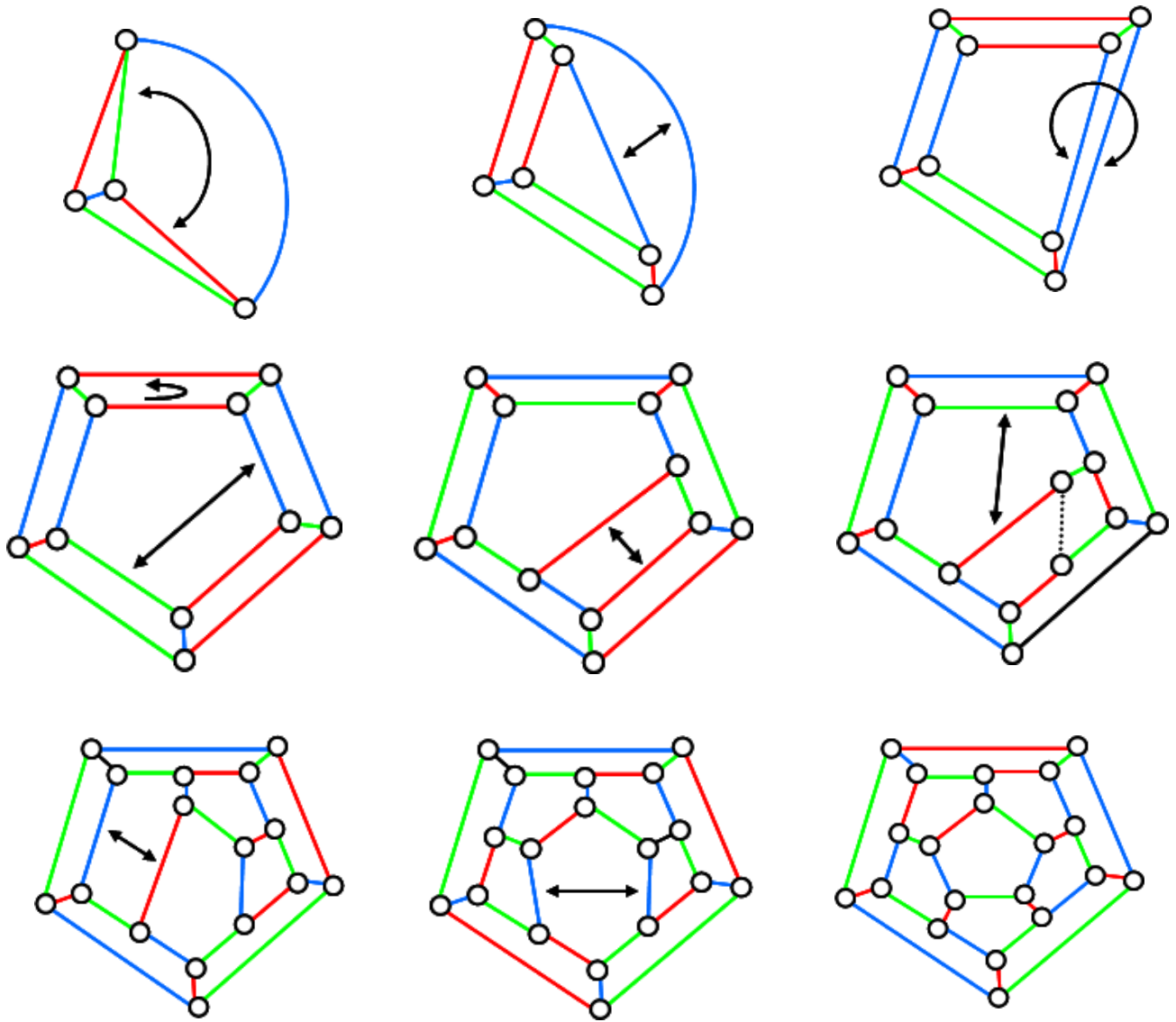

Рис. 1.21. Последовательный процесс введения нового ребра и перекраска ребер.

На рис. 1.20 представлен последовательный процесс удаления ребер из исходного кубического графа без мостов до получения известного правильно раскрашенного плоского кубического графа (здесь удаляемые ребра представлены красным цветом).

На рис. 1.21 представлен последовательный процесс правильной раскраски из известного правильно раскрашенного плоского кубического графа без мостов до получения правильной раскраски исходного кубического графа.

Так как по каждому из сцепленных ребер проходит два разных цветных диска, то обязательно имеется хотя бы один цветной 2-фактор принадлежащий сцепленный ребрам. И тогда возможны следующие комбинации:

- сцепленные ребра принадлежат одному цветному диску соответствующего цветного 2-фактора;
- сцепленные ребра принадлежат двум разным дискам соответствующего цветного 2-фактора.

Если сцепленные ребра одного цвета, то выбор цветного 2-фактора тривиален из двух других цветных 2-факторов. Если ребра различного цвета, то общий цветной 2-фактор будет третьего цвета.

В первом случае (тривиальный случай) вновь введенное ребро может быть легко раскрашено. Покажем, что оно может быть раскрашено и во втором случае.

## Выводы

Из рассмотрения свойств раскраски плоского кубического графа следует, что существует тесная связь между раскраской ребер и раскраской граней. Раскраска ребер в плоском кубическом графе индуцирует (порождает) раскраску граней, а раскраска граней, в свою очередь, индуцирует (порождает) раскраску ребер плоского кубического графа.

В раскрашенном кубическом графе обязательно присутствуют три цветных 2-фактора и сложение цветов осуществляется по законам преобразования группы Клейна. Цветной 2-фактор в раскрашенном кубическом графе состоит из дисков четной длины. Для раскрашенных кубических графов введена операция – ротация цветных дисков и установлена связь между раскраской ребер и раскраской независимой системы изометрических циклов графа. Показан рекурсивный принцип построения кубических графов. Рассмотрены вопросы построения произвольного плоского кубического графа без мостов и процесс введения (удаления) ребра в предыдущем графе для получения нового последующего плоского кубического графа. Показано, что для раскраски плоского кубического графа необходимо существовании цветного диска проходящего по двум сцепленным ребрам, принадлежащим изометрическому циклу плоского кубического графа.

## Глава 2. Теорема о существовании цветного диска, проходящего по сцепленным ребрам

### 2.1. Теорема о существовании цветного диска, проходящего по сцепленным ребрам

Таким образом нужно доказать теорему о существовании цветного диска, проходящего по сцепленным ребрам принадлежащим базисному циклу в произвольном плоском кубическом раскрашенном графе без мостов. Будем обозначать предыдущий граф до введения ребра $G_{pr}$, а последующий граф после введения ребра $G_{se}$.

Если мы докажем эту теорему, то следствием данной теоремы будет являться известная теорема о четырех красках для максимально плоского графа [1,17]. Гипотеза о четырех красках была подтверждена работами К.Аппеля и В.Хакена [18,19], проведенными ими комбинаторными методами с помощью вычислительных машин. Мы же приведем дедуктивный способ доказательства теоремы о существовании цветного диска, проходящего по сцепленным ребрам и, как следствие, получим утверждение о четырех красках [9].

Рассмотрим процесс последовательного введения ребра в плоский кубический граф $G_{pr}$ для получения нового плоского кубического графа $G_{se}$. Пусть, на каком-то шаге рекурсивного последовательного процесса для введения нового ребра, имеется раскрашенный кубический граф. Введем новое ребро и определим сцепленные ребра. Если по сцепленным ребрам проходит один цветной диск, то раскраска ребер в последующем графе тривиальна (см. рис. 2.1).

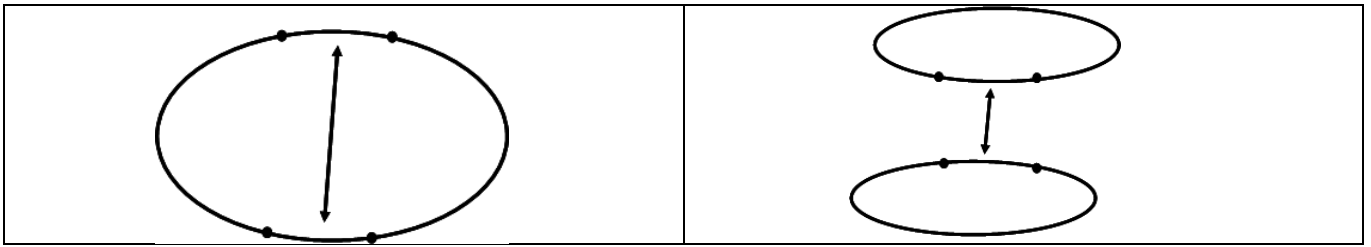

Рис. 2.1. Тривиальный случай раскраски ребер и особый случай раскраски ребер.

*Лемма 2.1*. *Если в предыдущем раскрашенном кубическом графе без мостов $G_{pr}$ существует цветной диск, проходящий по двум сцепленным ребрам, то последующий кубический граф $G_{se}$ раскрашен.*

*Доказательство*. По теореме Петерсена последующий кубический граф $G_{se}$ состоит из 1-фактора и 2-фактора. Если в предыдущем раскрашенном кубическом графе $G_{pr}$ существует цветной диск, проходящий по сцепленным ребрам, то соответствующий цветной диск, проходящий по вершинам вновь введенного ребра в графе $G_{se}$, увеличивается на две единицы и может быть раскрашен в цвет соответствующего диска в графе $G_{pr}$. Естественно, что вновь введенное ребро в $G_{se}$ принадлежит 1-фактору и тоже раскрашивается в цвет диска, проходящего по вершинам введенного ребра. Следовательно, в этом случае, граф $G_{se}$ – раскрашиваем.

Рассмотрим нетривиальный случай, когда сцепленные рёбра принадлежат двум одинаково раскрашенным дискам (см. рис. 2.2).

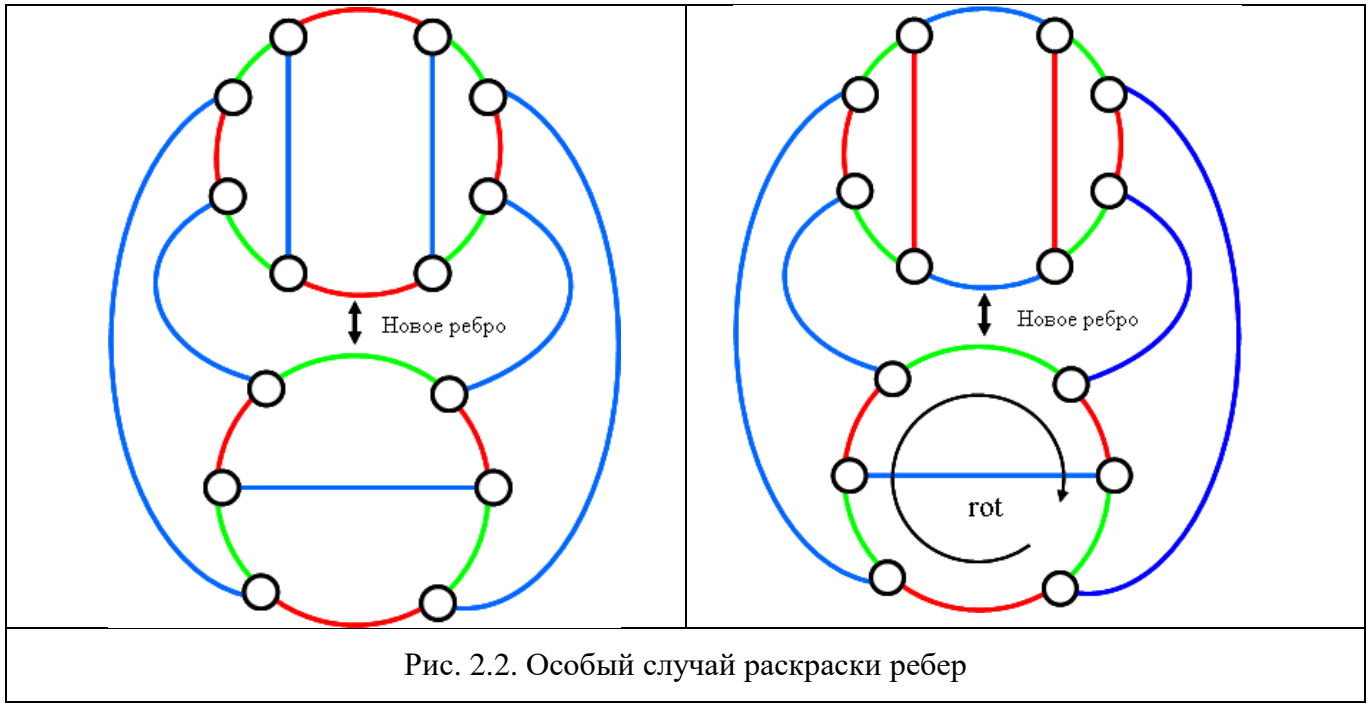

Рис. 2.2. Особый случай раскраски ребер

*Лемма 2*. *Если в предыдущем раскрашенном кубическом графе без мостов $G_{pr}$ не существует цветного диска, проходящего по двум сцепленным рёбрам, то последующий кубический граф $G_{se}$ без мостов не плоский и не раскрашиваемый.*

*Доказательство*. Предположим, что в предыдущем раскрашенном кубическом графе $G_{pr}$ не существует цветного диска, проходящего по сцепленным рёбрам. Тогда существует два одинаково раскрашенных цветных диска, проходящих по сцепленным ребрам. После введения нового ребра в последующий граф $G_{se}$, в соответствующем диске какое-то рядом стоящее ребро должно быть окрашено в белый цвет для соответствия правилам раскраски. Так как в этом случае каждый диск, проходящий по вершине нового ребра, увеличится в длину на одну единицу и станет диском нечетной длины (рис. 2.3). Но тогда ребро, окрашенное в белый цвет, должно индуцироваться тремя разноцветными базисными циклами для удовлетворения группового преобразования Клейна:

$$W = (R + B + G).$$

Но последнее высказывание нарушает условие планарности (напомним, что по теореме Маклейна, граф планарен тогда и только тогда, когда существует базис циклов такой, что по каждому ребру проходит не более двух циклов). Следовательно, рассматриваемый граф $G_{pr}$ не плоский.

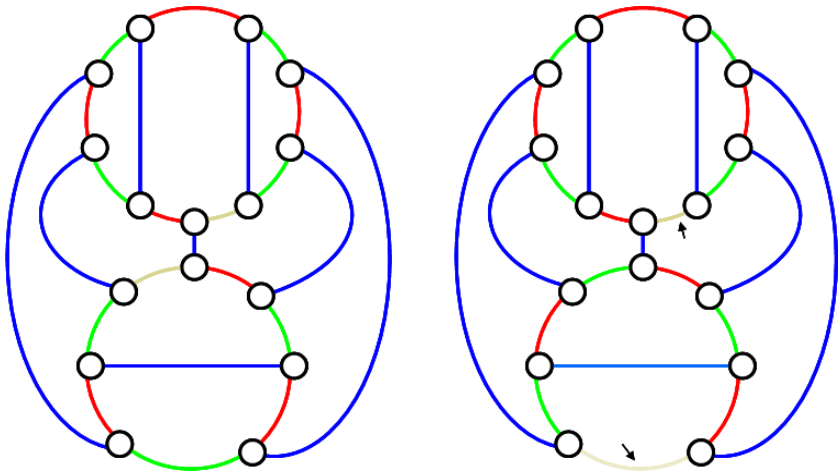 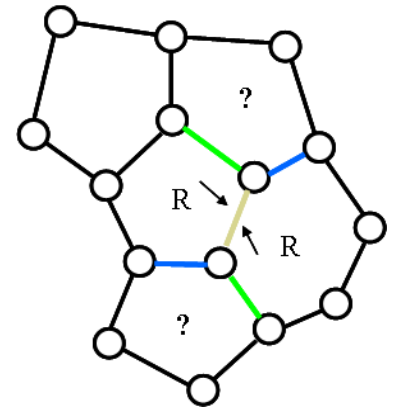

Рис. 2.3. Введение белых ребер в последующий плоский граф (на рисунке белые ребра указаны стрелками).

Рис. 2.4. Случай раскраски смежных циклов в один цвет.

Рассмотрим случай, когда $W = R + R = B + B = G + G$.

Если предположить, что белое ребро индуцировано двумя одинаково раскрашенными смежными циклами (гранями), то тогда возникает противоречие в раскраске циклов (граней) смежных с раскрашенными, так как эти циклы должны со стороны одного ребра быть раскрашены в один цвет, а со стороны другого ребра - в другой цвет (эти циклы на рис. 2.4 обозначены вопросом).

Мы пришли к противоречию.

Следовательно, наше предположение об отсутствии цветного диска проходящего по сцепленным ребрам не верно. Таким образом, лемма доказана.

Данная лемма утверждает, что появление ребра белого цвета в последующем кубическом графе $G_{se}$ характеризует не плоский и не раскрашенный кубический граф.

Теперь докажем следующую теорему.

**Теорема 2.3**[9]. *В плоском раскрашенном кубическом графе без мостов $G_{pr}$, существует цветной диск, проходящий по сцепленным ребрам.*

**Доказательство**. Если вновь введенное ребро принадлежит базисному циклу предыдущего раскрашенного плоского кубического графу, то последующий кубический граф $G_{se}$ является плоским по построению. Но тогда в плоском графе не возможно существование белого ребра по теореме Маклейна, не возможно существование трех циклов проходящих по ребру графа. Следовательно, ребра графа могут быть раскрашены только в три цвета R,B,G. А это значит, что в предыдущем плоском раскрашенном кубическом графе без мостов $G_{pr}$ существует цветной диск, проходящий по сцепленным ребрам. Понятие предыдущего раскрашенного плоского кубического графа можно распространить на все множество раскрашенных плоских кубических графов.

Теорема доказана.

Из данной теоремы вытекает следующие следствия:

***Следствие 2.3.1.*** *Вновь введенное ребро для последующего плоского кубического графа без мостов раскрашиваемо.*

Согласно доказанной теоремы вновь введенное ребро всегда можно раскрасить в цвет диска проходящего по сцепленным ребрам.

***Следствие 2.3.2.*** *Хроматический класс плоского кубического графа без мостов равен трем.*

Согласно доказанной теоремы из рекурсивности процесса построения раскраски, последующий граф раскрашивается путем введения и раскраски нового ребра с последующей перекраской ребер цветного диска проходящего по сцепленным ребрам.

***Следствие 2.3.3.*** *В плоском кубическом графе без мостов обязательно существует 2-фактор с дисками четной длины.*

Это следует из принципа процесса раскраски и доказанной теоремы**.**

***Следствие 2.3.4.*** *Хроматическое число максимально плоского графа равно четырем.*

Это следует из свойства 1.3.2 и теоремы 1.1 (5) и того факта, что изометрические циклы и обод плоского кубического графа H дуальны вершинам максимально плоского графа $\mathbf{G}'$.

В качестве демонстрации результатов наших рассуждений проведем пример для правильно раскрашенного непланарного графа представленного на рис. 2.5.

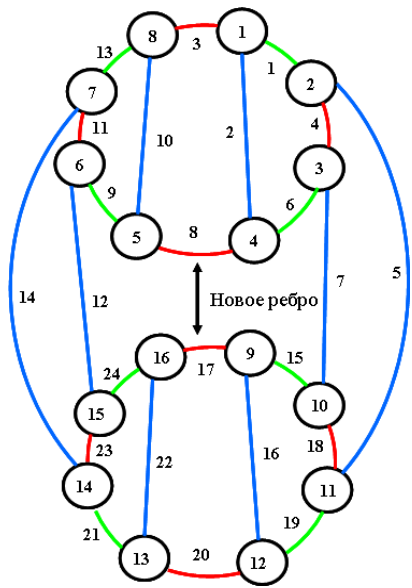
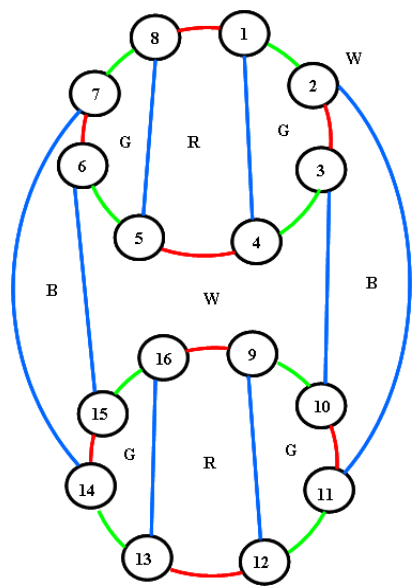

Рис. 2.5. Два диска одного цвета и введение нового ребра.　　　Рис. 2.6. Раскраска базисных циклов.

**Пример 2.1**. На рисунке 2.5 представлена раскраска кубического графа, где концы вновь введенного ребра принадлежат двум дискам синего цвета и двум дискам зеленого цвета.

Рассмотрим способы раскраски в этом особом случае.

Для данного кубического графа H (см. рис. 2.5) выберем базисные циклы

$c_1 = \{e_2, e_3, e_8, e_{10}\}$, $c_2 = \{e_1, e_2, e_4, e_6\}$, $c_3 = \{e_9, e_{10}, e_{11}, e_{13}\}$, $c_4 = \{e_4, e_5, e_7, e_{18}\}$,
$c_5 = \{e_{11}, e_{12}, e_{14}, e_{23}\}$, $c_6 = \{e_6, e_7, e_8, e_9, e_{12}, e_{15}, e_{17}, e_{24}\}$, $c_7 = \{e_{16}, e_{17}, e_{20}, e_{22}\}$,
$c_8 = \{e_{15}, e_{16}, e_{18}, e_{19}\}$, $c_9 = \{e_{21}, e_{22}, e_{23}, e_{24}\}$.

Формируем обод, $c_0 = \{e_1, e_3, e_5, e_{13}, e_{14}, e_{19}, e_{20}, e_{21}\}$.

Тогда цветные 2-факторы можно записать в виде:

$R_c = (c_2 \oplus c_4 \oplus c_8) \oplus (c_3 \oplus c_5 \oplus c_9) = (\{e_1, e_2, e_4, e_6\} \oplus \{e_4, e_5, e_7, e_{18}\} \oplus \{e_{15}, e_{16}, e_{18}, e_{19}\}) \oplus$
$\oplus (\{e_9, e_{10}, e_{11}, e_{13}\} \oplus \{e_{11}, e_{12}, e_{14}, e_{23}\} \oplus \{e_{21}, e_{22}, e_{23}, e_{24}\}) =$
$= \{e_1, e_2, e_5, e_6, e_7, e_9, e_{10}, e_{12}, e_{13}, e_{14}, e_{15}, e_{16}, e_{19}, e_{21}, e_{22}, e_{24}\};$
$B_c = (c_1 \oplus c_2 \oplus c_3) \oplus (c_7 \oplus c_8 \oplus c_9) = (\{e_2, e_3, e_8, e_{10}\} \oplus \{e_1, e_2, e_4, e_6\} \oplus \{e_9, e_{10}, e_{11}, e_{13}\}) \oplus$
$\oplus (\{e_{16}, e_{17}, e_{20}, e_{22}\} \oplus \{e_{15}, e_{16}, e_{18}, e_{19}\} \oplus \{e_{21}, e_{22}, e_{23}, e_{24}\}) =$
$= \{e_1, e_3, e_4, e_6, e_8, e_9, e_{11}, e_{13}, e_{15}, e_{17}, e_{18}, e_{19}, e_{20}, e_{21}, e_{23}, e_{24}\};$
$G_c = c_1 \oplus c_4 \oplus c_5 \oplus c_7 = \{e_2, e_3, e_8, e_{10}\} \oplus \{e_4, e_5, e_7, e_{18}\} \oplus \{e_{11}, e_{12}, e_{14}, e_{23}\} \oplus \{e_{16}, e_{17}, e_{20}, e_{22}\} =$
$= \{e_2, e_3, e_4, e_5, e_7, e_8, e_{10}, e_{11}, e_{12}, e_{14}, e_{16}, e_{17}, e_{18}, e_{20}, e_{22}, e_{23}\};$
$W_c = c_6 \oplus c_0 = \{e_6, e_7, e_8, e_9, e_{12}, e_{15}, e_{17}, e_{24}\} \oplus \{e_1, e_3, e_5, e_{13}, e_{14}, e_{19}, e_{20}, e_{21}\} =$
$= \{e_1, e_3, e_5, e_6, e_7, e_8, e_9, e_{12}, e_{13}, e_{14}, e_{15}, e_{17}, e_{19}, e_{20}, e_{21}, e_{24}\};$
$\overline{R}_c = c_1 \oplus c_6 \oplus c_7 \oplus c_0 = \{e_2, e_3, e_8, e_{10}\} \oplus \{e_6, e_7, e_8, e_9, e_{12}, e_{15}, e_{17}, e_{24}\} \oplus \{e_{16}, e_{17}, e_{20}, e_{22}\} \oplus$
$\oplus \{e_1, e_3, e_5, e_{13}, e_{14}, e_{19}, e_{20}, e_{21}\} = \{e_1, e_2, e_5, e_6, e_7, e_9, e_{10}, e_{12}, e_{13}, e_{14}, e_{15}, e_{16}, e_{19}, e_{21}, e_{22}, e_{24}\};$
$\overline{B}_c = c_4 \oplus c_5 \oplus c_6 \oplus c_0 = \{e_4, e_5, e_7, e_{18}\} \oplus \{e_{11}, e_{12}, e_{14}, e_{23}\} \oplus \{e_6, e_7, e_8, e_9, e_{12}, e_{15}, e_{17}, e_{24}\} \oplus$
$\oplus \{e_1, e_3, e_5, e_{13}, e_{14}, e_{19}, e_{20}, e_{21}\} = \{e_1, e_3, e_4, e_6, e_8, e_9, e_{11}, e_{13}, e_{15}, e_{17}, e_{18}, e_{19}, e_{20}, e_{21}, e_{23}, e_{24}\};$
$\overline{G}_c = c_2 \oplus c_3 \oplus c_6 \oplus c_8 \oplus c_0 = \{e_1, e_2, e_4, e_6\} \oplus \{e_9, e_{10}, e_{11}, e_{13}\} \oplus \{e_6, e_7, e_8, e_9, e_{12}, e_{15}, e_{17}, e_{24}\} \oplus$
$\oplus \{e_{15}, e_{16}, e_{18}, e_{19}\} \oplus \{e_1, e_3, e_5, e_{13}, e_{14}, e_{19}, e_{20}, e_{21}\} =$
$= \{e_2, e_3, e_4, e_5, e_7, e_8, e_{10}, e_{11}, e_{12}, e_{14}, e_{16}, e_{17}, e_{18}, e_{20}, e_{22}, e_{23}\}.$

Осуществим операцию ротации дисков. Будем вращать два простых диска зеленого цвета $c_1 = \{e_2, e_3, e_8, e_{10}\}$ и $c_7 = \{e_{16}, e_{17}, e_{20}, e_{22}\}$.

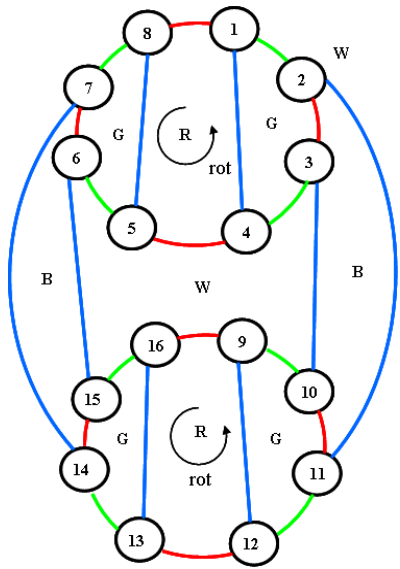

Рис. 2.7. Новая раскраска ребер и базисных циклов.

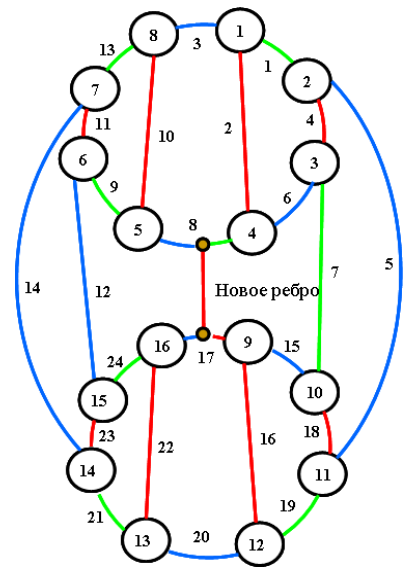

Рис. 2.8. Новая раскраска ребер и базисных циклов.

Новая раскраска характеризуется другим набором цветных 2-факторов и 1-факторов и позволяет раскрасить вновь введенное ребро красным цветом.

$R_c = c_1 \oplus c_2 \oplus c_3 \oplus c_4 \oplus c_5 \oplus c_7 \oplus c_8 \oplus c_9 = \{e_2,e_3,e_8,e_{10}\} \oplus \{e_1,e_2,e_4,e_6\} \oplus \{e_9,e_{10},e_{11},e_{13}\} \oplus$
$\oplus \{e_4,e_5,e_7,e_{18}\} \oplus \{e_{11},e_{12},e_{14},e_{23}\} \oplus \{e_{16},e_{17},e_{20},e_{22}\} \oplus \{e_{15},e_{16},e_{18},e_{19}\} \oplus$
$\oplus \{e_{21},e_{22},e_{23},e_{24}\} = \{e_1,e_3,e_5,e_6,e_7,e_8,e_9,e_{12},e_{13},e_{14},e_{15},e_{17},e_{19},e_{20},e_{21},e_{24}\};$

$B_c = c_2 \oplus c_3 \oplus c_8 \oplus c_9 = \{e_1,e_2,e_4,e_6\} \oplus \{e_9,e_{10},e_{11},e_{13}\} \oplus \{e_{15},e_{16},e_{18},e_{19}\} \oplus \{e_{21},e_{22},e_{23},e_{24}\}) =$
$= \{e_1,e_2,e_4,e_6,e_9,e_{10},e_{11},e_{13},e_{15},e_{16},e_{18},e_{19},e_{21},e_{22},e_{23},e_{24}\};$
$G_c = c_1 \oplus c_4 \oplus c_5 \oplus c_7 = \{e_2,e_3,e_8,e_{10}\} \oplus \{e_4,e_5,e_7,e_{18}\} \oplus \{e_{11},e_{12},e_{14},e_{23}\} \oplus \{e_{16},e_{17},e_{20},e_{22}\} =$
$= \{e_2,e_3,e_4,e_5,e_7,e_8,e_{10},e_{11},e_{12},e_{14},e_{16},e_{17},e_{18},e_{20},e_{22},e_{23}\}.$

Применение операции ротации двух зеленных дисков изменяет цвета базисных циклов с красного на синий (см. рис. 2.7), и тогда вновь вводимое ребро можно окрасить в красный цвет (см. рис. 2.8).

**Пример 2.2.** Введем новое ребро для сцепленных ребер ($e_2,e_{12}$) с целью получения последующего кубического графа. В данном непланарном графе не существует цветного диска, проходящего по данным сцепленным ребрам. Поэтому раскрасим ребро $e_2$ и ребро $e_{15}$ нового графа белым цветом (см. рис. 2.9).

Множество изометрических циклов для нового графа состоит из двенадцати циклов длиной пять:

$c_1 = \{e_1,e_2,e_4,e_6,e_8\}$, $c_2 = \{e_1,e_2,e_5,e_{10},e_{14}\}$, $c_3 = \{e_1,e_3,e_4,e_7,e_{11}\}$, $c_4 = \{e_1,e_3,e_5,e_{12},e_{13}\}$,
$c_5 = \{e_2,e_3,e_{10},e_{11},e_{15}\}$, $c_6 = \{e_2,e_3,e_8,e_9,e_{12}\}$, $c_7 = \{e_4,e_5,e_6,e_9,e_{13}\}$, $c_8 = \{e_4,e_5,e_7,e_{14},e_{15}\}$,
$c_9 = \{e_6,e_7,e_9,e_{11},e_{12}\}$, $c_{10} = \{e_6,e_7,e_8,e_{10},e_{15}\}$, $c_{11} = \{e_8,e_9,e_{10},e_{13},e_{14}\}$, $c_{12} = \{e_{11},e_{12},e_{13},e_{14},e_{15}\}$.

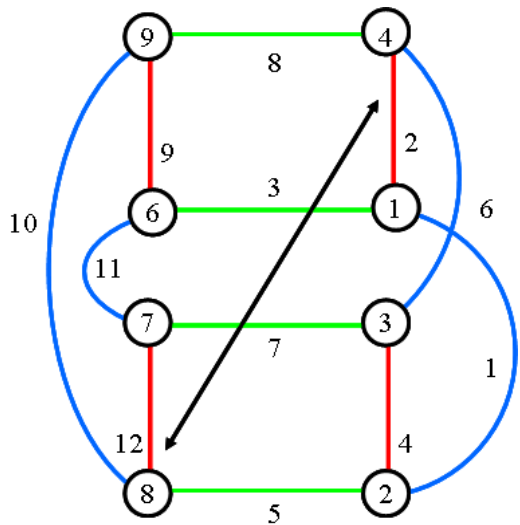

Рис. 2.9. Раскрашенный непланарный граф и вновь введенное ребро.

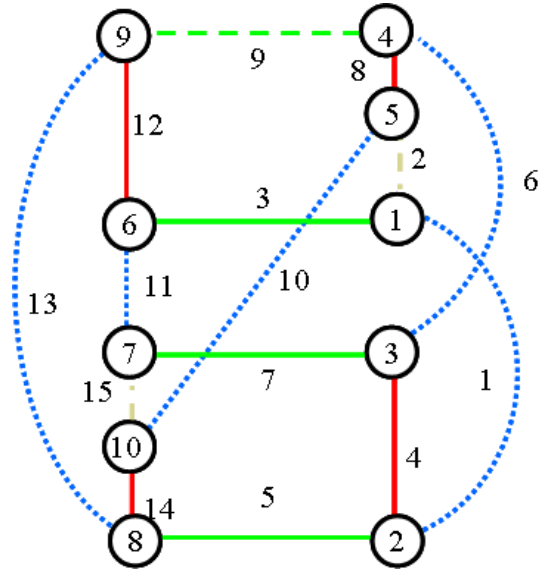

Рис. 2.10. Раскраска ребер $e_2$ и $e_{15}$ в белый цвет.

Если в качестве базиса выбрать следующие циклы:

$c_1 = \{e_1, e_2, e_4, e_6, e_8\}$, $c_3 = \{e_1, e_3, e_4, e_7, e_{11}\}$, $c_5 = \{e_2, e_3, e_{10}, e_{11}, e_{15}\}$, $c_6 = \{e_2, e_3, e_8, e_9, e_{12}\}$,

$c_8 = \{e_4, e_5, e_7, e_{14}, e_{15}\}$, $c_{12} = \{e_{11}, e_{12}, e_{13}, e_{14}, e_{15}\}$,

то получим, что по ребрам $e_2$ и $e_{15}$ проходит три изометрических цикла (см. рис. 2.10).

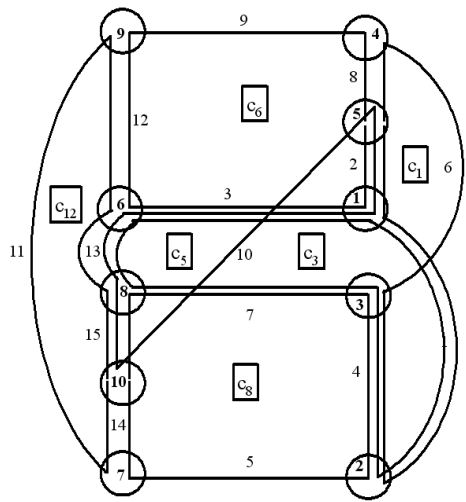

Рис. 2.11. Базисные изометрические циклы непланарного графа.

Если перерисовать рисунок данного графа, то мы получим более известный рисунок этого графа Петерсена (см. рис. 1.13,а), у которого хроматический класс не равен трем [17]. Это построение демонстрирует доказанную теорему 2.3.

**Выводы**

Таким образом, утверждается, что для плоского кубического графа раскраска ребер и раскраска граней является единым целым. Раскраска ребер порождает раскраску граней и наоборот, раскраска граней порождает раскраску ребер. Причем, раскраску ребер можно осуществ-лять рекурсивно, последовательно добавляя ребро в предыдущий раскрашенный плоский кубический граф. Как правило, в большинстве случаев по сцепленным ребрам проходит один цветной диск, и тогда раскраска ребра производится в цвет диска. Однако при рекурсивном построении может оказаться, что по сцепленным ребрам проходят два цветных диска окрашенных в один цвет. В этом случае, также возможна раскраска вновь введенного ребра, используя операцию ротации цветного диска. В общем случае поиск таких необходимых цветных дисков является искусством и осуществляется визуально при рассмотрении рисунка графа. Поэтому стоит рассмотреть алгоритмы перекраски плоского кубического графа для осуществления процесса раскраски вновь введенного ребра.

# Глава 3. Алгоритмы для рекурсивной раскраски ребер плоского кубического графа
## 3.1. Рекурсивный метод раскараски

**Пример 3.1.** Рассмотрение процесса раскраски ребер плоских кубических графов начнем с графа Тата.

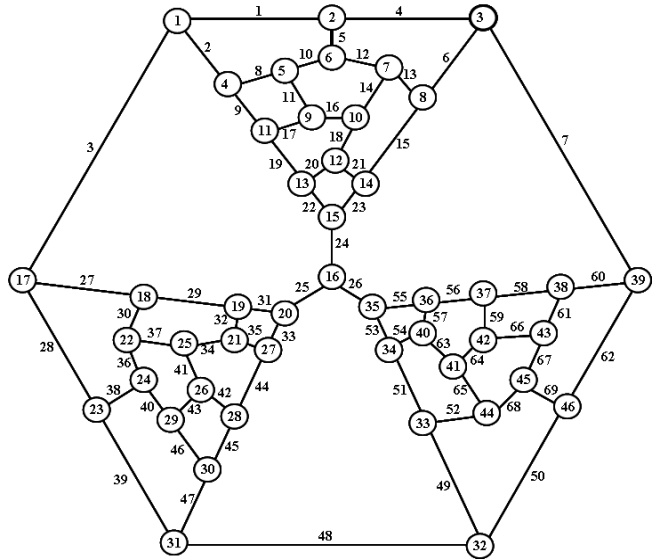

Рис. 3.1. Граф Тата.

Начнем последовательно удалять ребра из исходного графа. Удаление ребер будем производить таким образом, что бы ни в коем случае не получить сепарабельный граф с мостом.

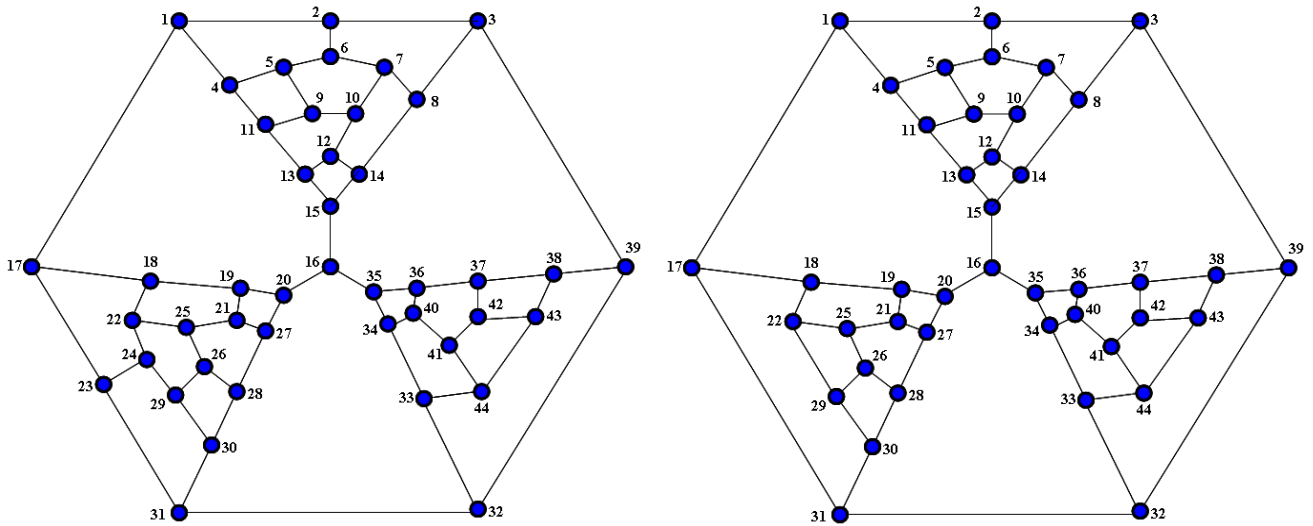

Рис. 3.3. Удаление из графа ребра ($v_{23}, v_{24}$).

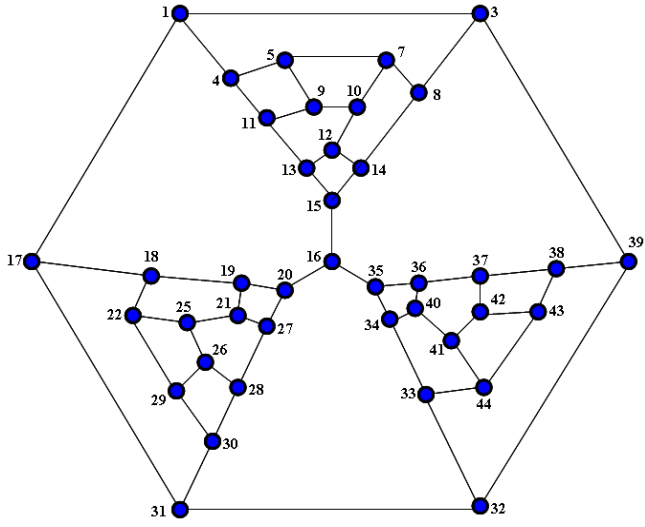
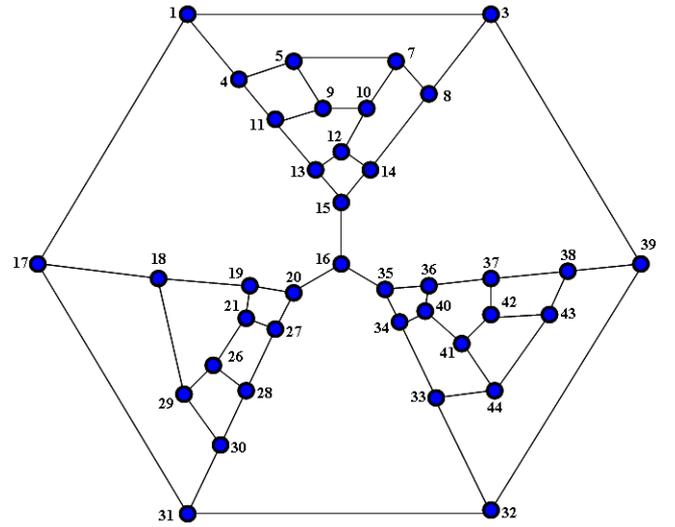
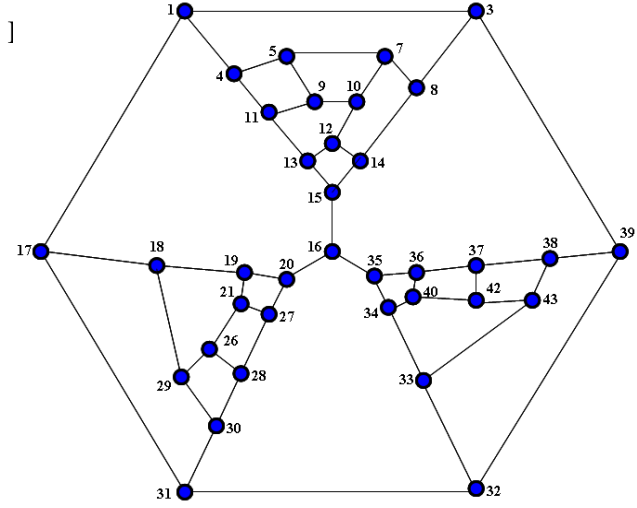
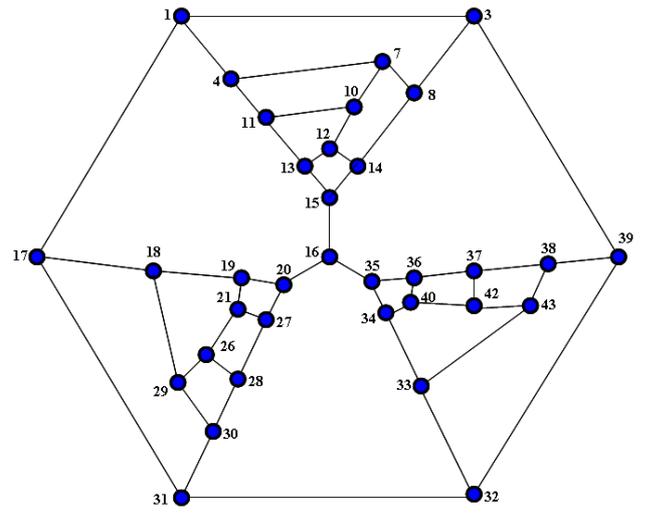
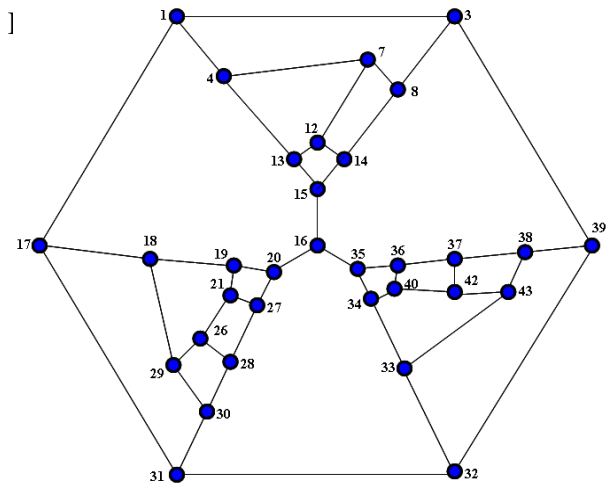
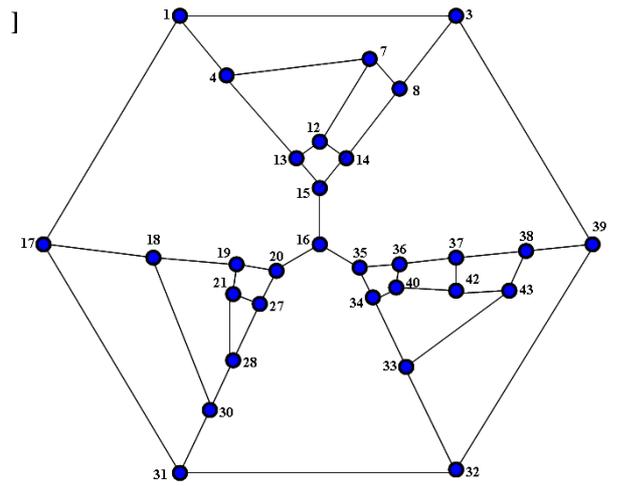

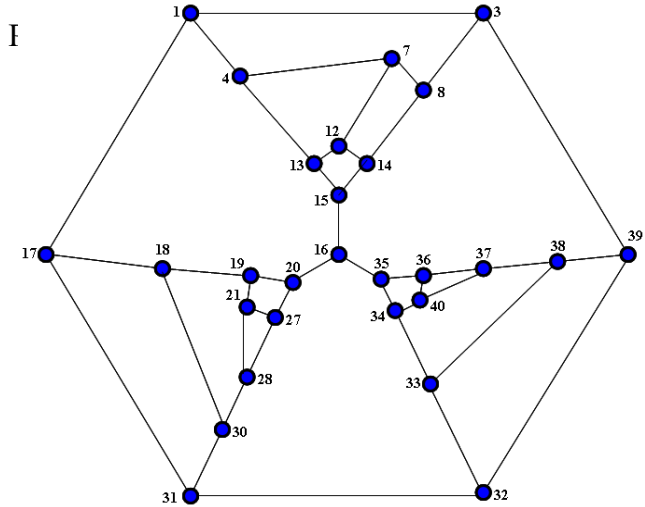
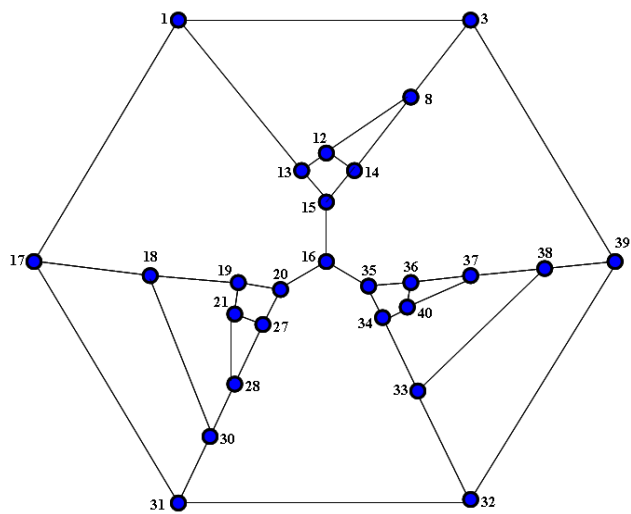

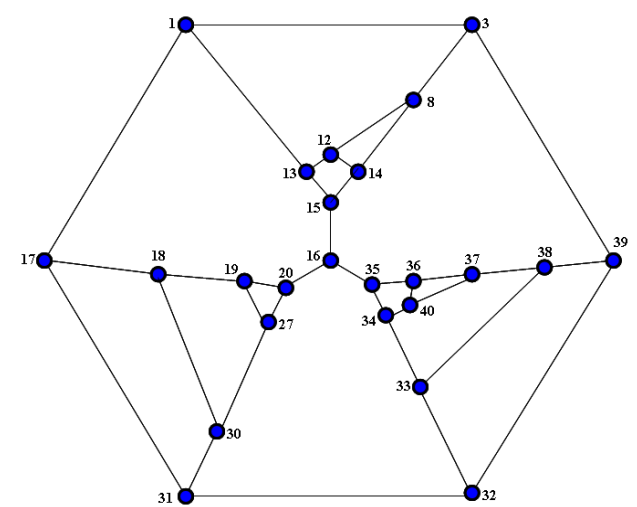
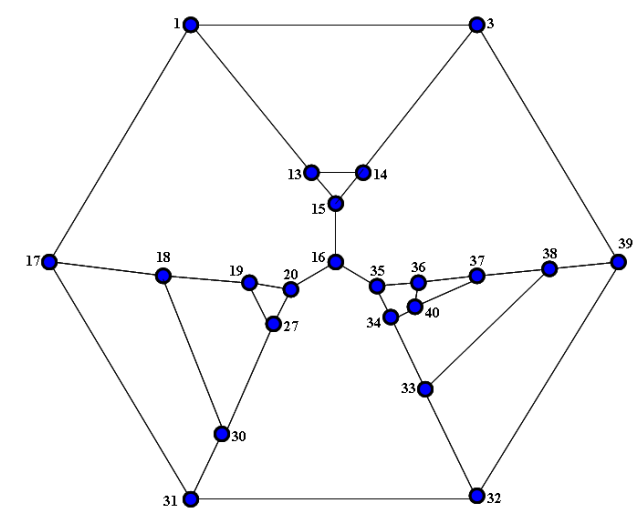

Рис. 3.12. Удаление из графа ребра ($v_{21}, v_{28}$). Рис. 3.13. Удаление из графа ребра ($v_8, v_{12}$).

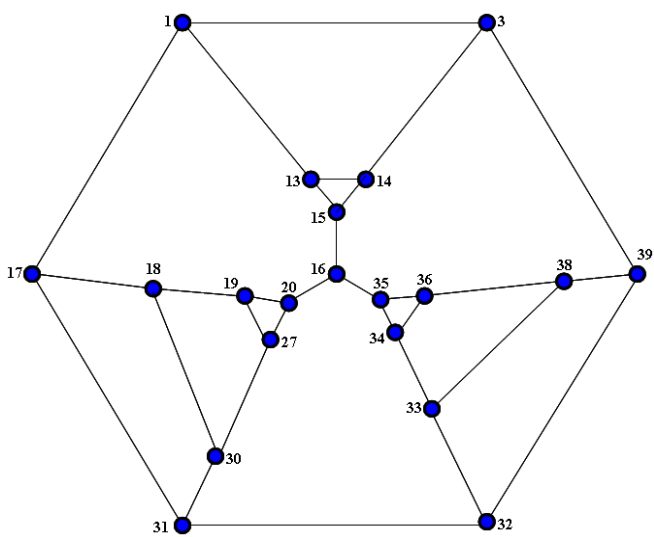
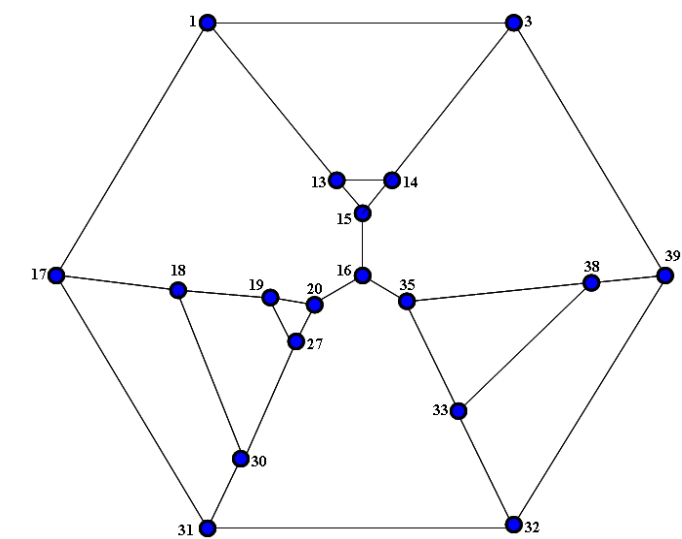

Рис. 3.14. Удаление из графа ребра ($v_{37}, v_{40}$).

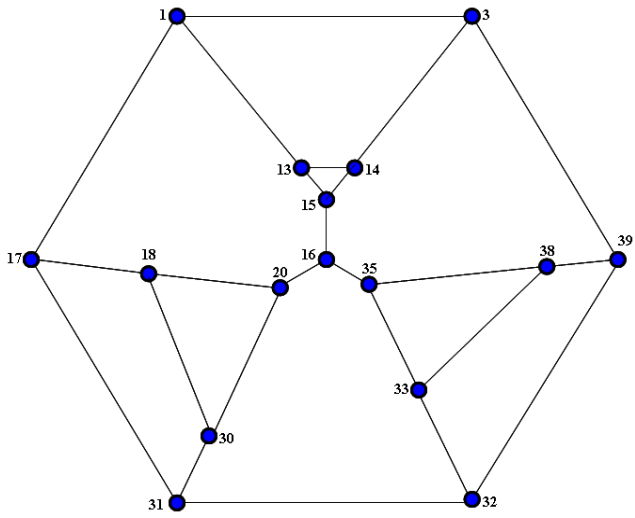

Рис. 3.16. Удаление из графа ребра ($v_{19}, v_{27}$).

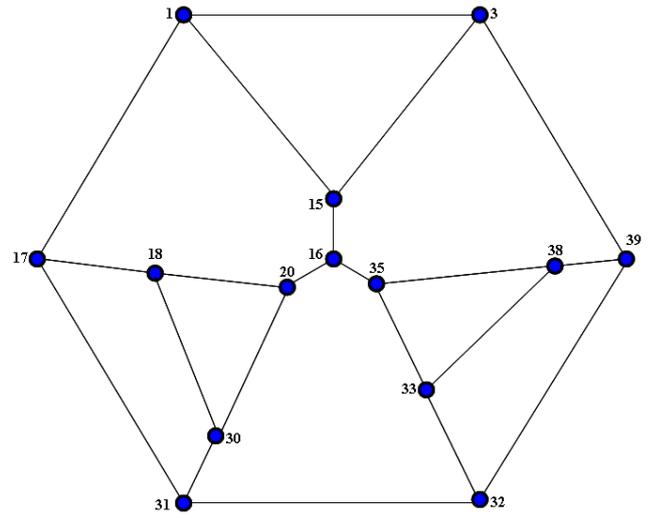

Рис. 3.17. Удаление из графа ребра ($v_{13}, v_{14}$).

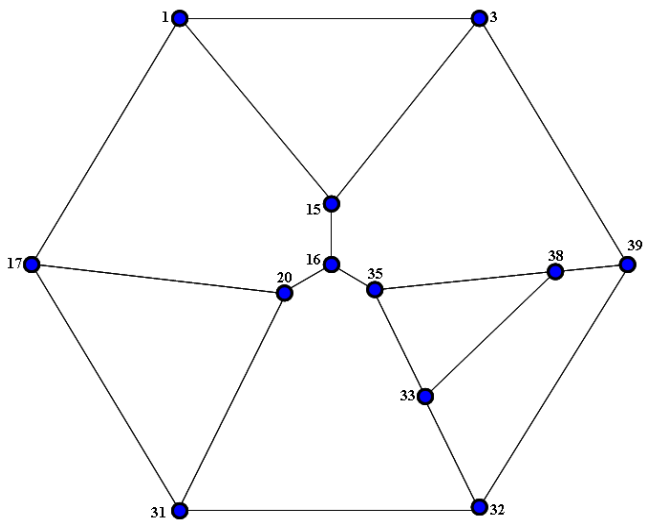

Рис. 3.18. Удаление из графа ребра ($v_{18}, v_{30}$).

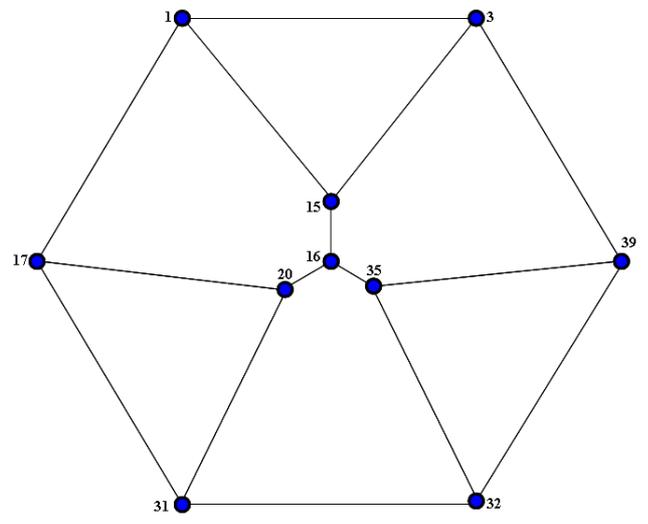

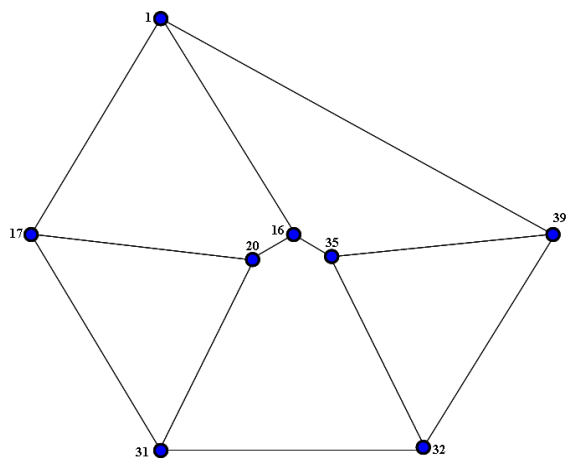

Рис. 3.20. Удаление из графа ребра ($v_3, v_{15}$).

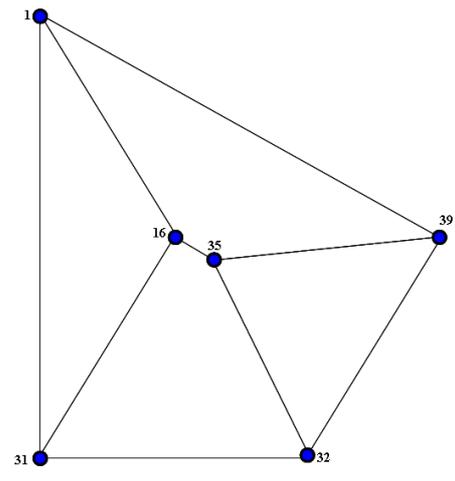

Рис. 3.21. Удаление из графа ребра ($v_{17}, v_{20}$).

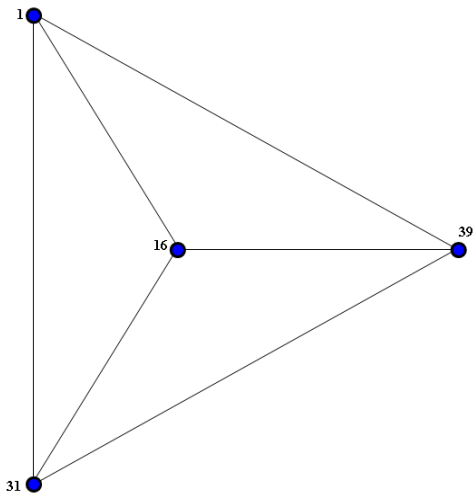 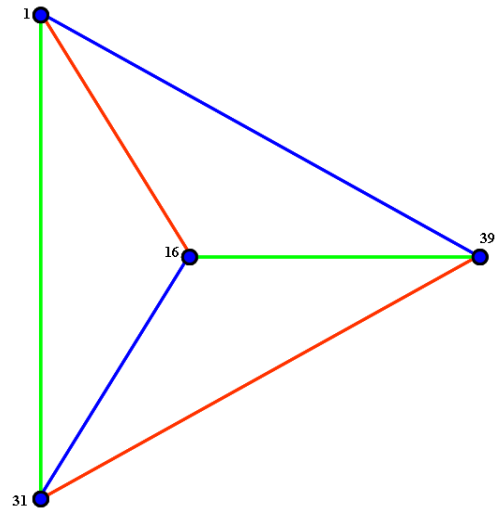

Рис. 3.22. Удаление из графа ребра $(v_{32}, v_{35})$.      Рис. 3.23. Раскраска ребер.

Процесс последовательного удаления ребер из графа представлен на рис. 3.2-3.22. На рис. 3.23 представлена раскраска ребер для известного плоского кубического графа $К_4$ без мостов.

Теперь начнем раскраску ребер согласно ранее построенной цепочки удаленных ребер.

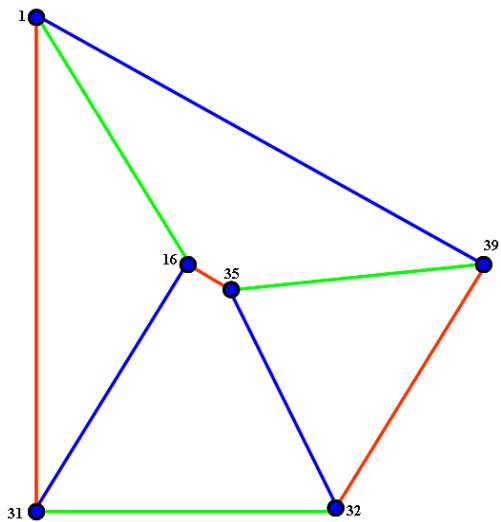 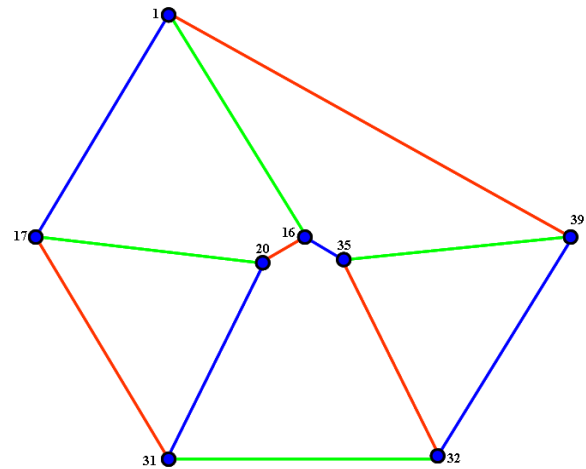

Рис. 3.24. Раскраска ребра $(v_{32}, v_{35})$.      Рис. 3.25. Раскраска ребра $(v_{17}, v_{20})$.

Раскрашиваем ребро $(v_{32}, v_{35})$ в синий цвет, так как сцепленные ребра принадлежат синему цветному диску (см. рис. 3.24).

Раскрашиваем ребро $(v_{17}, v_{20})$ в зеленый цвет, так как сцепленные ребра принадлежат зеленному цветному диску (см. рис. 3.25).

Раскрашиваем ребро $(v_3, v_{15})$ в синий цвет, так как сцепленные ребра принадлежат синему цветному диску (см. рис. 3.26).

Раскрашиваем ребро $(v_{33}, v_{38})$ в синий цвет, так как сцепленные ребра принадлежат синему цветному диску (см. рис. 3.27).

Рис. 3.26. Раскраска ребра ($v_3, v_{15}$).

Рис. 3.27. Раскраска ребра ($v_{33}, v_{38}$).

Рис. 3.28. Раскраска ребра ($v_{18}, v_{30}$).

Рис. 3.29. Раскраска ребра ($v_{13}, v_{14}$).

Раскрашиваем ребро ($v_{18}, v_{30}$) в красный цвет, так как сцепленные ребра принадлежат красному цветному диску (см. рис. 3.28).

Раскрашиваем ребро ($v_{13}, v_{14}$) в синий цвет, так как сцепленные ребра принадлежат синему цветному диску (см. рис. 3.29).

Продолжаем последовательный процесс раскраски ребер.

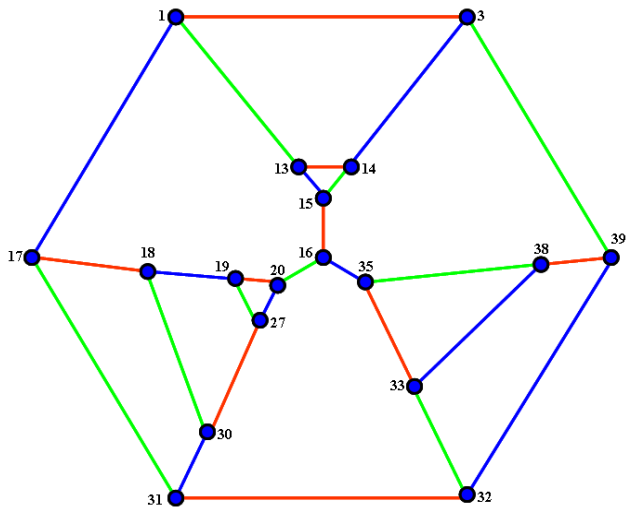

Рис. 3.30. Раскраска ребра ($v_{19}, v_{27}$).

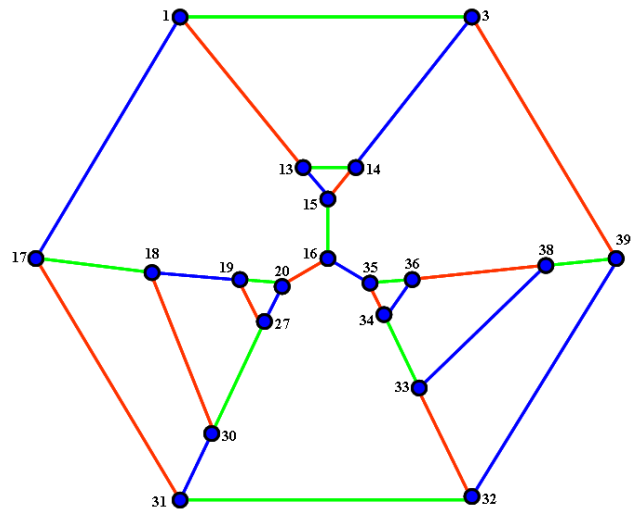

Рис. 3.31. Раскраска ребра ($v_{34}, v_{36}$).

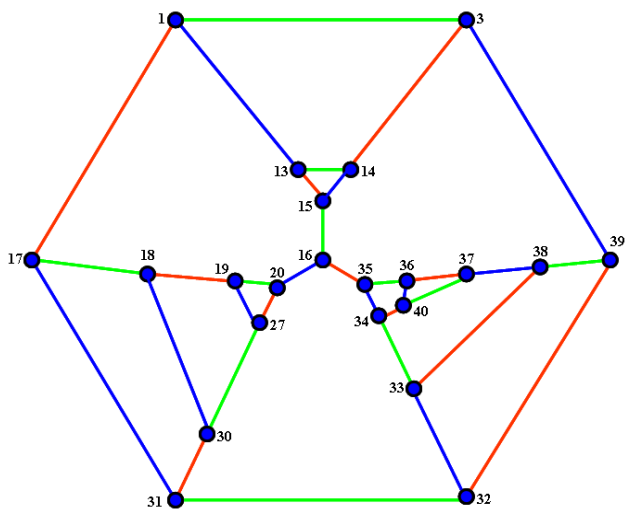

Рис. 3.32. Раскраска ребра ($v_{37}, v_{40}$).

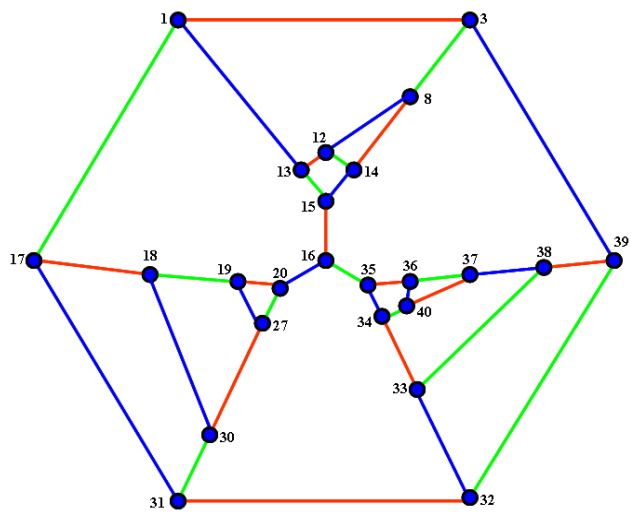

Рис. 3.33. Раскраска ребра ($v_8, v_{12}$).

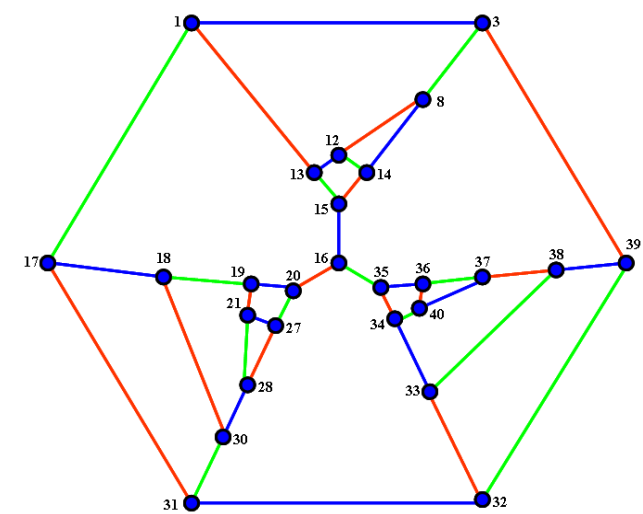

Рис. 3.34. Раскраска ребра ($v_{21}, v_{28}$).

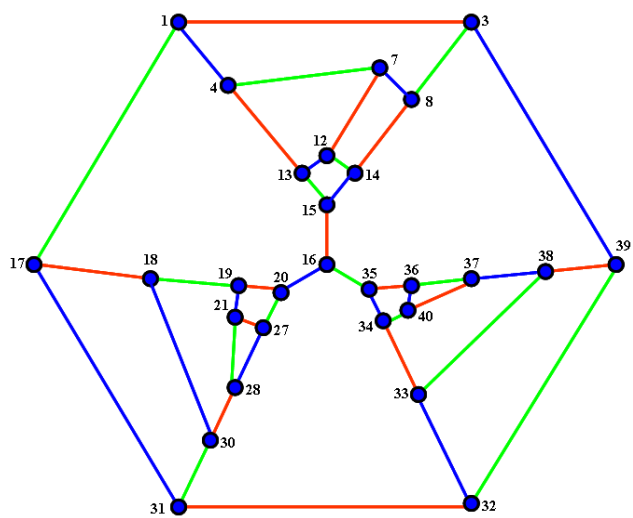

Рис. 3.35. Раскраска ребра ($v_4, v_7$).

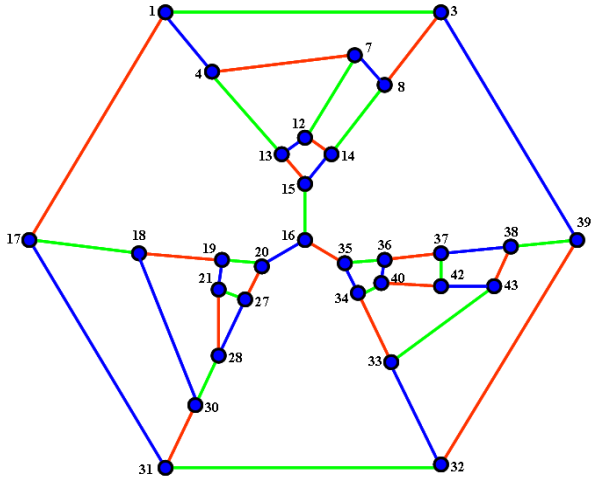

Рис. 3.36. Раскраска ребра ($v_{42}, v_{43}$).

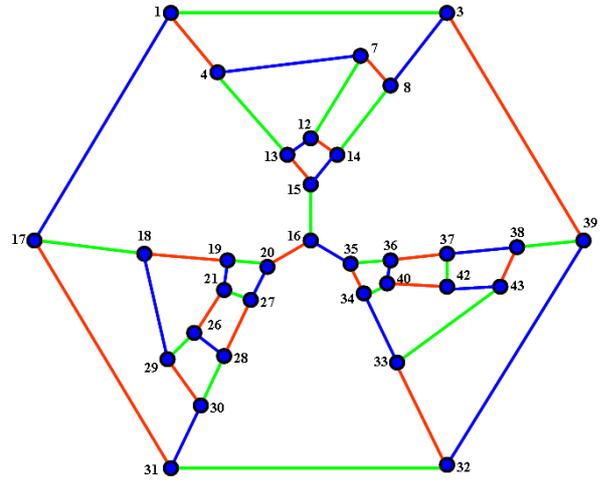

Рис. 3.37. Раскраска ребра ($v_{26}, v_{29}$).

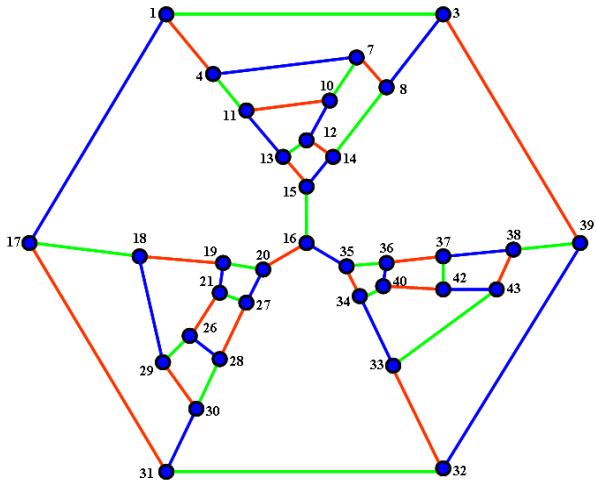

Рис. 3.38. Раскраска ребра ($v_{10}, v_{11}$).

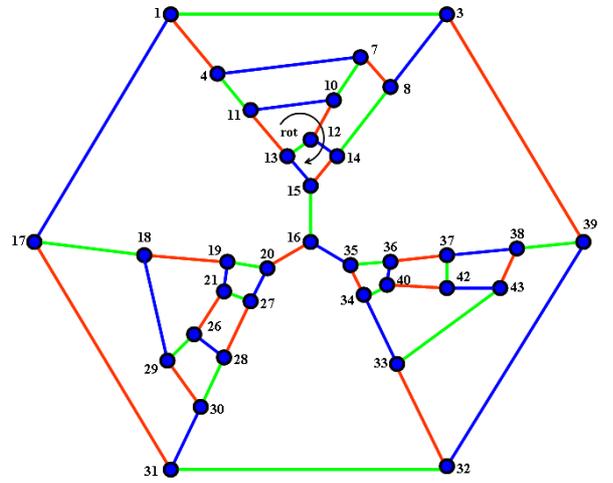

Рис. 3.39. Ротация диска ($v_{11}, v_{10}, v_{12}, v_{14}, v_{15}, v_{13}$).

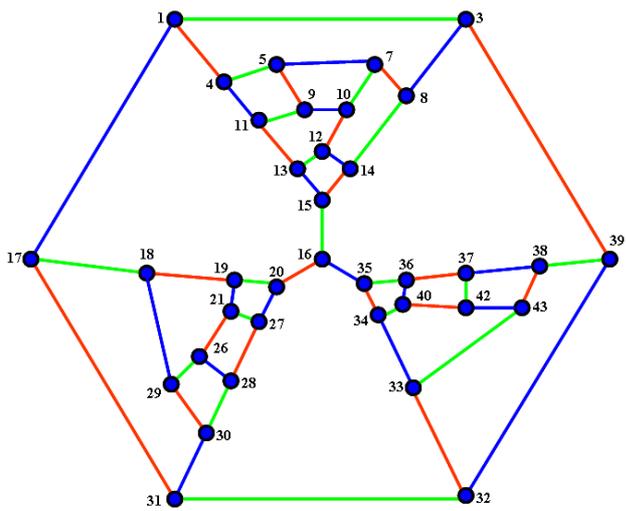

Рис. 3.40. Раскраска ребра ($v_5, v_9$).

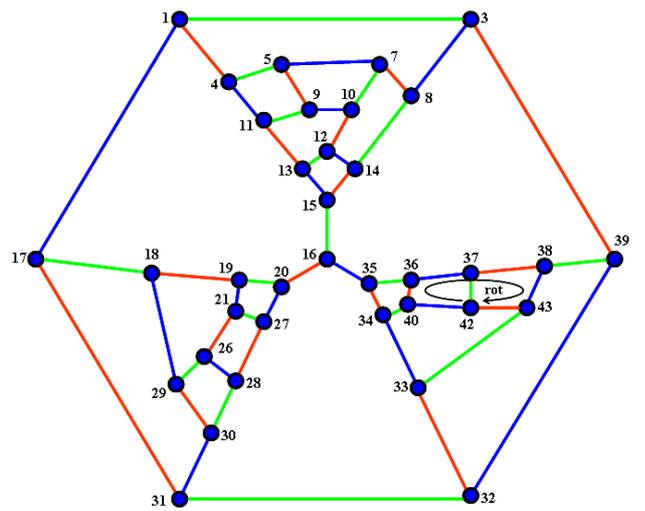

Рис. 3.41. Ротация диска ($v_{36}, v_{37}, v_{38}, v_{43}, v_{42}, v_{40}$).

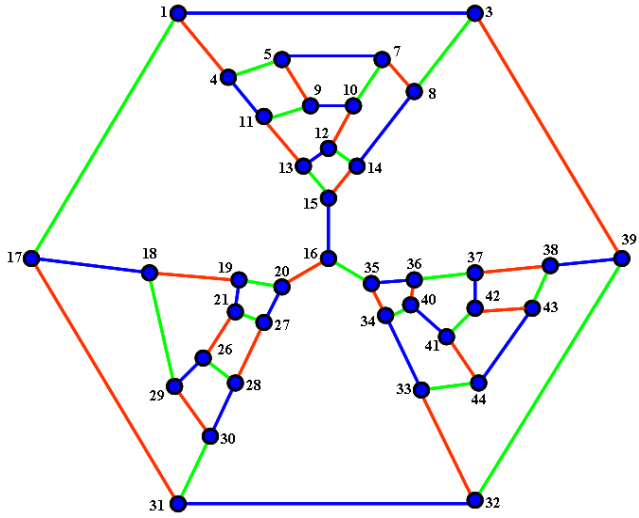

Рис. 3.42. Раскраска ребра ($v_{41}, v_{44}$).

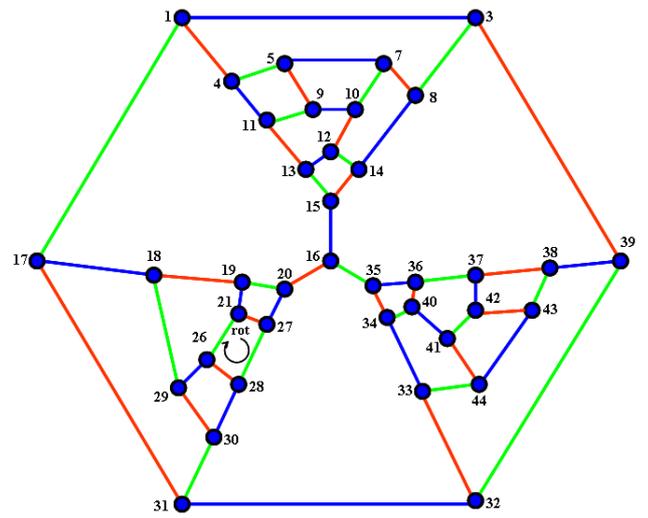

Рис. 3.43. Ротация диска ($v_{21}, v_{27}, v_{28}, v_{26}$).

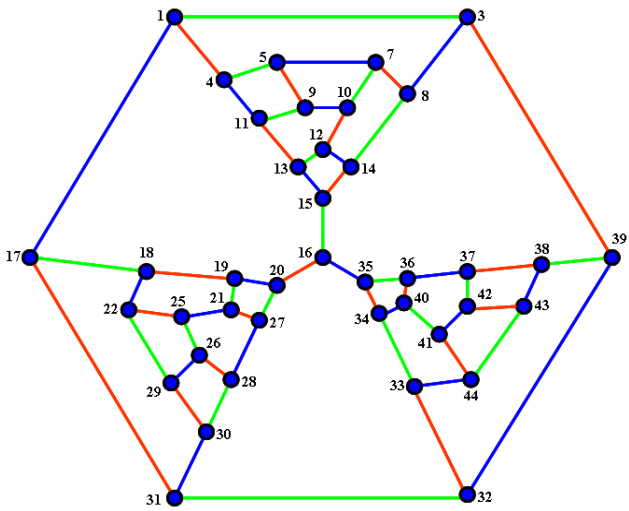

Рис. 3.44. Раскраска ребра ($v_{22}, v_{25}$).

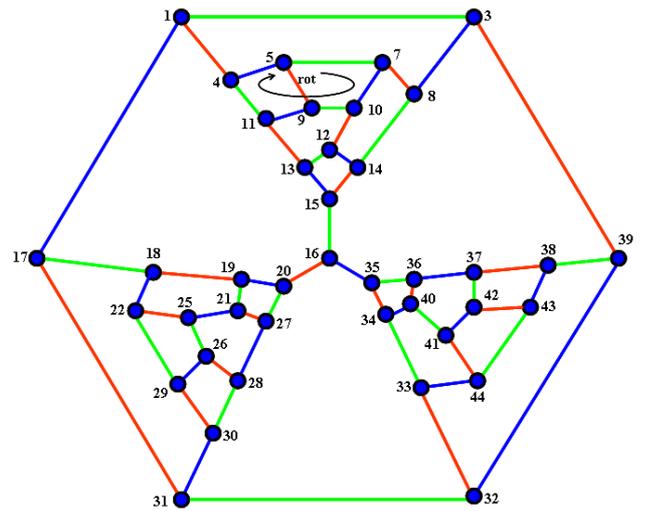

Рис. 3.45. Ротация диска ($v_4, v_5, v_7, v_{10}, v_9, v_{11}$).

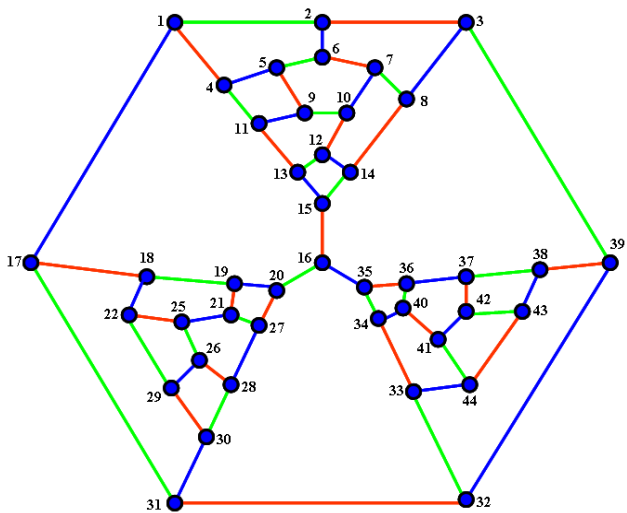

Рис. 3.46. Раскраска ребра ($v_2, v_6$).

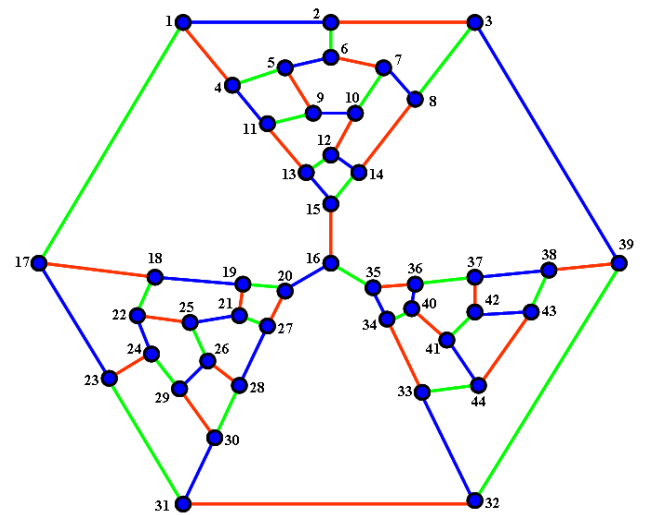

Рис. 3.47. Раскраска ребра ($v_{23}, v_{24}$).

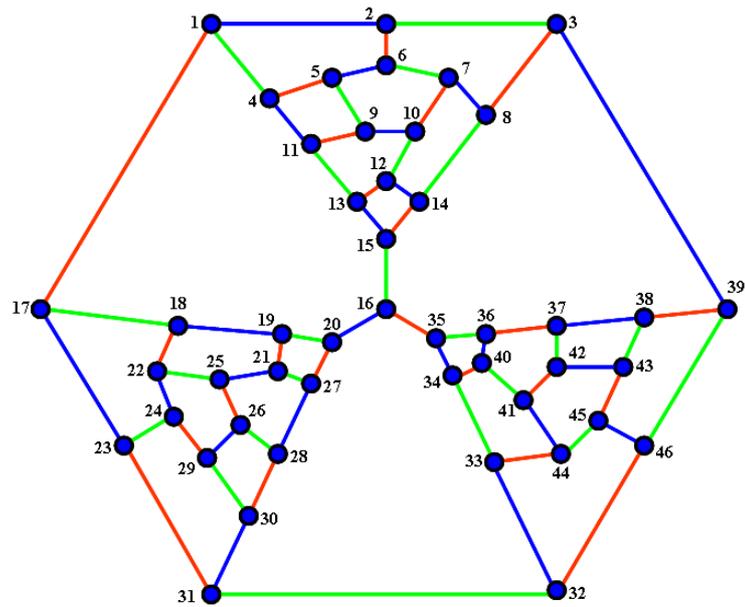

Рис. 3.48. Раскраска ребра ($v_{45}, v_{46}$) и окончательная раскраска рёбер графа Татта.

**Пример 3.2**. Рассмотрим раскраску рёбер для плоского кубического негамильтонового графа $G_{44}$ (см. рис. 3.49).

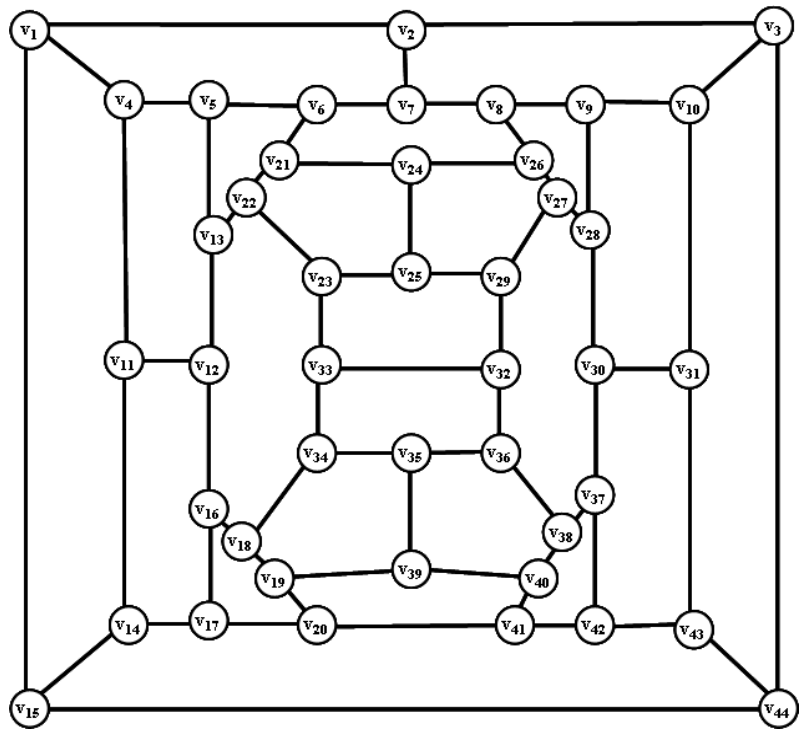

Рис. 3.49. Плоский кубический негамильтонов граф $G_{44}$.

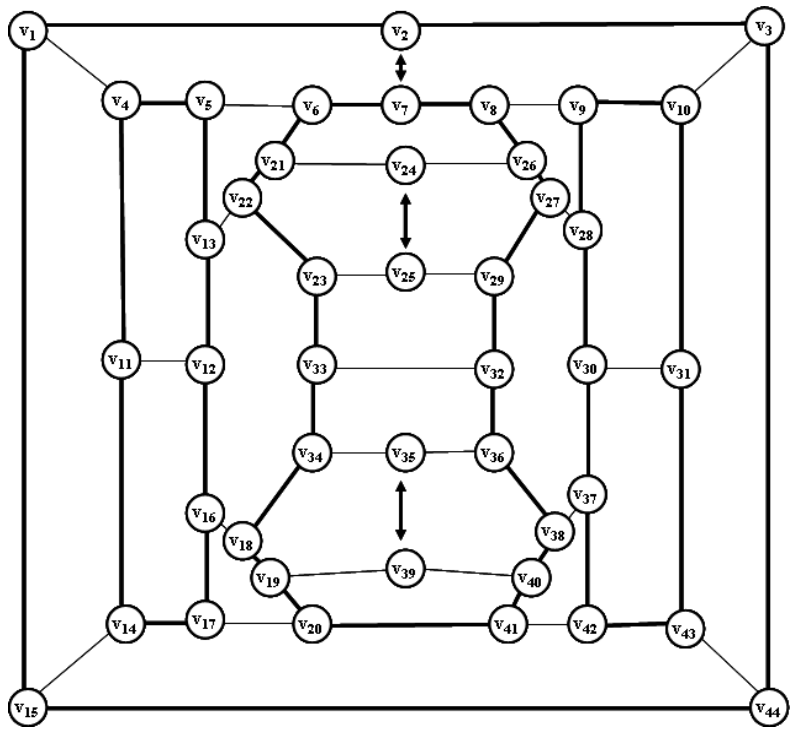

Рис. 3.50. Удаление рёбер в плоском кубическом негамильтоновом графе $G_{44}$.

С целью получения совокупности циклов чётной длины (дисков) удаляем из графа рёбра $(v_2,v_7),(v_{24},v_{25}),(v_{35},v_{39})$.

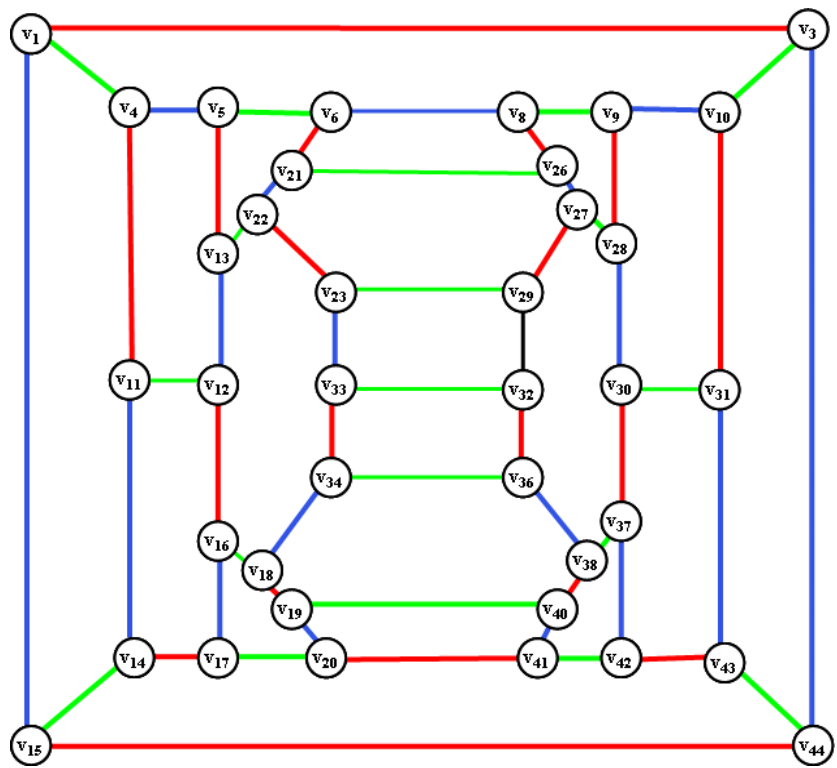

Рис. 3.51. Раскраска рёбер в выделенных дисках (пусть чётные диски будут зелёные).

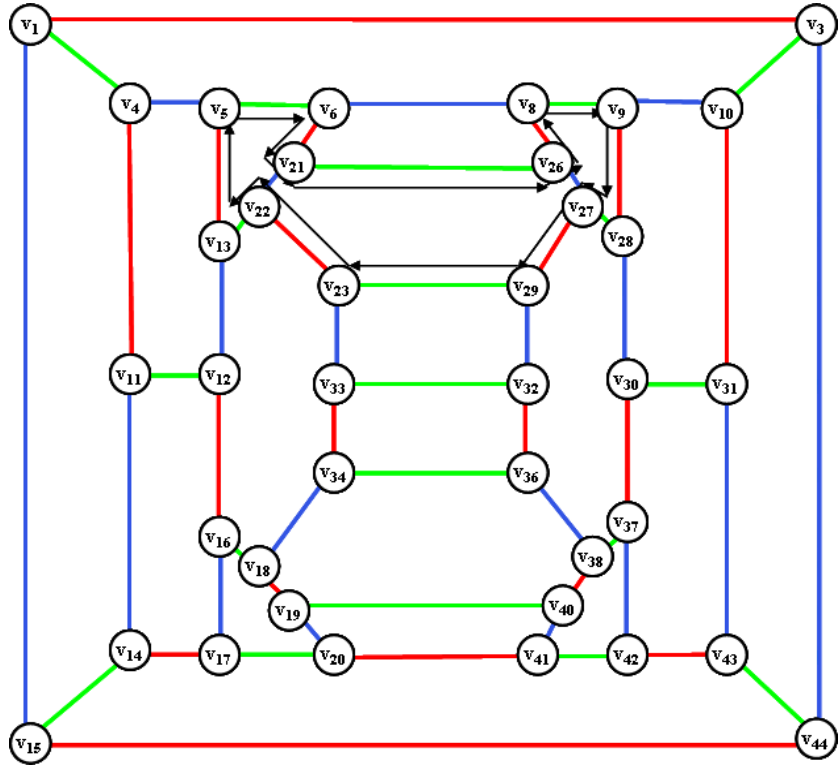

Рис. 3.52. Определение цветного диска для введения ребра ($v_{24}, v_{25}$).

Визуально, находим цветной диск для введения ребра ($v_{24}, v_{25}$). Это диск $<v_5, v_{13}, v_{22}, v_{23}, v_{29}, v_{27}, v_{28}, v_9, v_8, v_{26}, v_{21}, v_6>$.

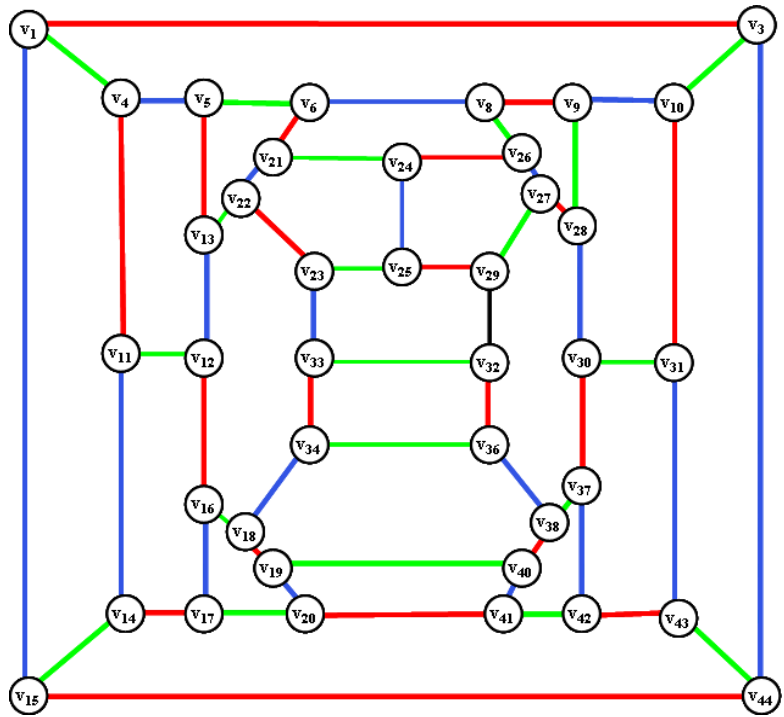

Рис. 3.53. Введение и раскраска ребра ($v_{24}, v_{25}$).

Введение и раскраска ребра ($v_{24}, v_{25}$) представлена на рис. 3.53. Осуществляется перекраска ребер цветного диска $<v_5, v_{13}, v_{22}, v_{23}, v_{25}, v_{29}, v_{27}, v_{28}, v_9, v_8, v_{26}, v_{24}, v_{21}, v_6>$.

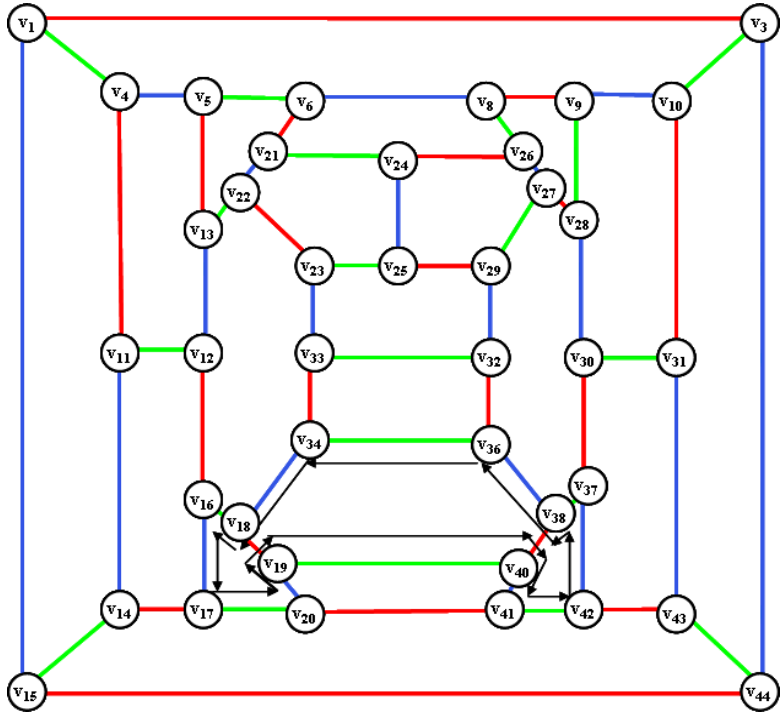

Рис. 3.54. Определение цветного диска для введения ребра ($v_{35},v_{39}$).

Визуально, определяем диск <$v_{16},v_{17},v_{20},v_{19},v_{40},v_{41},v_{42},v_{37},v_{38},v_{36},v_{34},v_{18}$> (см. рис. 3.54) для введения ребра ($v_{35},v_{39}$).

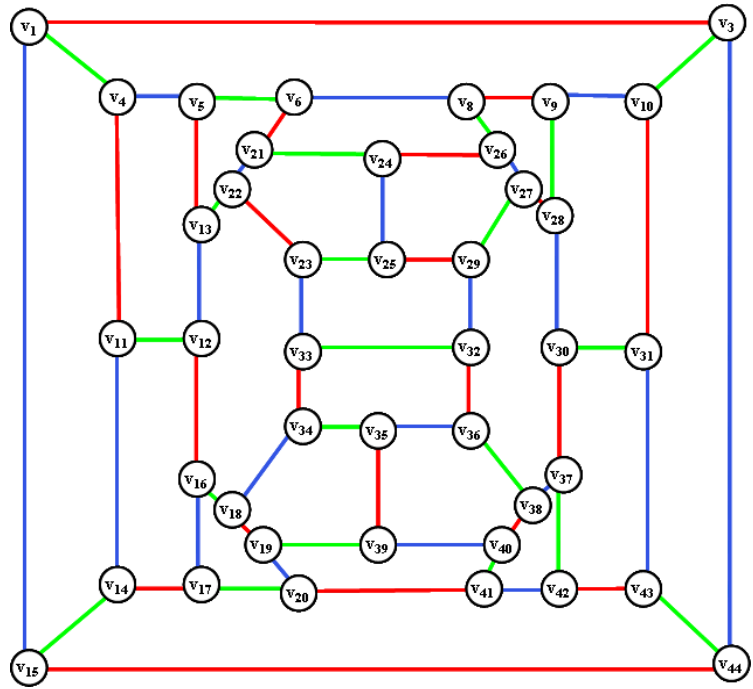

Рис. 3.55. Введение и раскраска ребра ($v_{35},v_{39}$).

Вводим и раскрашиваем ребро ($v_{35},v_{39}$), перекрашиваем ребра в цветном диске (см. рис. 3.55) <$v_{16},v_{17},v_{20},v_{19},v_{39},v_{40},v_{41},v_{42},v_{37},v_{38},v_{36},v_{35},v_{34},v_{18}$>.

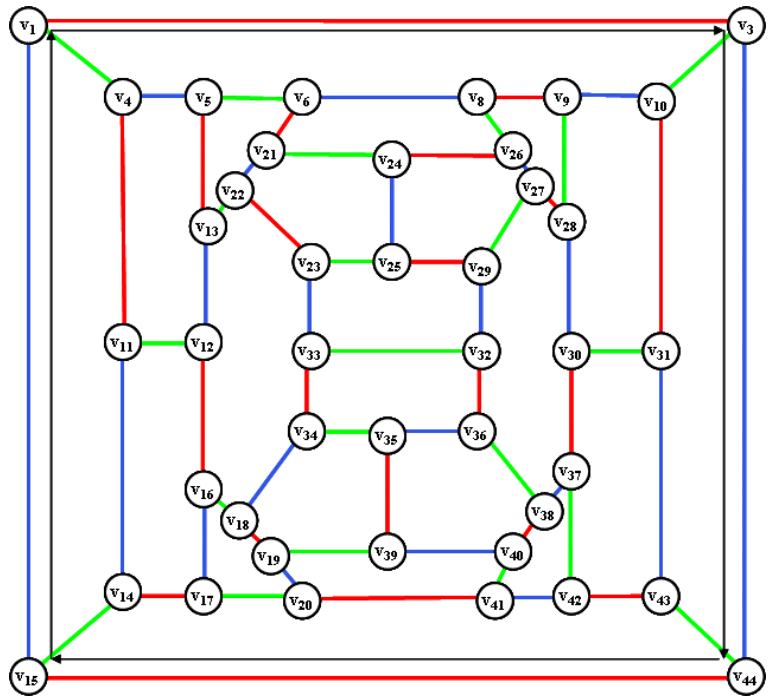

Рис. 3.56. Определение ротации цветного диска для введения ребра ($v_2,v_7$).

Осуществляем ротацию цветного диска <$v_{15},v_1,v_3,v_{44}$> с целью введения ребра ($v_2,v_7$).

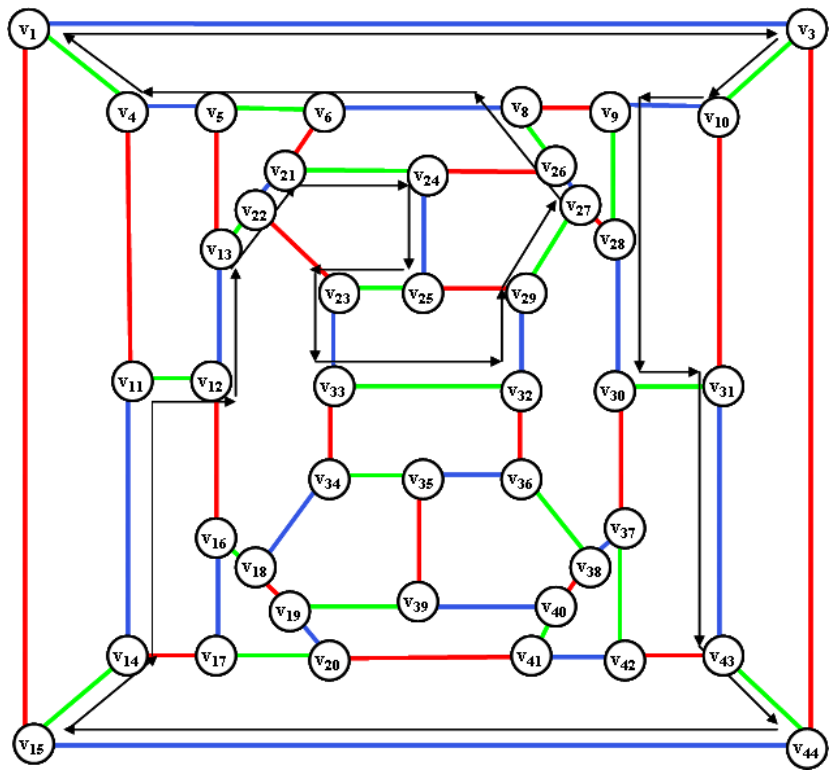

Рис. 3.57. Нахождение цветного диска для введения ребра ($v_2,v_7$) после ротации.

Находим цветной диск

<$v_1,v_3,v_{10},v_9,v_{28},v_{30},v_{31},v_{43},v_{44},v_{15},v_{14},v_{11},v_{12},v_{13},v_{22},v_{21},v_{24},v_{25},v_{23},v_{33},v_{32},v_{29},v_{27},v_{26},v_8,v_6,v_5,v_4$> для введения ребра ($v_2,v_7$). Окончательная раскраска негамильтонового графа $G_{44}$ представлена на рис. 3.58.

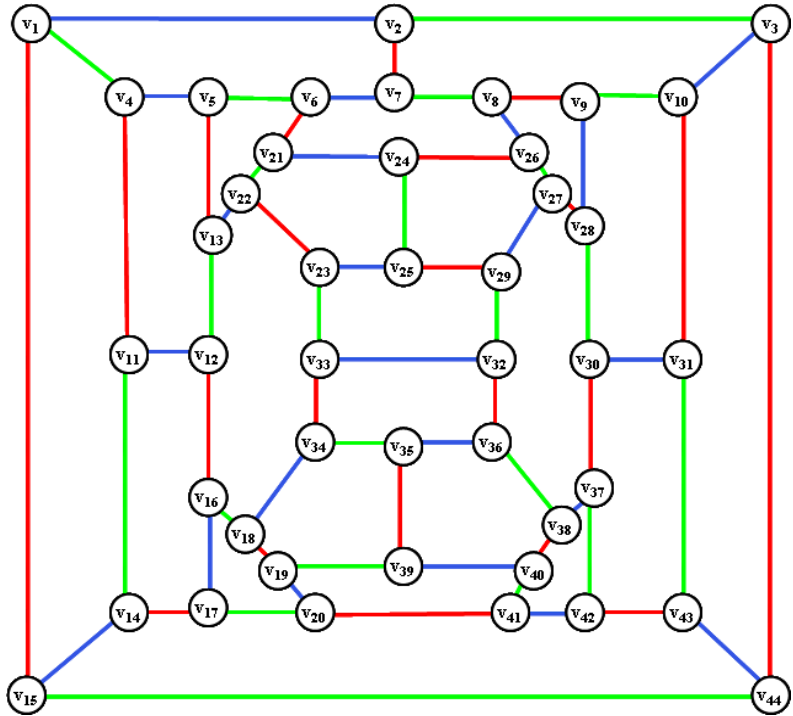

Рис. 3.58. Окончательная раскраска негамильтонового графа $G_{44}$.

**Пример 3.3.** Произведем раскраску ребер в негамильтоновом графе $G_{46}$.

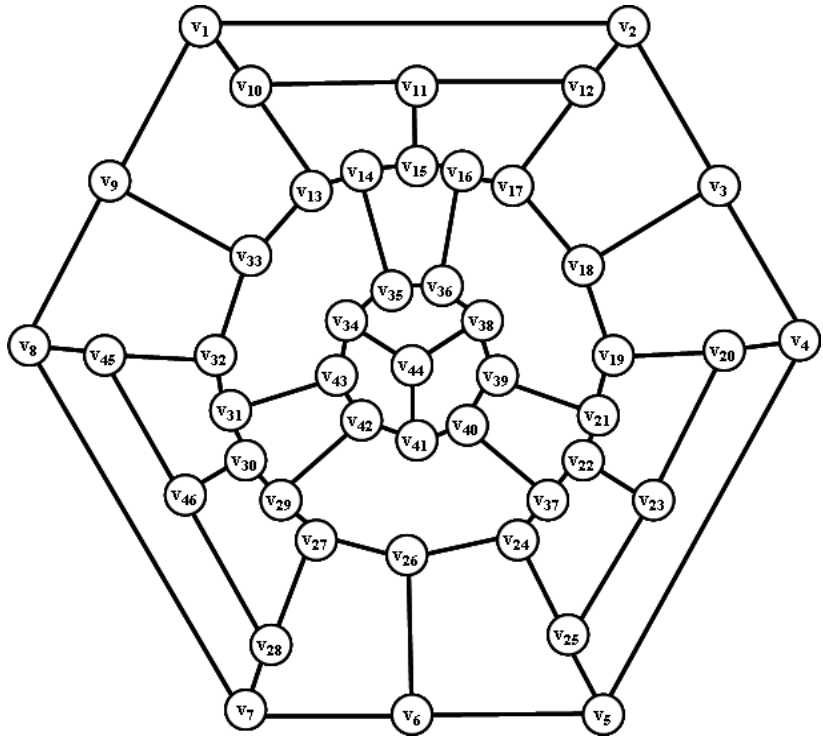

Рис. 3.59. Кубический негамильтонов граф $G_{46}$.

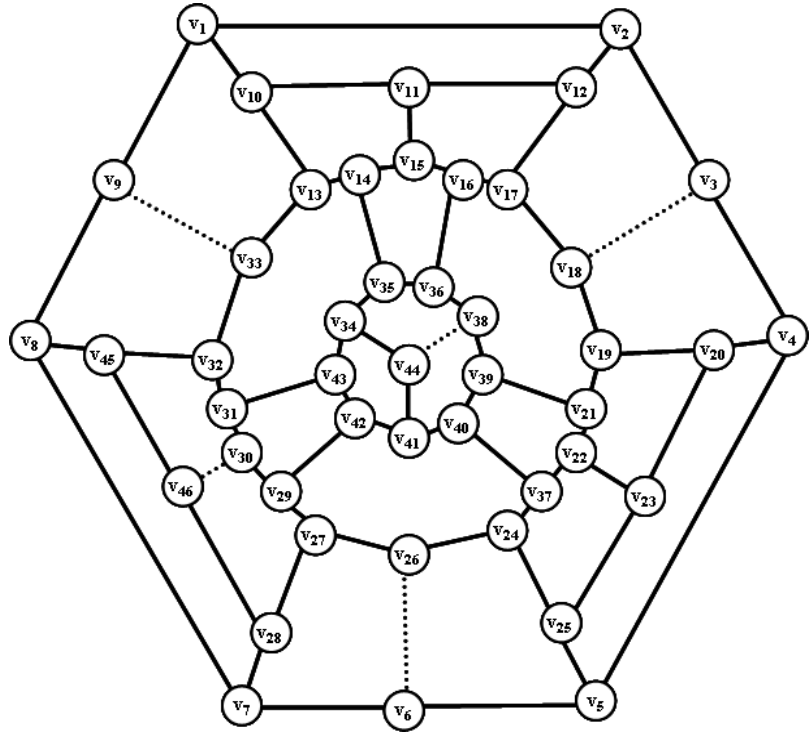

Рис. 3.60. Удаление рёбер с целью получения дисков чётной длины.

Удаляем из графа рёбра $(v_9,v_{33})$, $(v_3,v_{18})$, $(v_{30},v_{46})$, $(v_6,v_{26})$, $(v_{38},v_{44})$ с целью получения дисков чётной длины.

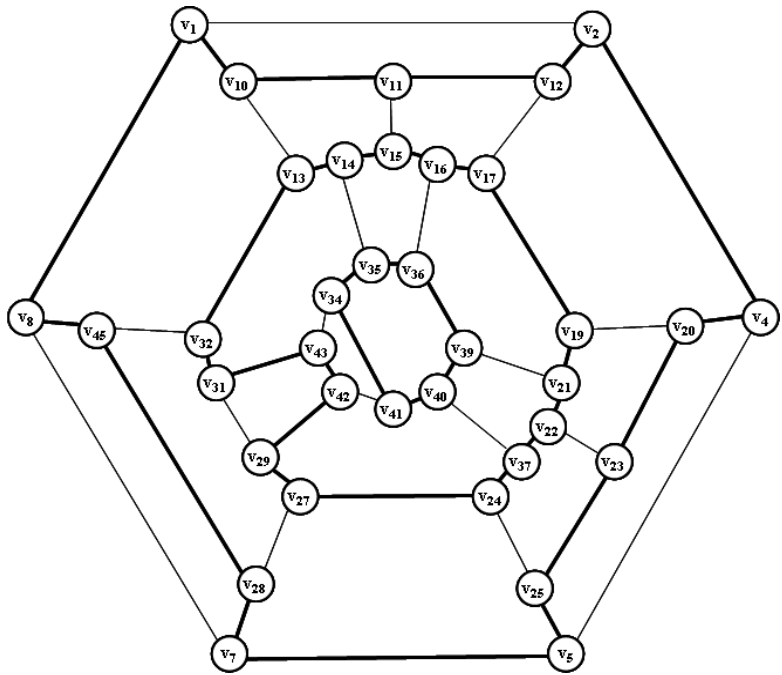

Рис. 3.61. Диски чётной длины.

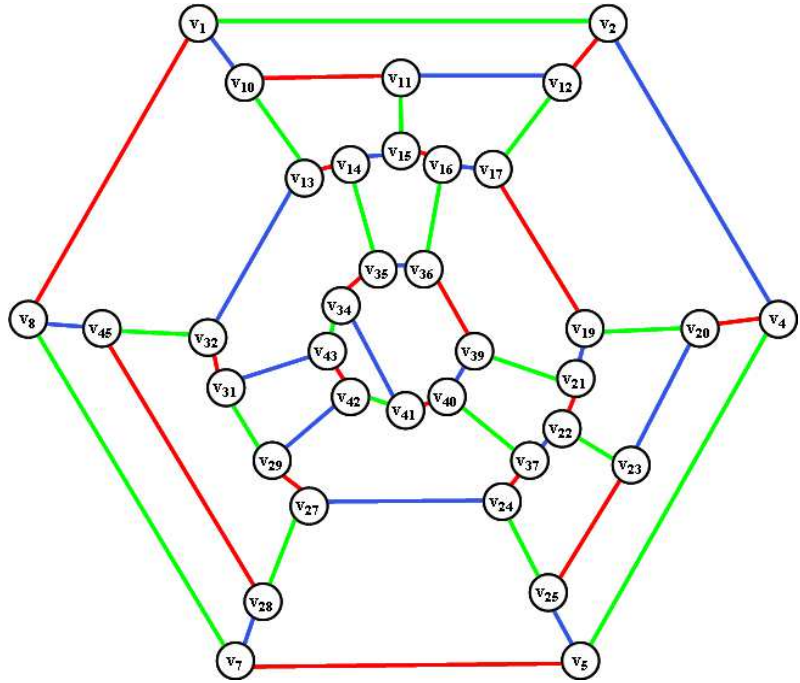

Рис. 3.62. Раскраска дисков четной длины.

Раскрашиваем диски четной длины в зеленый цвет. Тогда ребра данного цветного диска раскрашиваются в красный и синий цвета. Ребра, принадлежащие зеленому 1-фактору графа, также раскрашиваем в зеленый цвет.

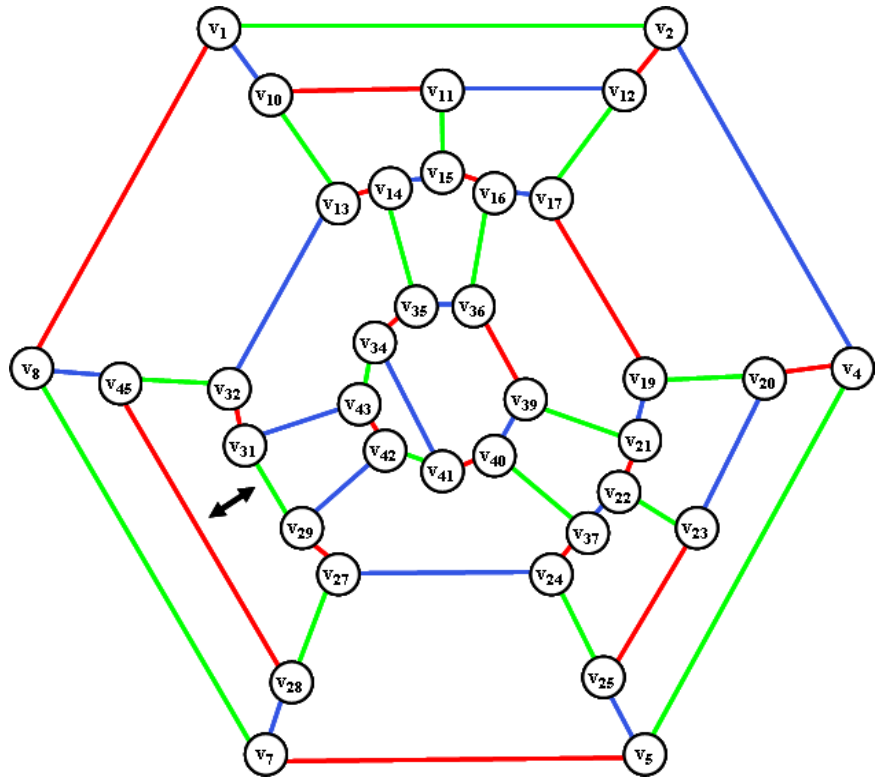

Рис. 3.63. Введение ребра ($v_{30}$, $v_{46}$).

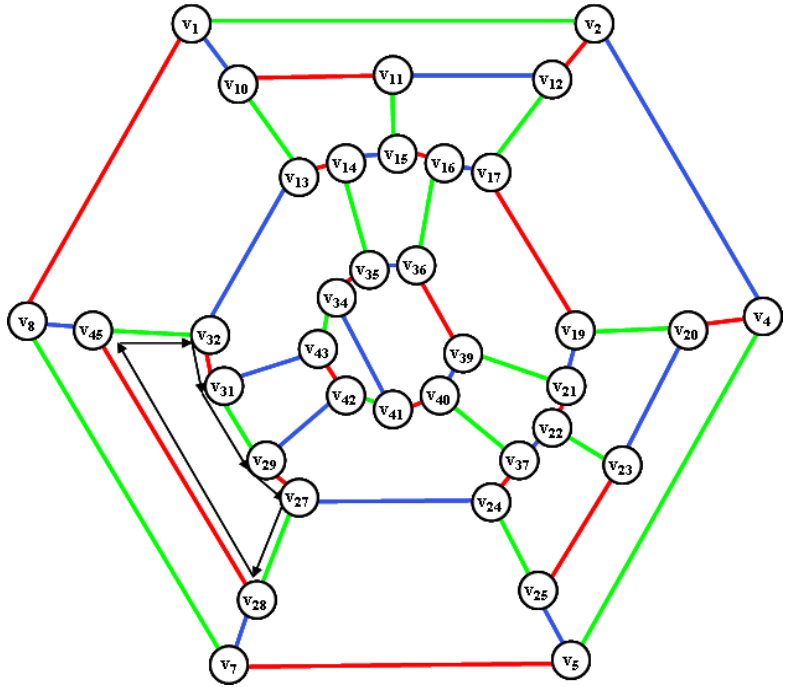

Рис. 3.64. Поиск диска для введения ребра ($v_{30}, v_{46}$).

Осуществляем поиск цветного диска для введения ребра ($v_{30}, v_{46}$). Это синий диск <$v_{45}, v_{32}, v_{31}, v_{29}, v_{27}, v_{28}$>.

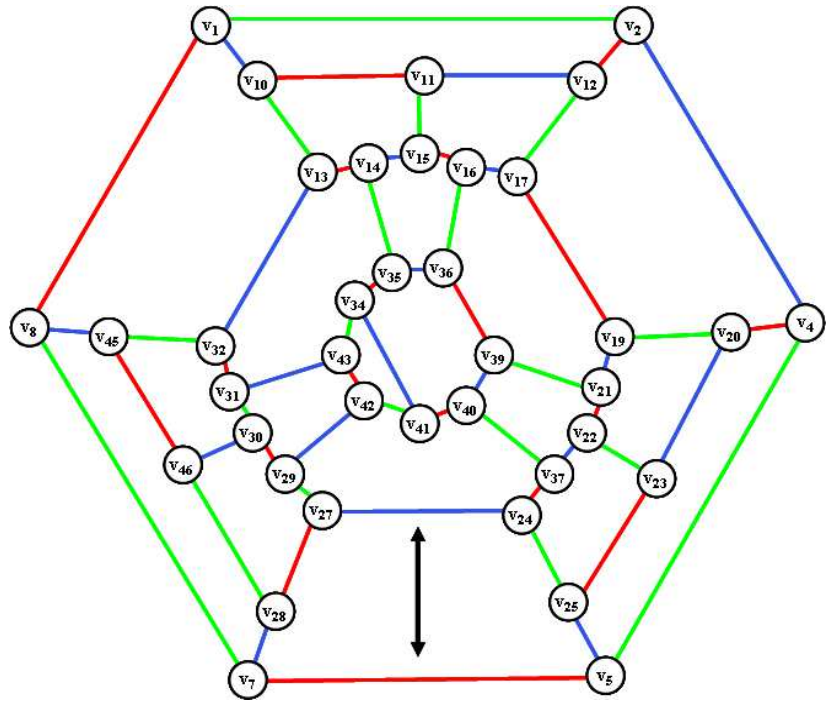

Рис. 3.65. Введение ребра ($v_{30}, v_{46}$) и перекраска ребер диска.

Вводим ребро ($v_{30}, v_{46}$) и осуществляем перекраску ребер диска

<$v_{45}, v_{32}, v_{31}, v_{30}, v_{29}, v_{27}, v_{28}, v_{46}$>. Определяем цветной диск

<$v_{27}, v_{28}, v_7, v_5, v_{25}, v_{23}, v_{20}, v_4, v_2, v_{12}, v_{11}, v_{10}, v_1, v_8, v_{45}, v_{46}, v_{30}, v_{29}, v_{42}, v_{43}, v_{31}, v_{32}, v_{13}, v_{14}, v_{15}, v_{16}, v_{17}, v_{19}, v_{21}, v_{22}, v_{37}, v_{24}$> с сцепленными ребрами для введения ребра ($v_6, v_{26}$))

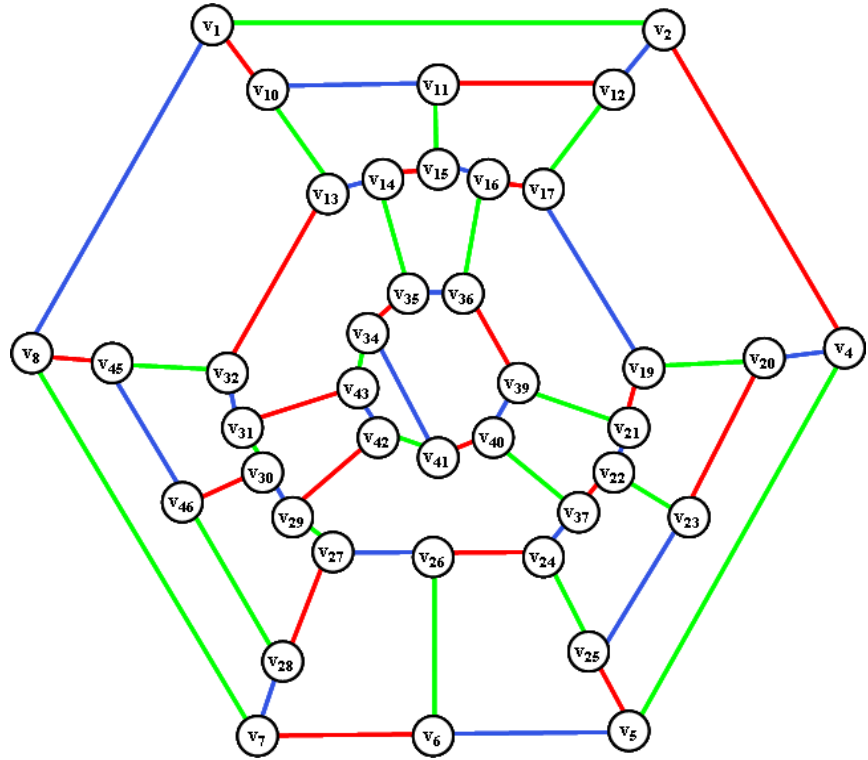

Рис. 3.66. Введение ребра ($v_6, v_{26}$) и перекраска ребер диска.

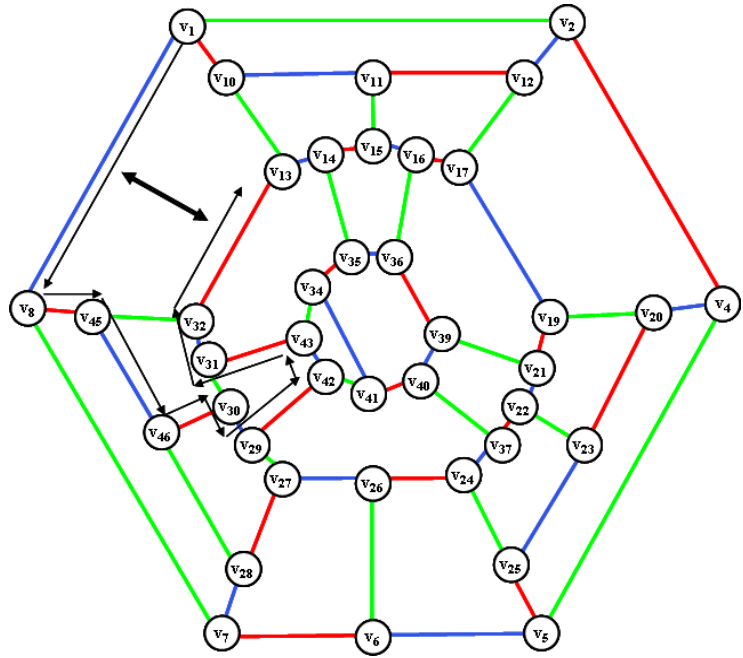

Рис. 3.67. Поиск цветного диска для введения ребра ($v_9, v_{33}$).

Осуществляем поиск диска для введения ребра ($v_9, v_{33}$). Это диск

<$v_{17}, v_{19}, v_{21}, v_{22}, v_{37}, v_{24}, v_{26}, v_{27}, v_{28}, v_7, v_6, v_5, v_{25}, v_{23}, v_{20}, v_4, v_2, v_{12}, v_{11}, v_{10}, v_1, v_9, v_8, v_{45}, v_{46}, v_{30}, v_{29}, v_{42}, v_{43}, v_{31}, v_{32}, v_{33}, v_{13}, v_{14}, v_{15}, v_{16}$>.

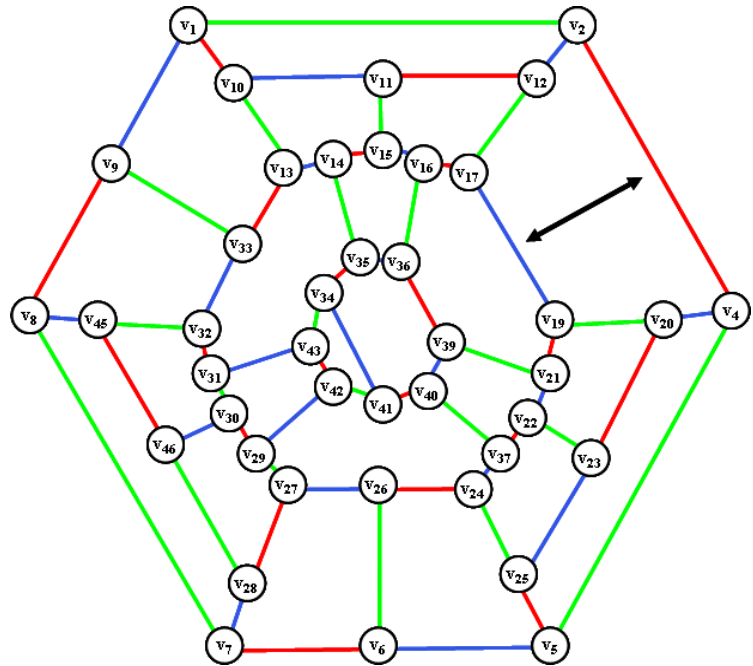

Рис. 3.68. Введение ребра ($v_9$,$v_{33}$) и перекраска ребер цветного диска.

Введение ребра ($v_9$,$v_{33}$) и перекраска ребер цветного диска. Определение диска для введения ребра ($v_3$,$v_{18}$).

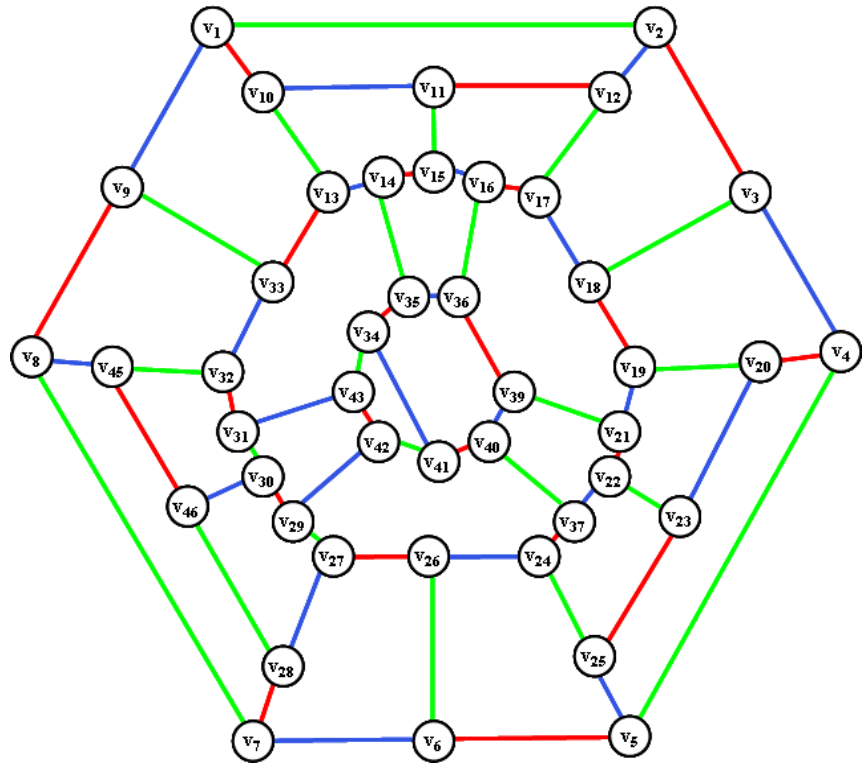

Рис. 3.69. Введение ребра ($v_3$,$v_{18}$) и перекраска ребер цветного диска.

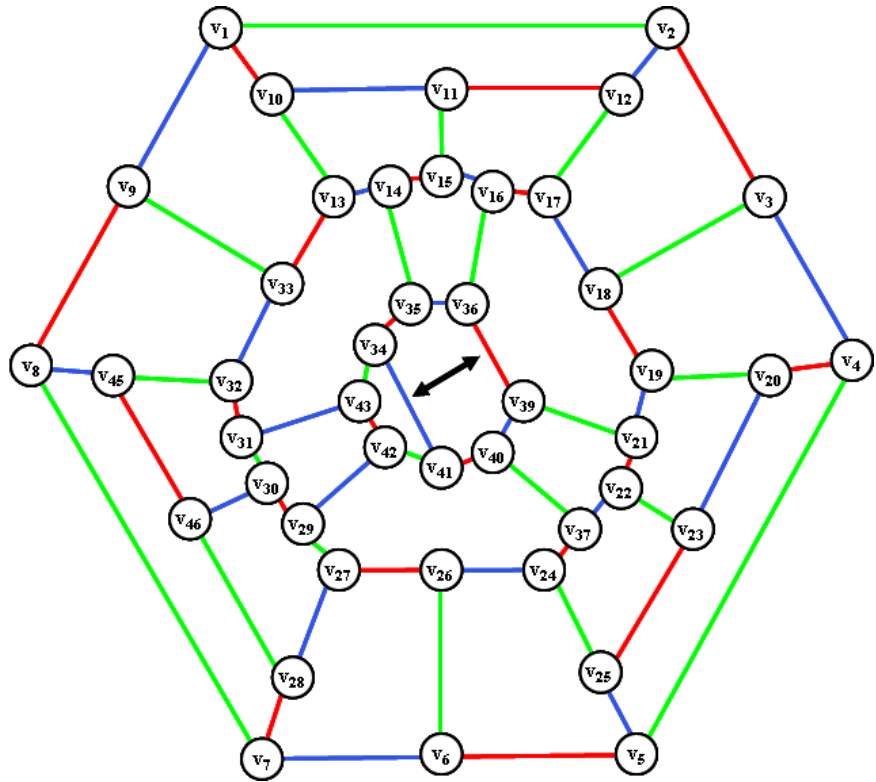

Рис. 3.70. Определение местоположения для введения ребра ($v_{38}, v_{44}$).

Нахождение цветного диска <$v_{34}, v_{35}, v_{36}, v_{39}, v_{40}, v_{41}$> для введения ребра ($v_{38}, v_{44}$).

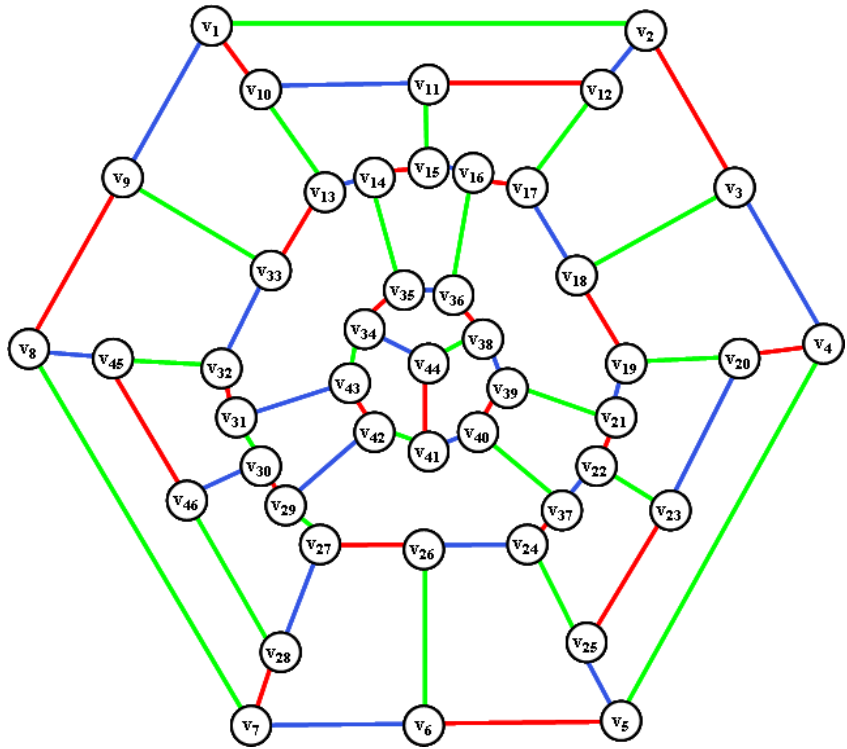

Рис. 3.71. Окончательная раскраска кубического негамильтонового графа $G_{46}$.

**Пример 3.4.** Раскрасим кубический негамильтонов граф $G_{42}$.

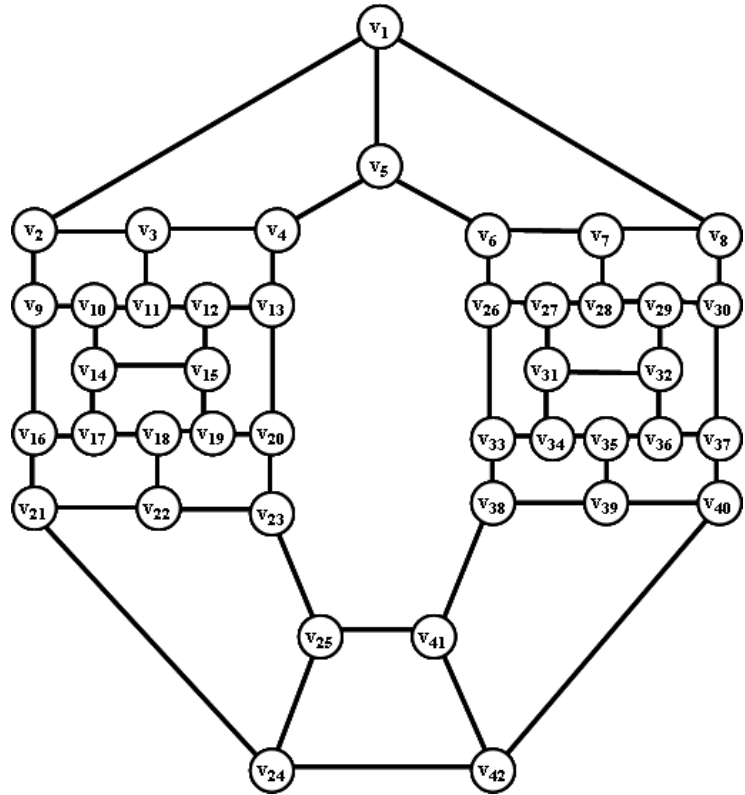

Рис. 3.72. Кубический негамильтонов граф $G_{42}$.

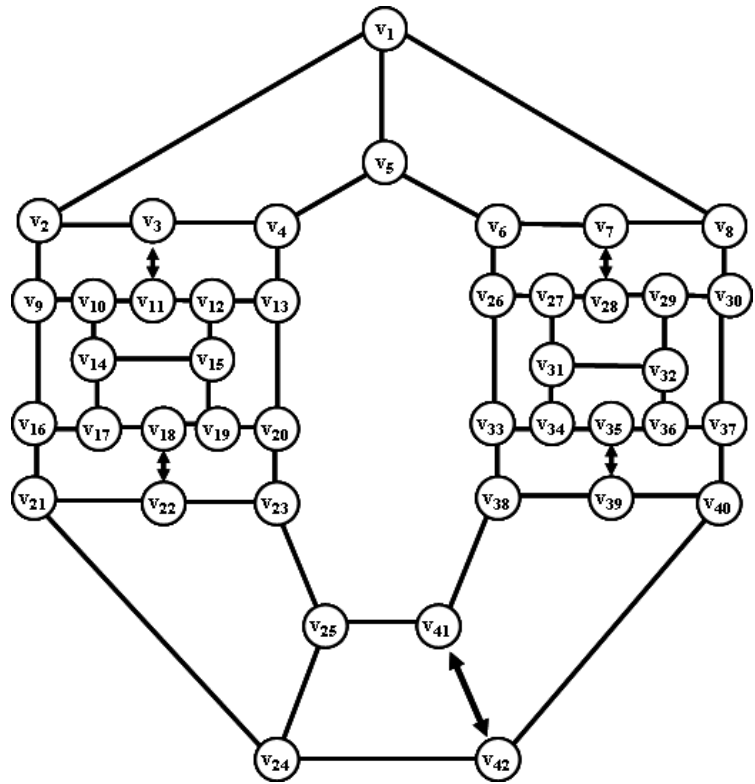

Рис. 3.73. Удаление ребер для получения дисков четной длины.

Удаляем из графа ребра $(v_3,v_{11}),(v_{18},v_{22}),(v_7,v_{28}),(v_{35},v_{39}),(v_{41},v_{42})$ с целью получения дисков четной длины.

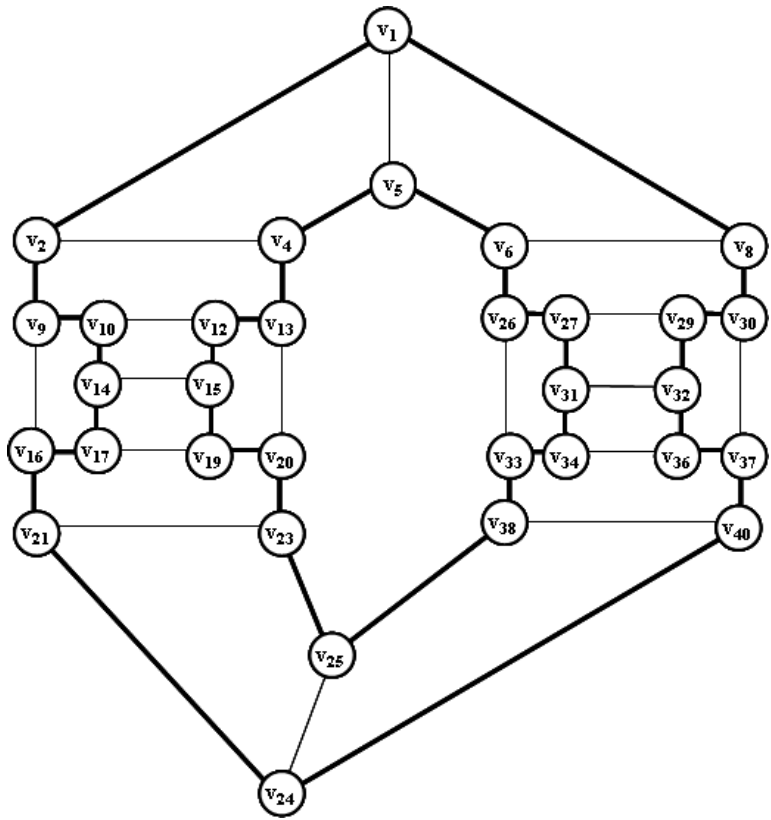

Рис. 3.74. Построение дисков четной длины.

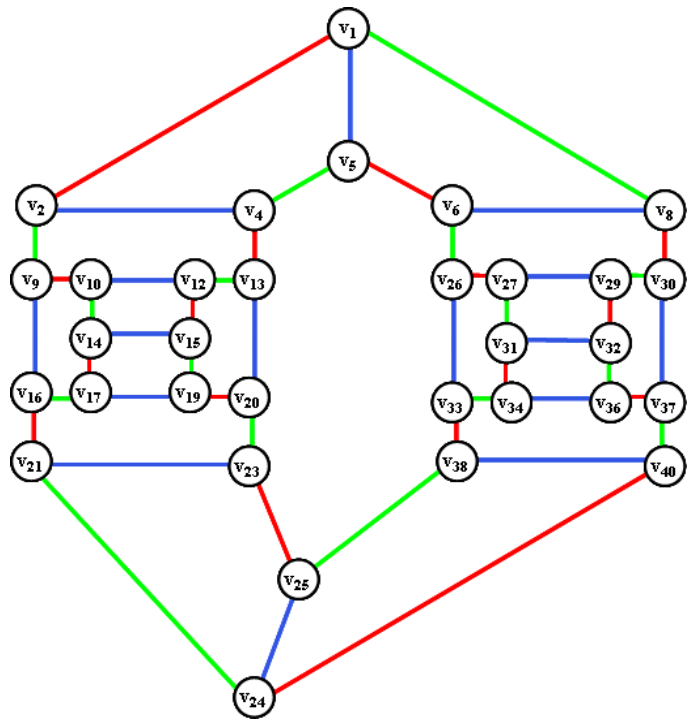

Рис. 3.75. Раскраска дисков четной длины (синий цвет).

Производим раскраску дисков четной длины в синий цвет и производим раскраска соответствующего 1-фактора также в синий цвет.

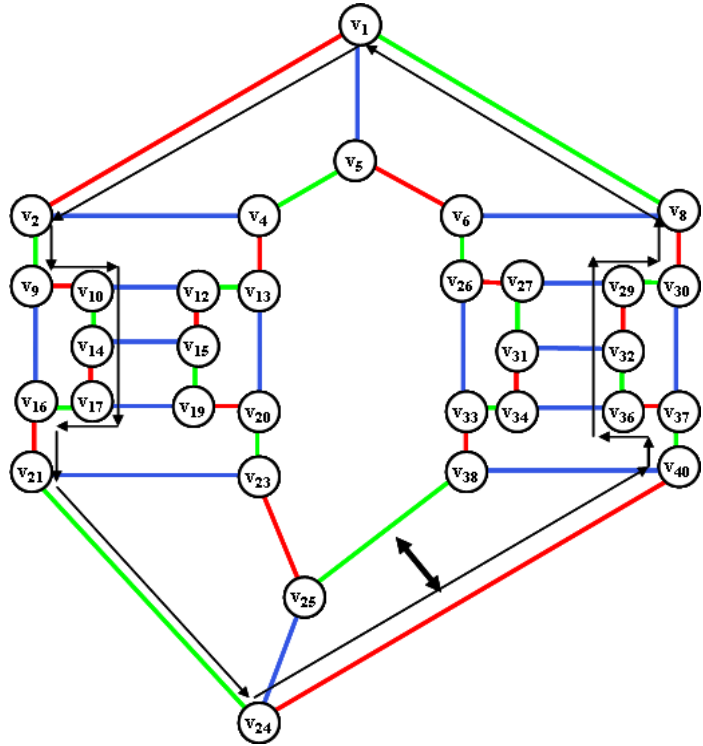

Рис. 3.76. Ротация цветного диска для введения ребра ($v_{41}$,$v_{42}$).

Осуществляем ротацию цветного диска <$v_{21}$,$v_{16}$,$v_{17}$,$v_{14}$,$v_{10}$,$v_{9}$,$v_{2}$,$v_{1}$,$v_{8}$,$v_{30}$,$v_{29}$,$v_{32}$,$v_{36}$,$v_{37}$,$v_{40}$,$v_{24}$> для введения ребра ($v_{41}$,$v_{42}$).

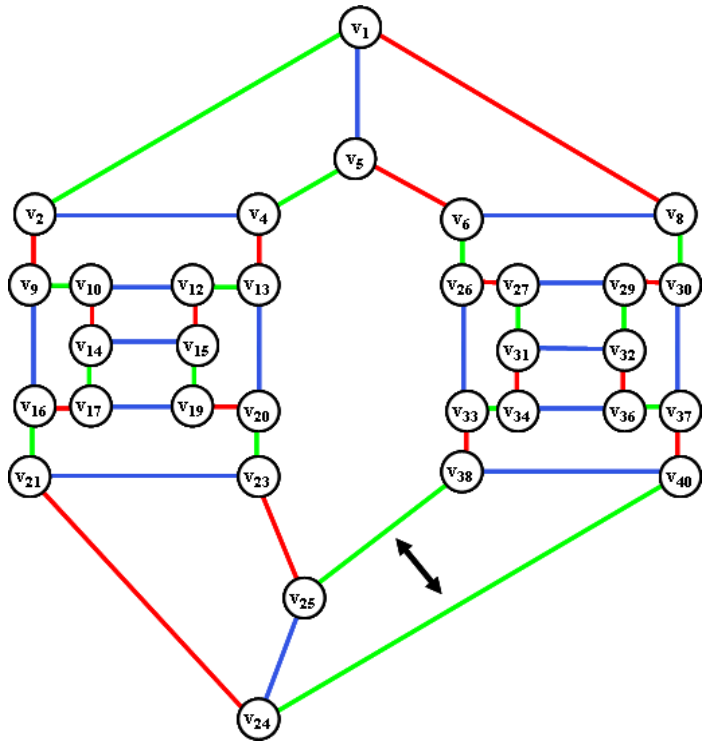

Рис. 3.77. Перекраска ребер в процессе ротации цветного диска.

После проведения ротации, определяем цветной диск <$v_{25}$,$v_{38}$,$v_{40}$,$v_{24}$> для проведения ребра ($v_{41}$,$v_{42}$).

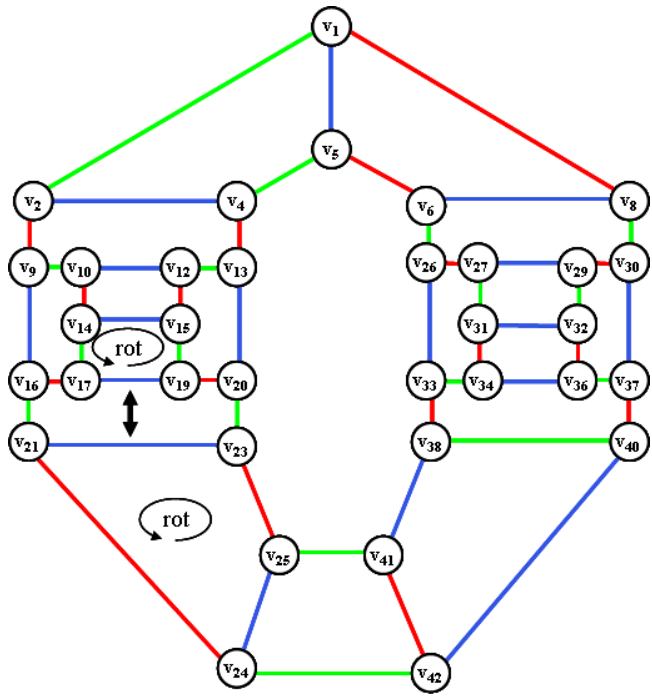

Рис. 3.78. Введение ребра ($v_{41}$,$v_{42}$) и перекраска ребер цветного диска.

После введения ребра ($v_{41}$,$v_{42}$) осуществляем ротацию цветных дисков <$v_{21}$,$v_{23}$,$v_{25}$,$v_{24}$> и <$v_{14}$,$v_{15}$,$v_{19}$,$v_{17}$> для введения ребра ($v_{18}$,$v_{22}$).

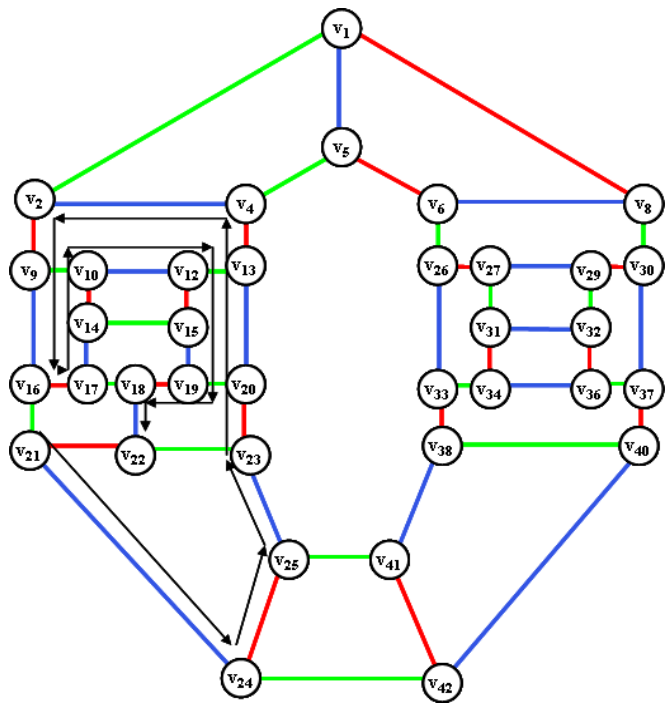

Рис. 3.79. Раскраска ребер цветного диска после введения ребра ($v_{18}$,$v_{22}$).

После введения ребра ($v_{18}$,$v_{22}$) перекрашиваем ребра цветного диска <$v_{16}$,$v_{17}$,$v_{18}$,$v_{19}$,$v_{20}$,$v_{23}$,$v_{22}$,$v_{21}$>. Далее производим ротацию диска <$v_{21}$,$v_{24}$,$v_{25}$,$v_{23}$,$v_{20}$,$v_{13}$,$v_{4}$,$v_{2}$,$v_{9}$,$v_{16}$,$v_{17}$,$v_{14}$,$v_{10}$,$v_{12}$,$v_{15}$,$v_{19}$,$v_{18}$,$v_{22}$> для введения ребра ($v_{3}$,$v_{11}$).

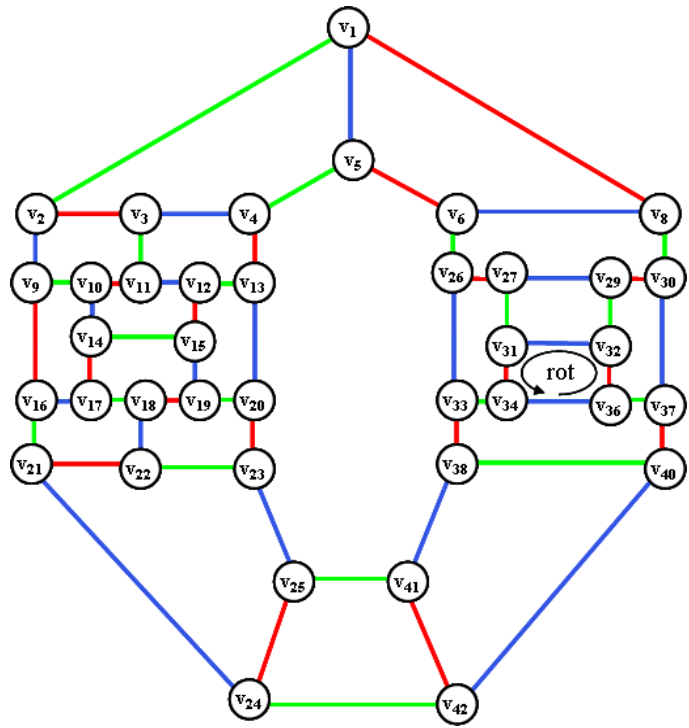

Рис. 3.80. Введение ребра ($v_3, v_{11}$) и перекраска ребер в цветном диске.

Для введения ребра ($v_{35}, v_{39}$) производим ротацию цветного диска <$v_{31}, v_{32}, v_{36}, v_{34}$>.

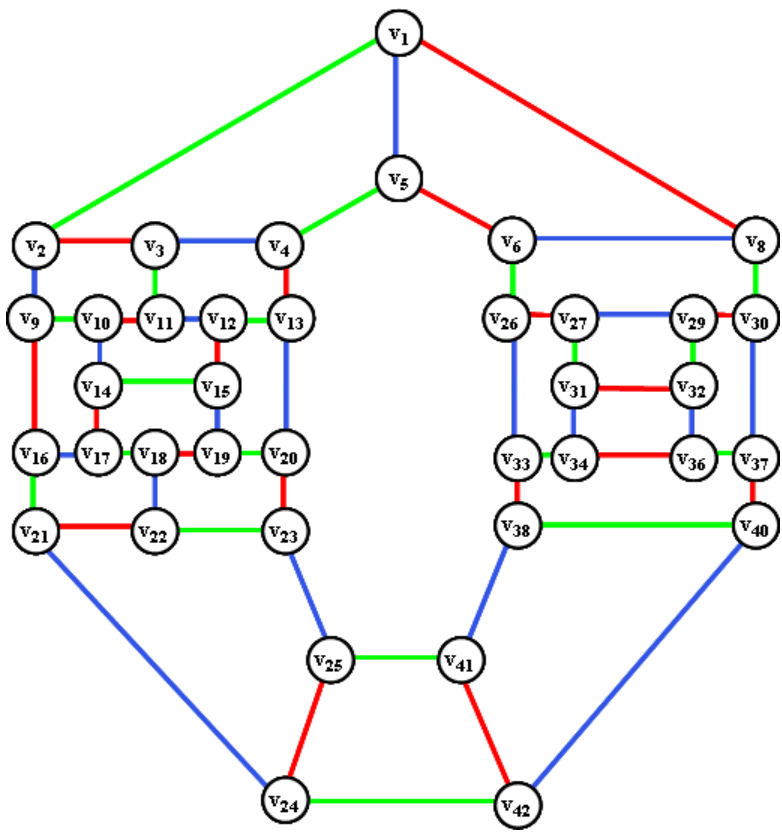

Рис. 3.81. Перекраска ребер после операции ротации цветного диска.

Производим перекраску ребер цветного диска <$v_{31}, v_{32}, v_{36}, v_{34}$>.

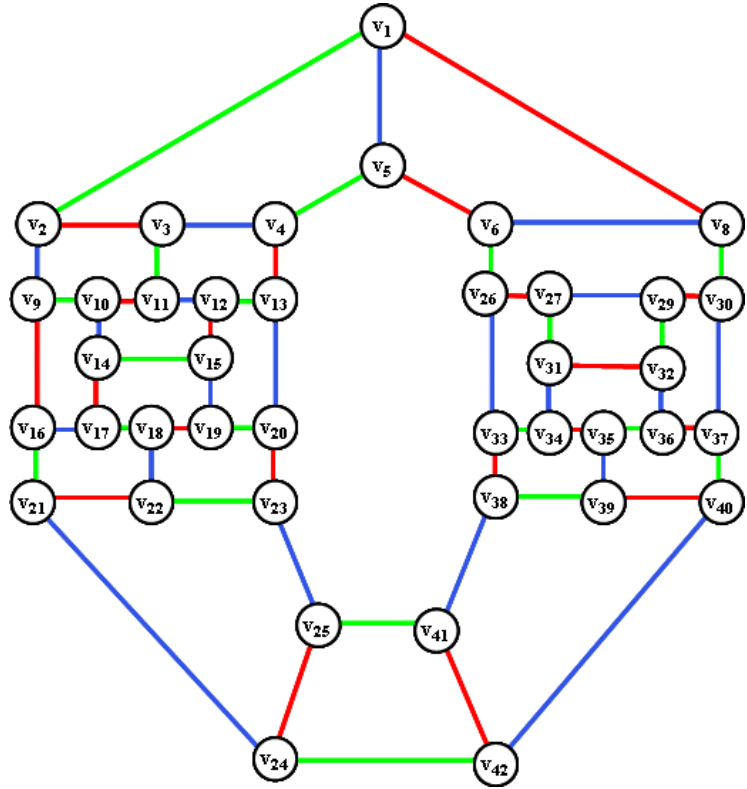

Рис. 3.82. Введение ребра ($v_{35}, v_{39}$) и перекраска ребер цветного диска.

Вводим ребро ($v_{35}, v_{39}$) и перекрашиваем ребра цветного диска $<v_{33}, v_{34}, v_{35}, v_{36}, v_{37}, v_{40}, v_{39}, v_{38}>$

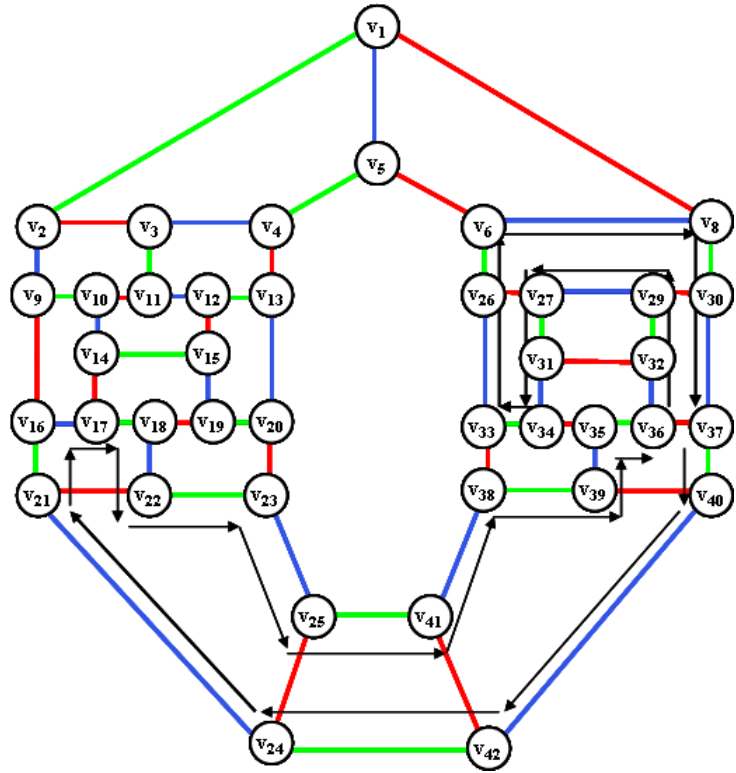

Рис. 3.83. Определение цветного диска для введения ребра ($v_7, v_{28}$).

Осуществляем ротацию диска

$<v_6, v_8, v_{30}, v_{37}, v_{40}, v_{42}, v_{24}, v_{21}, v_{16}, v_{17}, v_{18}, v_{22}, v_{23}, v_{25}, v_{41}, v_{38}, v_{39}, v_{35}, v_{36}, v_{32}, v_{29}, v_{27}, v_{31}, v_{34}, v_{33}, v_{26}>$ для вве-

дения ребра ($v_7, v_{28}$). Введение ребра ($v_7, v_{28}$) и окончательная раскраска негамильтонового графа $G_{42}$:

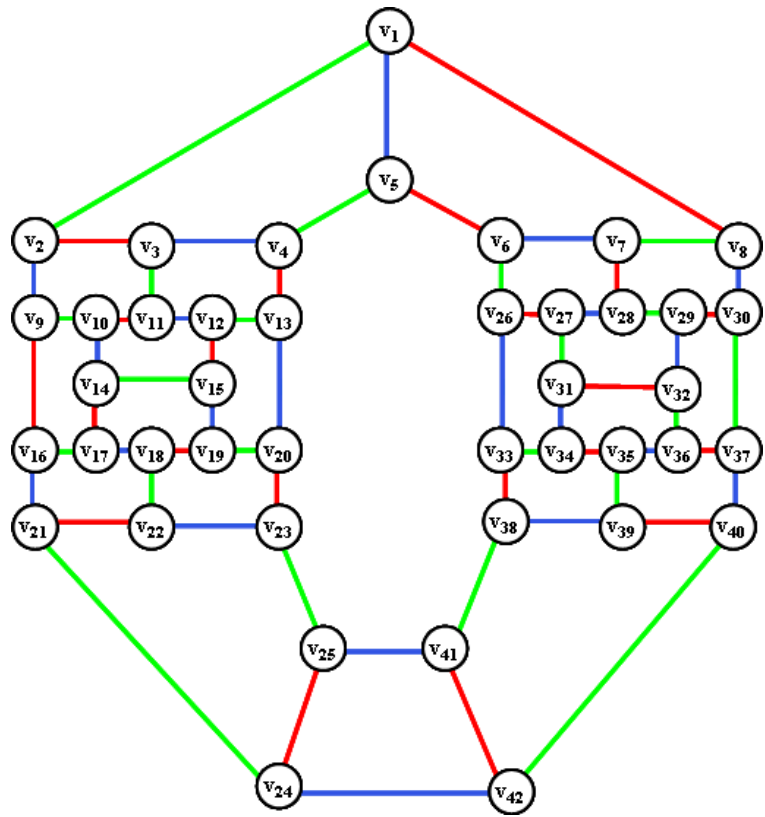

Рис. 3.84. Окончательная раскраска негамильтонового графа $G_{42}$.

Рассматривая приведенные примеры, заметим, что раскраска вновь введенных ребер определяется методами поиска цветных дисков для определения цвета введенного ребра. Данные методы явно носят искусственный характер. Требуется определенное искусство для нахождения таких цветных дисков.

| 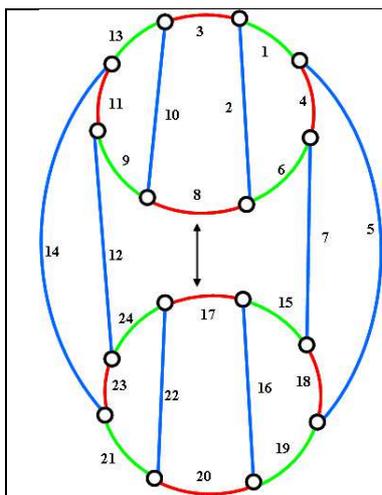 | 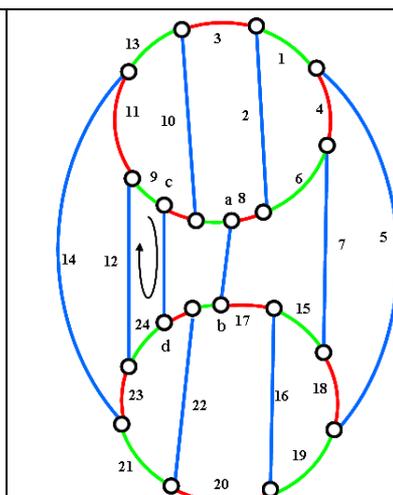 | 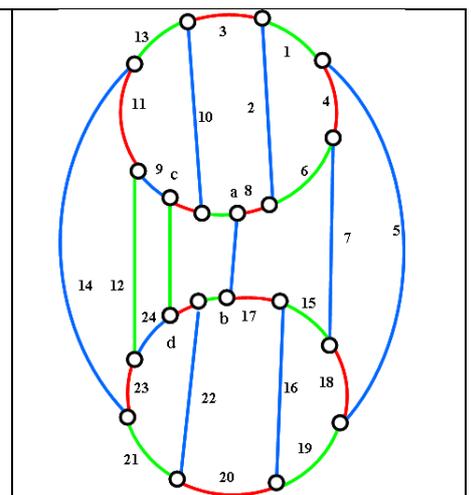 |
|---|---|---|
| Рис. 3.85. Введение ребра. | Рис. 3.86. Введение двух ребер. | Рис. 3.87. Перекраска квадрата. |

Рассмотрим вопрос создания стабильного алгоритма для раскраски, вновь введенного ребра.

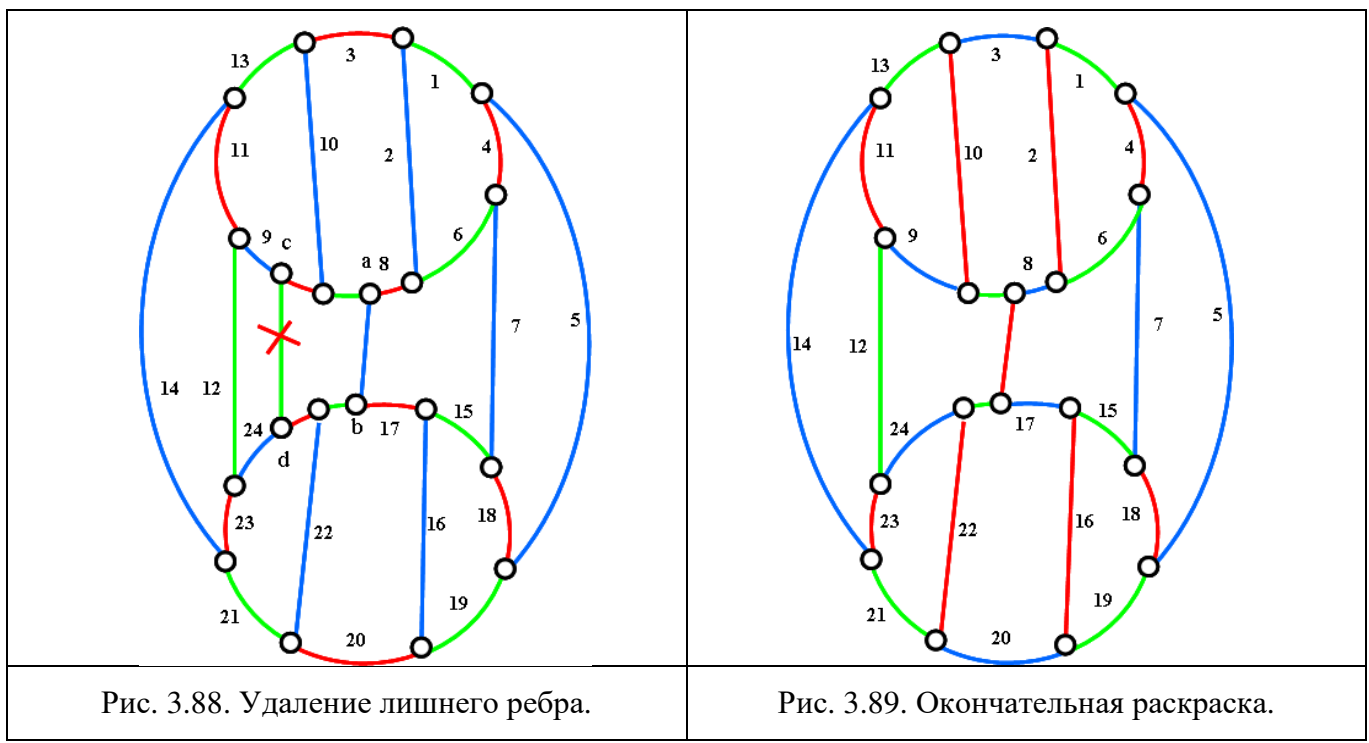

| Рис. 3.88. Удаление лишнего ребра. | Рис. 3.89. Окончательная раскраска. |

Предположим, что нужно вставить ребро (a,b) между ребрами $e_8$ и $e_{17}$. Предварительный цвет данного ребра определим как цвет дисков проходящих по этим сцепленным ребрам – синий (см. рис. 3.85). Будем вставлять не одно ребро, а два ребра (a,b) и (c,d), располагая их в соответствующим простом цикле (в нашем случае - это цикл $\{e_6, e_7, e_8, e_9, e_{12}, e_{15}, e_{17}, e_{24}\}$) где находится сцепленные ребра (см. рис. 3.86). Причем второе ребро должно располагаться так, что бы образовать одноцветный цикл длиной четыре с ребром цикла того же цвета (в нашем случае цикл $\{e_9, e_{12}, e_{24}, (c,d)\}$). Стоит заметить, что введение этих двух ребер не изменяет четность цветных дисков. Но тогда вращение (ротация) диска $\{e_9, e_{12}, e_{24}, (c,d)\}$ изменяет цвет второго ребра, и создает возможность для удаления лишнего (c,d) ребра, как принадлежащего цветному зеленому диску {(см. рис. 3.87). Удаляем лишнее ребро (c,d) и осуществляем перекраску зеленого цветного диска (см. рис. 3.88). Теперь можно раскрасить вновь введенное (a,b) ребро (см. рис. 3.89).

**3.2. Метод раскраски ребер плоского кубического графа с циклами четной длины**

Рассмотрим другой способ раскраски ребер в плоском кубическом графе без мостов. Данный способ основан на свойствах наличия четных изометрических циклов в плоском двудольном кубическом графе. Если задан плоский кубический граф без мостов, то путем введения определенного числа новых ребер плоский кубический граф, сводится к ближайшему двудольному плоскому кубическому графу. Ребра полученного таким путем двудольного кубического графа, в изометрических циклах четной длины, раскрашиваются только в два цвета. Данный способ очень удобен для раскраски ребер в полученном двудольном плоском кубическом графе.

Рассмотрим метод построения плоского кубического графа с изометрическими циклами четной длины для произвольного плоского кубического графа без мостов. Это вызвано возможностью применения простого способа раскраски ребер такого четного цикла в два цвета [9,11].

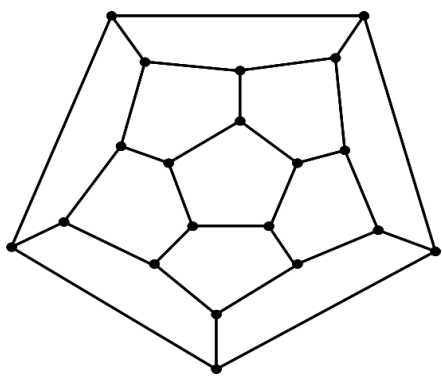 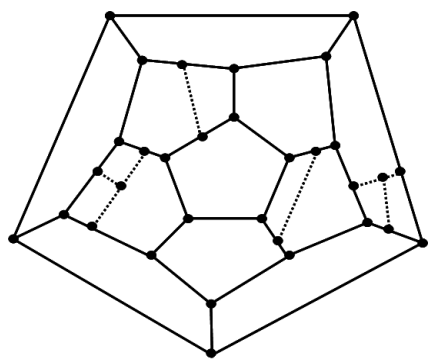 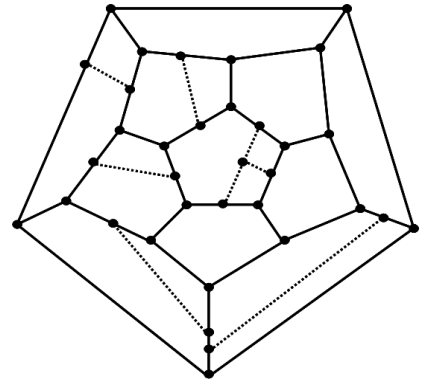

Рис. 3.90. Плоский кубический граф H с изометрическими циклами нечетной длины.

Рис. 3.91. Плоский кубический граф с изометрическими циклами четной длины H′.

После построения плоского кубического графа с изометрическими циклами четной длины производится удаление вновь введенных ребер. Естественно, что удаление ребра можно производить только в случае если ребро принадлежит диску соответствующего цвета.

Опишем данный алгоритм раскраски ребер плоского кубического графа.

**Шаг 0**. [Инициализация]. Задан плоский кубический граф без мостов.

**Шаг 1**. [Соединение двух изометрических циклов]. Соединяем два изометрических цикла нечетной длины ребрами. Ребра проводятся таким образом, чтобы получить изометрические циклы четной длины, такое построение всегда возможно согласно теоремы 1.3. Помечаем введенные ребра.

**Шаг 2**. [Построение плоского кубического графа с изометрическими циклами четной длины]. Если существуют в графе изометрические циклы нечетной длины, то идем на шаг 1. Иначе идем на шаг 3.

**Шаг 3**. [Удаление помеченного ребра]. Применяя операцию ротации для дисков, добиваемся построения одного цветного диска проходящего по концевым вершинам помеченного к удалению ребра. После получения такого цветного диска удаляем помеченное ребро. Осуществляем перекраску ребер в цветном диске.

**Шаг 4**. [Раскраска исходного плоского кубического графа]. Выбираем следующее, помеченное к удалению ребро, и идем на шаг 3. Если помеченных к удалению ребер нет, то конец рабо-ты алгоритма. В результате получим плоский кубический граф без мостов хроматический класс, которого равен трем.

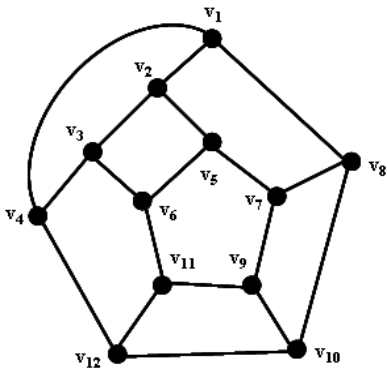 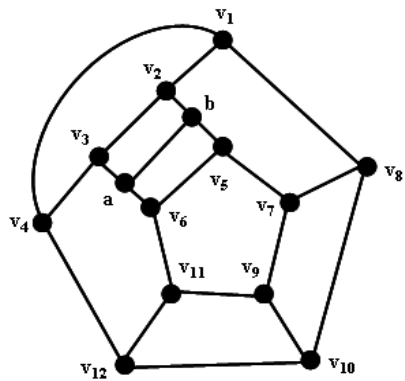 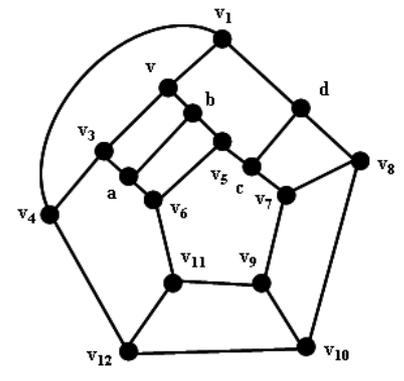

Рис. 3.92. Плоский кубический граф H.   Рис. 3.93. Проведение ребра (a,b).   Рис. 3.94. Проведение ребра (c,d).

**Пример 3.5.** Рассмотрим следующий плоский кубический граф H представленный на рис. 3.92. Преобразуем данный граф в граф с четными изометрическими циклами. На рис. 3.93. показано введение ребра (a,b), а на рис. 3.94 показано введение ребра (c,d). После этого получим кубический граф с изометрическими циклами четной длины, который легко раскрасить тремя цветами.

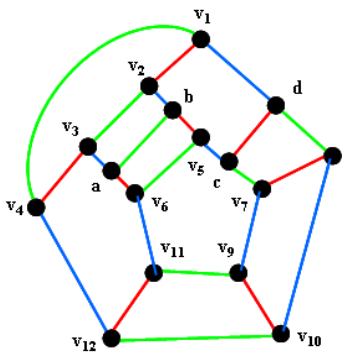 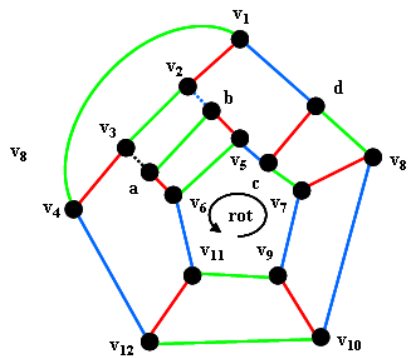 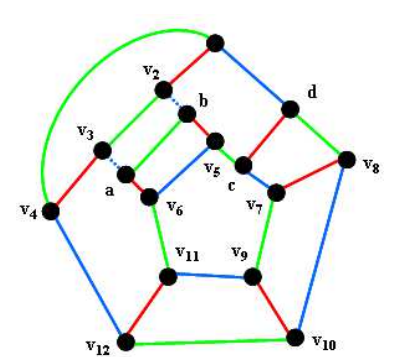

Рис. 3.95. Раскраска плоского кубического графа H′.   Рис. 3.96. Определение диска для ротации.   Рис. 3.97. Ротация диска $\langle v_c, v_5, v_6, v_{11}, v_9, v_7 \rangle$.

Будем удалять ребро (a,b). Концевые вершины данного ребра, принадлежат двум разным зеленым дискам $\langle b, v_2, v_1, d, c, v_5 \rangle$ и $\langle a, v_3, v_4, v_{12}, v_{11}, v_6 \rangle$ (см. рис. 3.95). Для получения только одного зеленого диска проходящего по концевым вершинам ребра (a,b) необходимо произвести ротацию диска $\langle c, v_5, v_6, v_{11}, v_9, v_7 \rangle$ (см. рис. 3.96). После проведения ротации образовался один зеленый диск $\langle a, v_3, v_4, v_{12}, v_{11}, v_9, v_{10}, v_8, v_7, c, d, v_1, v_2, b, v_5, v_6 \rangle$ (см. рис.3.97). Удаляем ребро (a,b) и производим перекраску ребер (см. рис. 3.98).

Будем удалять ребро (c,d). Концевые вершины данного ребра принадлежат двум разным красным дискам $\langle c, v_5, v_2, v_3, v_6, v_{11}, v_9, v_7 \rangle$ и $\langle d, v_1, v_4, v_{12}, v_{10}, v_8 \rangle$ (см. рис. 3.98). Для получения толь-ко одного красного диска проходящего по концевым вершинам ребра (c,d) необходимо произвести ротацию диска <$v_1,v_2,v_3,v_4$> (см. рис. 3.99). После проведения ротации образовался один красный диск <$v_1,v_2,v_5,c,v_7,v_9,v_{11},v_6,v_3,v_4,v_{12},v_{10},v_8,d$> (см. рис. 3.100). Удаляем ребро (c,d) и производим перекраску ребер (см. рис. 3.101), тем самым получили раскраску ребер исходного графа **H**.

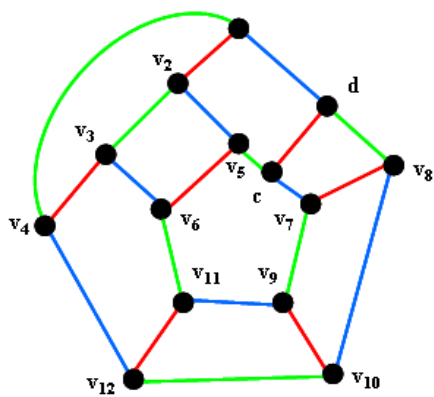 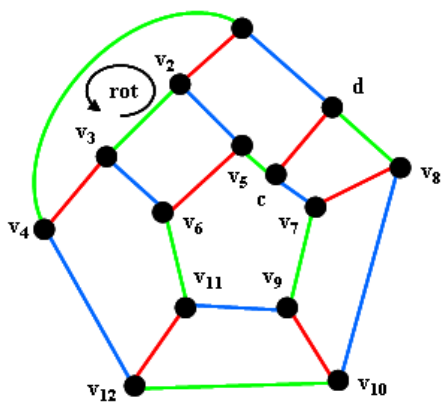 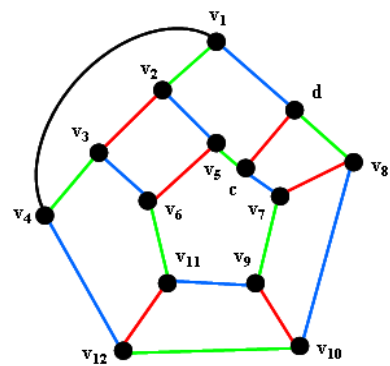

| Рис. 3.98. Удаление ребра (a,b) и перекраска ребер. | Рис. 3.99. Определение диска для ротации. | Рис. 3.100. Ротация диска {$v_2,v_3,v_5,v_6$}. |
|---|---|---|

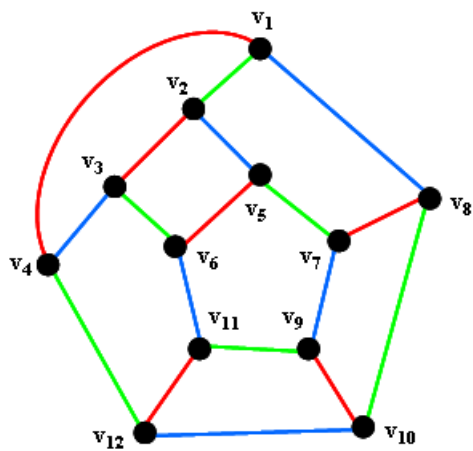

Рис. 3.101. Раскраска ребер исходного плоского кубического графа **H**.

Для применения операции удаления или вставки ребра, нужно обязательно найти цветной диск, проходящий по сцепленным ребрам. Однако в общем случае, поиск такого диска может вызывать затруднения и представляет собой определенное искусство.

### Выводы

Определим вычислительную сложность рекурсивного алгоритма раскраски. Очевидно, что для получения минимального раскрашиваемого кубического графа $K_4$ из исходного графа, в общем нужно затратить операций удаления ребер равное количеству ребер графа. Затем, для получения раскраски одного ребра, может потребоваться перекраска всех ребер графа. В об-щем случае потребуется $m+m^2$ операций удаления и перекраски ребер. Таким образом, вычислительная сложность алгоритма определится как $o(m^2)$. Вычислительная сложность алгоритма раскраски ребер для плоских кубических графов с циклами четной длины, очевидно того же порядка. При этом нужно еще добавить время на построение в произвольном плоском кубическом графе циклов четной длины.

## Глава 4. Раскраска плоских биквадратных графов
### 4.1. Биквадратный граф и его свойства

Плоский биквадратный граф $B_2$ можно построить как рёберный граф к кубическому графу H. Это проиллюстрировано на рис. 4.1, где вершины закрашены красным цветом, а ребра точечными линиями.

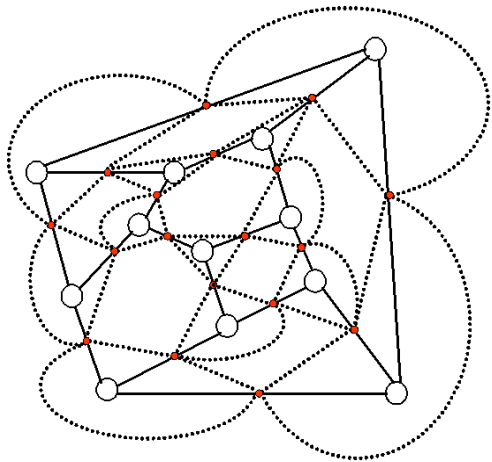

Рис. 4.1. Построение биквадратного графа $B_2$ как дуального реберного графа к кубическому графу H.

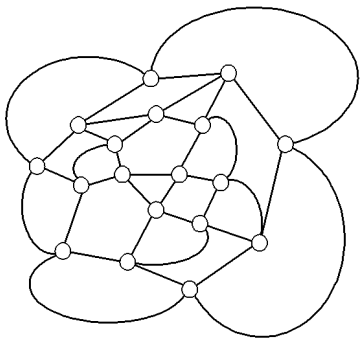

Рис. 4.2. Плоский биквадратный граф $B_2$.

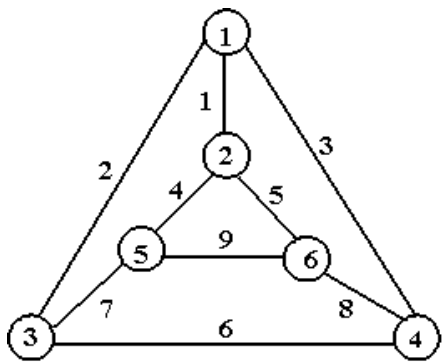

Рис. 4.3. Плоский кубический граф H.

**Пример 4.1.** Рассмотрим плоский кубический граф H представленный на рис. 4.3. Матрица инциденций данного графа имеет вид:

| ребра \ вершины | $e_1$ | $e_2$ | $e_3$ | $e_4$ | $e_5$ | $e_6$ | $e_7$ | $e_8$ | $e_9$ |
|---|---|---|---|---|---|---|---|---|---|
| $v_1$ | 1 | 1 | 1 | | | | | | |
| $v_2$ | 1 | | | 1 | 1 | | | | |
| $v_3$ | | 1 | | | | 1 | 1 | | |
| $v_4$ | | | 1 | | | 1 | | 1 | |
| $v_5$ | | | | 1 | | | 1 | | 1 |
| $v_6$ | | | | | 1 | | | 1 | 1 |

Матрица смежностей биквадратного графа $B_2$ образуется как произведение транспонированной матрицы инциденций и матрицы инциденций для кубического графа H, где номера ребер графа H соответствуют номерам вершин графа $B_2$:

$$A = B \times B^T \quad (4.1)$$

|   | $v_1$ | $v_2$ | $v_3$ | $v_4$ | $v_5$ | $v_6$ | $v_7$ | $v_8$ | $v_9$ |
|---|---|---|---|---|---|---|---|---|---|
| $v_1$ | | 1 | 1 | 1 | 1 | | | | |
| $v_2$ | 1 | | 1 | | | 1 | 1 | | |
| $v_3$ | 1 | 1 | | | | 1 | | 1 | |
| $v_4$ | 1 | | | | 1 | | 1 | | 1 |
| $v_5$ | 1 | | | 1 | | | | 1 | 1 |
| $v_6$ | | 1 | 1 | | | | 1 | 1 | |
| $v_7$ | | 1 | | 1 | | 1 | | | 1 |
| $v_8$ | | | 1 | | 1 | 1 | | | 1 |
| $v_9$ | | | | 1 | 1 | | 1 | 1 | |

(B = матрица слева)

Здесь сложение и умножение производится по законам булевой алгебры:

$1 + 1 = 0; 1 + 0 = 1; 0 + 0 = 0; 1 \times 1 = 1; 1 \times 0 = 0$.

Построенный таким образом биквадратный граф $B_2$ обладает следующими свойствами:

- количество вершин реберного биквадратного графа соответствует количеству ребер кубического графа H;
- количество базовых треугольных граней соответствует количеству вершин в кубическом графе H (см. рис. 4.4);

- количество не базовых граней и их размер соответствует количеству и размеру граней кубического графа H.

На рис. 4.4 затемнены базовые треугольные грани в плоском биквадратном графе $B_2$.

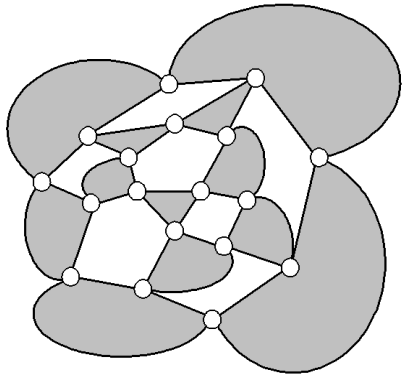

Рис. 4.4 Базовые треугольные грани в биквадратном графе $B_2$ дуальные вершинам в кубическом графе H.

Так как ребра максимально плоского графа однозначно соответствуют ребрам двойственного кубического графа H, то теорему Тэйта можно выразить в следующей эквивалентной форме.

**Теорема 4.1**. Для того чтобы хроматическое число максимально плоского графа $G'$ было равно четырем, необходимо и достаточно, чтобы его ребра допускали такую раскраску тремя цветами, при которой никакие два ребра одной грани не были бы окрашены одним цветом.

С другой стороны, ребра плоского биквадратного графа также можно раскрасить в три цвета, так как цвет ребра может быть определен операцией преобразования группы Клейна.

На рис. 4.5 представлен биквадратный граф $B_2$ и его раскраска вершин и ребер.

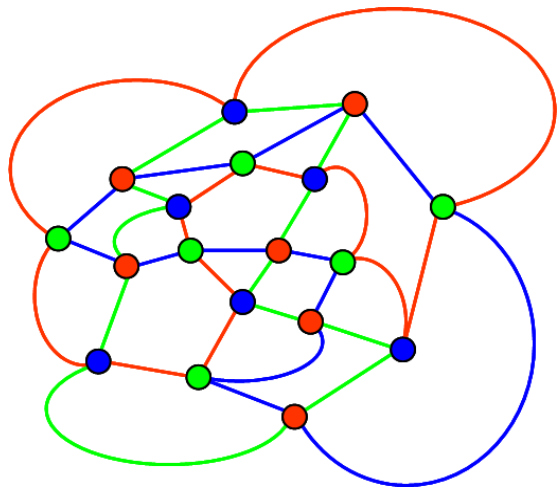

Рис. 4.5. Раскраска вершин и ребер в биквадратном графе $B_2$.

## 4.2. Свойства эйлерова цикла в биквадратном плоском графе

Эйлер доказал следующую теорему.

**Теорема 4.1** [16]. *Связный граф является эйлеровым тогда и только тогда, когда степени всех его вершин четны.*

Будем называть *диском эйлерова цикла* замкнутый ориентированный маршрут, исходящий из вершины с номером $i$ и входящий в туже вершину.

Согласно доказанной теореме в любом плоском биквадратном графе существуют эйлеров цикл [16]. Докажем следующую теорему:

**Теорема 4.2**. *Для того, чтобы хроматический класс графа H был равен трем, необходимо и достаточно существование в однородном реберном биквадратном графе $B_2$ эйлерова цикла, где длина каждого диска для произвольно выбранной вершины кратна трем.*

Доказательство.

**Необходимость.** Предположим, что существует раскраска вершин биквадратного графа $B_2$ тремя цветами, причем такая, что каждой вершине инцидентны ребра только двух цветов, и только попарно. Выделим в графе $B_2$ две треугольные грани соединённые вершиной. Соответственно, для этой вершины образуются два контура длиной три. Зададим ориентацию каждого контура относительно выбранной раскраски, допустим 1,2,3. Присоединим через вершину новую треугольную грань к какому-нибудь контуру. Тем самым задаётся новая ориентация контура длиной кратной трем с учетом выбранной раскраски. Так как граф **$B_2$** конечен, то рано или поздно, продолжая далее этот процесс присоединения треугольных граней к контуру, мы получим эйлеров цикл с двумя контурами длиной кратной трем для каждой вершины с ориентацией, соответствующей выбранной раскраске.

**Достаточность.** Пусть найден эйлеров цикл, в котором длина каждого контура для произвольно выбранной вершины кратна трем. Тогда возможно последовательно пронумеровать ребра такого эйлерова цикла начиная с произвольно выбранного ребра целыми положительными числами от 1 до m (m – число рёбер биквадратного графа $B_2$, m = 3n). Заметим, что количество ребер в биквадратном графе $B_2$ всегда кратно трем. Раскрасим ребра. Тогда, если для произвольно выбранной вершины исходящее по контуру эйлерова цикла ребро помечено числом кратным единице по mod 3, то входящее по контуру эйлерова цикла ребро будет помечено числом кратным нулю по mod 3. Это следует из того факта, что длина контура эйлерова цикла кратна трем. То же происходит и для другого контура эйлерова цикла относительно выбранной вершины. Данную вершину можно раскрасить во второй цвет. Если для произвольно выбран-ной вершины исходящее по контуру эйлерова цикла ребро будет помечено числом кратным двум по mod 3, то входящее ребро по контуру эйлерова цикла будет помечено числом кратным единице по mod 3. Тогда данную вершину можно раскрасить третьим цветом. Если для произвольно выбранной вершины исходящее по петле эйлерова цикла ребро будет помечено числом кратным нулю по mod 3, то входящее ребро по контуру эйлерова цикла будет помечено числом кратным двум по mod 3 и тогда данную вершину можно раскрасить первым цветом. Таким образом, хроматический класс однородного кубического графа Н будет равен трем.

Теорема 4.2 доказана.

В дальнейшем нам понадобятся методы описания рисунка плоского графа на плоскости [12].

Необходимым понятием для описания топологического плоского рисунка графа **G** является понятие о вращении вершин графа, введенное Г. Рингелем [12].

**Определение 4.1.** Для данного графа G вращение вершины А – это ориентированный циклический порядок (или циклическая перестановка) всех рёбер инцидентных вершине А.

Пометим вершины графа числами 1,2,...,n. Если вершина А имеет степень три, а 1,2,3 – три вершины смежные с вершиной А, то имеются две различные возможности для задания вращения вершины А. Будем описывать вращение вершины, указывая вместо рёбер инцидентных этой вершине циклический порядок вершин смежных с ней. Тогда указанные две возможности это:

( 1 2 3 ) = ( 2 3 1 ) = ( 3 1 2 )

и ( 3 2 1 ) = ( 2 1 3 ) = ( 1 3 2 ).

Вообще, число возможных вращений вершины степени n равно:

$$(n-1)! = 1 \cdot 2 \ldots n-1 \qquad (4.2)$$

**Определение 4.2.** Вращение $\hbar$ графа G – это вращение всех вершин графа G.

Запись (G, $\hbar$ ) будет обозначать граф G с некоторым вращением $\hbar$ , которое характеризует топологический рисунок графа на плоскости.

Граф G с вращением часто бывает удобно изображать на плоскости таким образом, чтобы читая рёбра инцидентные некоторой вершине или вершины смежные с этой вершиной по или против часовой стрелки мы получили вращение в этой вершине.

Вращение графа можно описывать и представлять следующим образом. Обозначим вершины числами 1,2,...,n. Затем выпишем циклическую перестановку соседей для каждой вершины $i$. Эта перестановка порождается вращением вершины $i$, которое является циклической перестановкой рёбер, инцидентных вершине $i$. Например, графу с вращением, показанному на рис. 4.6, соответствует схема, называемая диаграммой вращения вершин:

| | | | | |
|---|---|---|---|---|
| $v_1$: | $v_2$ | $v_3$ | $v_4$ | |
| $v_2$: | $v_5$ | $v_3$ | $v_1$ | |
| $v_3$: | $v_1$ | $v_2$ | $v_5$ | $v_4$ |
| $v_4$: | $v_1$ | $v_3$ | $v_5$ | |
| $v_5$: | $v_4$ | $v_3$ | $v_2$ | |

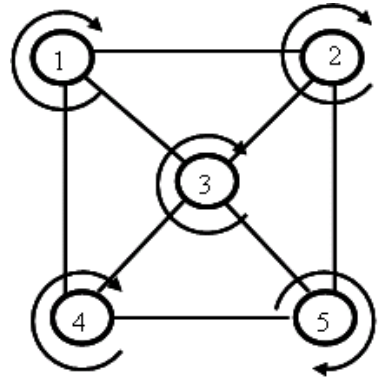

Рис. 4.6. Граф и его вращение вершин

Пусть $v_1$ – вершина, инцидентная ребру $e_1$ в графе G с вращением $(G, \hbar)$. Мы построим в графе G замкнутый маршрут

$$v_1, e_1, v_2, e_2, v_3, e_3, ..., \qquad (4.3)$$

где: вершина $v_2$ – второй конец ребра $e_1$, а ребро $e_2$ следует за ребром $e_1$ во вращении вершины $v_2$, определяемое вращением $\hbar$. Затем определяется $v_3$ как вершина инцидентная ребру $e_2$ и не равная $v_2$. После этого в качестве $e_3$ выбирается ребро, следующее за ребром $e_2$ во вращении вершины $v_3$ и т.д. Закончим процесс в точности перед тем моментом, когда должна повториться пара $v_1, e_1$. Она должна повториться, ибо граф G конечный, а наш процесс однозначно определен и в обратном направлении. А именно, если часть $v_{t-1}, e_t, v_t, ...$ известна, то ребро $e_{t-1}$ определяется вращением вокруг вершины $v_{t-1}$. Мы назовем такой замкнутый маршрут циклом, порожденным вершиной $v_1$ и ребром $e_1$, и индуцированным вращением.

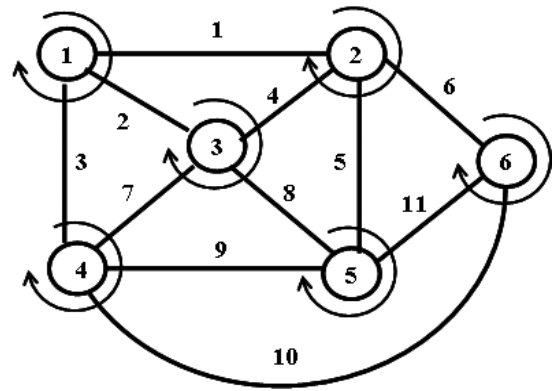

Рис. 4.7. Граф G и его вращение вершин $\hbar$.

Заметим, что во вращении каждое ребро появляется в точности дважды, второй раз – всегда в противоположном направлении.

Если граф планарен и имеется вращение, описывающее плоский рисунок, то циклы индуцированные вращением – суть простые циклы. Например, для плоского графа G с вращением, представленным на рис. 4.7, имеем следующую систему индуцированных циклов.

Диаграмма вращений 4.7 (для рис. 4.7, но запись произведена в реберном представлении).

| | | | | |
|---|---|---|---|---|
| $v_1$: | $e_1$ | $e_2$ | $e_3$ | |
| $v_2$: | $e_6$ | $e_5$ | $e_4$ | $e_1$ |
| $v_3$: | $e_2$ | $e_4$ | $e_8$ | $e_7$ |
| $v_4$: | $e_3$ | $e_7$ | $e_9$ | $e_{10}$ |
| $v_5$: | $e_9$ | $e_8$ | $e_5$ | $e_{11}$ |
| $v_6$: | $e_6$ | $e_{10}$ | $e_{11}$ | |

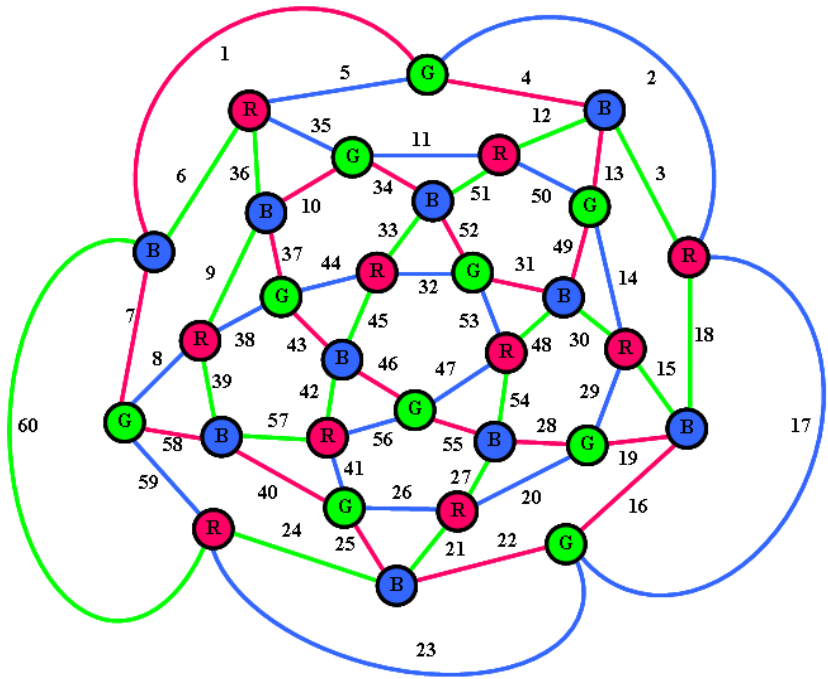

Рис. 4.8. Раскраска рёбер и вершин плоского биквадратного графа $B_2$.

Индуцированные циклы: $<v_1,v_3,v_2>$, $<v_1,v_4,v_3>$, $<v_2,v_3,v_5>$, $<v_5,v_3,v_4>$, $<v_2,v_5,v_6>$, $<v_6,v_5,v_4>$, $<v_2,v_6,v_4,v_1>$ для вершин. Или в виде множества рёбер: $\{e_1,e_2,e_4\}$, $\{e_2,e_3,e_7\}$, $\{e_7,e_8,e_9\}$, $\{e_4,e_5,e_8\}$, $\{e_5,e_6,e_{11}\}$, $\{e_9,e_{10},e_{11}\}$, $\{e_1,e_3,e_6,e_{10}\}$.

В общем случае вращение вершин графа $\hbar$ позволяет описывать топологический рисунок графа на плоскости.

Выделим эйлеров цикл [16] в плоском биквадратном графе $B_2$. На рис. 4.8 эйлеров цикл, выделенный согласно условиям теоремы 4.2, представлен в виде последовательности пронумерованных рёбер.

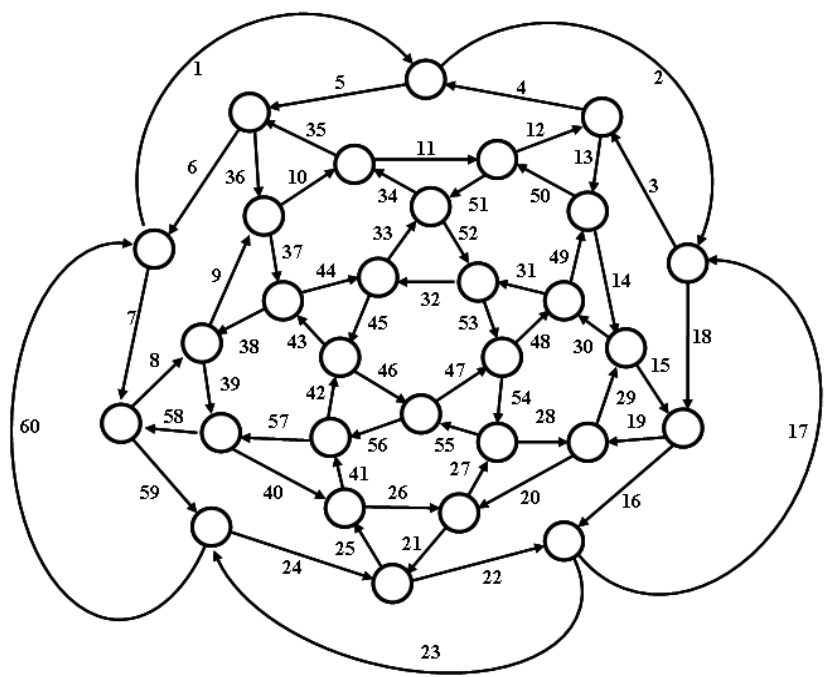

Рис. 4.9. Эйлеров цикл плоского биквадратного графа.

Выделенную таким образом последовательность рёбер можно раскрасить тремя цветами, соблюдая порядок цветов, например: R – B – G (см. рис. 4.9).

В общем случае поиск такого эйлерового цикла с треугольными дисками является трудно решаемой переборной задачей [5,8,11]. Однако если известна раскраска вершин биквадратного графа в три цвета, то поиск такого пути значительно облегчается.

Введем определение вращения базовой треугольной грани биквадратного графа $B_2$ для соответствующей вершины графа H.

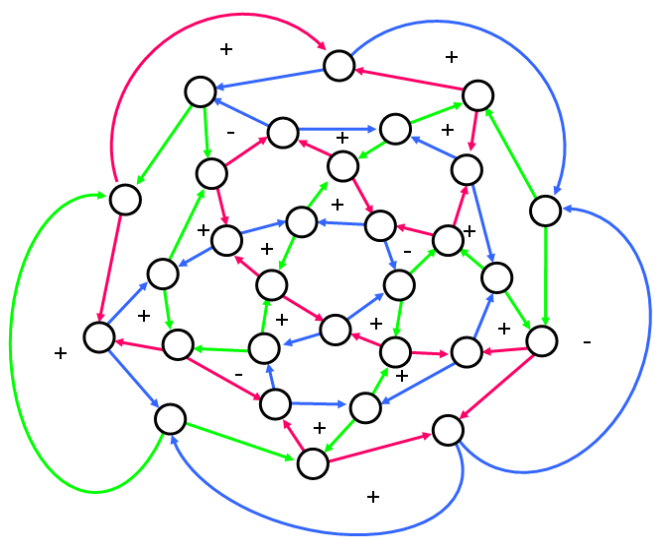

Рис. 4.10. Вращение базовых треугольных циклов в биквадратном графе $B_2$.

**Определение 4.3.** Вращение базовой треугольной грани биквадратного графа B$_2$ для соответствующей вершины однородного кубического графа H – это замкнутый ориентированный маршрут состоящий из трёх дуг треугольной грани направленный по или против часовой стрелки относительно рисунка графа на плоскости.

Будем говорить, что вращение базовой треугольной грани имеет знак «+», если вращение происходит по часовой стрелке, а знак «–», если вращение происходит против часовой стрелки (см. рис. 4.10). Эйлеров цикл, у которого все диски имеют длину кратную трём, порождает вращение базовых треугольных граней и определяет их знак вращения.

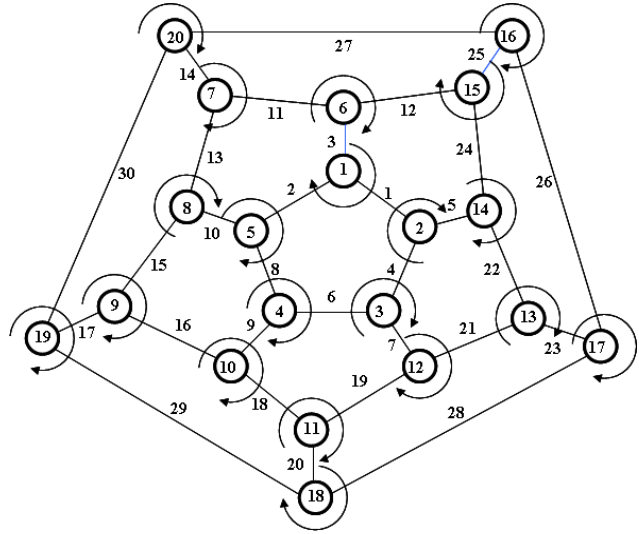

Рис. 4.11. Вращение вершин для топологического рисунка графа H.

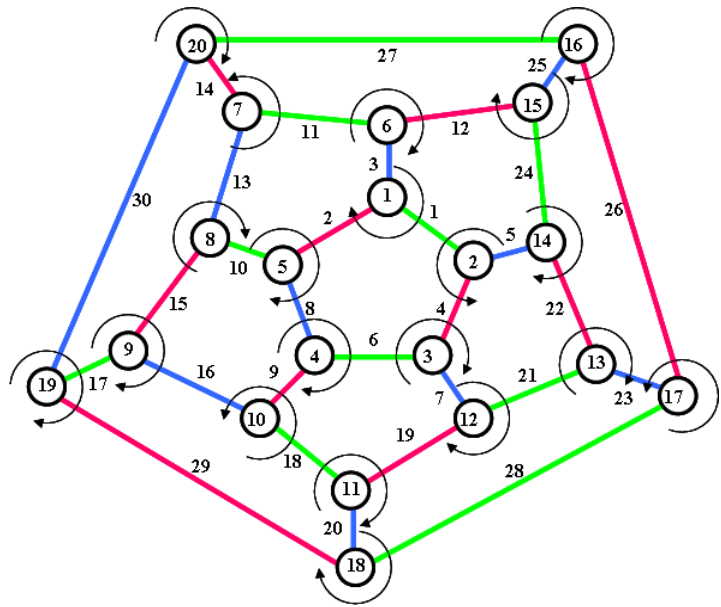

Рис. 4.12. Цветное вращение вершин кубического графа H.

С другой стороны, вращению базовых треугольных граней не соответствует вращение вершин $\hbar$ в плоском топологическом рисунке кубического графа (см. рис. 4.11). Цветному вращению базовых треугольных граней (см. рис. 4.12) соответствует вращение вершин, которое можно записать в виде диаграммы:

| | | | |
|---|---|---|---|
| $v_1$: | $v_2$ | $v_5$ | $v_6$ |
| $v_2$: | $v_1$ | $v_3$ | $v_{14}$ |
| $v_3$: | $v_2$ | $v_{12}$ | $v_4$ |
| $v_4$: | $v_5$ | $v_3$ | $v_{10}$ |
| $v_5$: | $v_4$ | $v_8$ | $v_1$ |
| $v_6$: | $v_7$ | $v_{15}$ | $v_1$ |
| $v_7$: | $v_8$ | $v_6$ | $v_{20}$ |
| $v_8$: | $v_7$ | $v_5$ | $v_9$ |
| $v_9$: | $v_8$ | $v_{10}$ | $v_{19}$ |
| $v_{10}$: | $v_4$ | $v_9$ | $v_{11}$ |
| $v_{11}$: | $v_{10}$ | $v_{12}$ | $v_{18}$ |
| $v_{12}$: | $v_3$ | $v_{13}$ | $v_{11}$ |
| $v_{13}$: | $v_{14}$ | $v_{17}$ | $v_{12}$ |
| $v_{14}$: | $v_2$ | $v_{15}$ | $v_{13}$ |
| $v_{15}$: | $v_6$ | $v_{16}$ | $v_{14}$ |
| $v_{16}$: | $v_{20}$ | $v_{17}$ | $v_{15}$ |
| $v_{17}$: | $v_{13}$ | $v_{18}$ | $v_{16}$ |
| $v_{18}$: | $v_{11}$ | $v_{17}$ | $v_{19}$ |
| $v_{19}$: | $v_{20}$ | $v_9$ | $v_{18}$ |
| $v_{20}$: | $v_{19}$ | $v_{16}$ | $v_7$ |

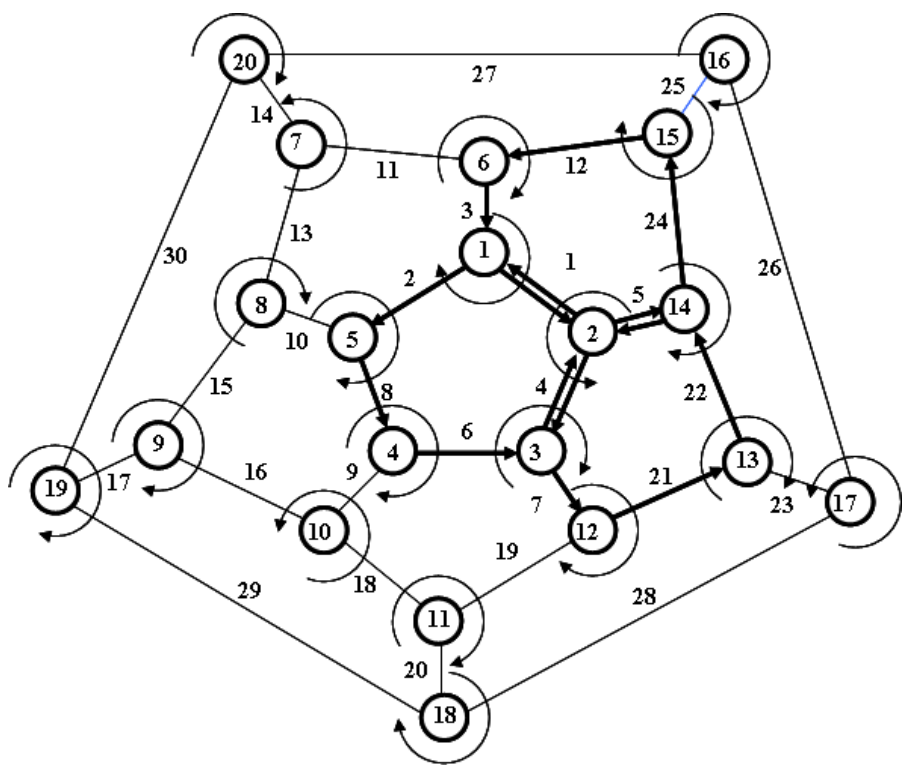

Рис. 4.13. Цикл $c_1$ индуцированный цветным вращением вершин.

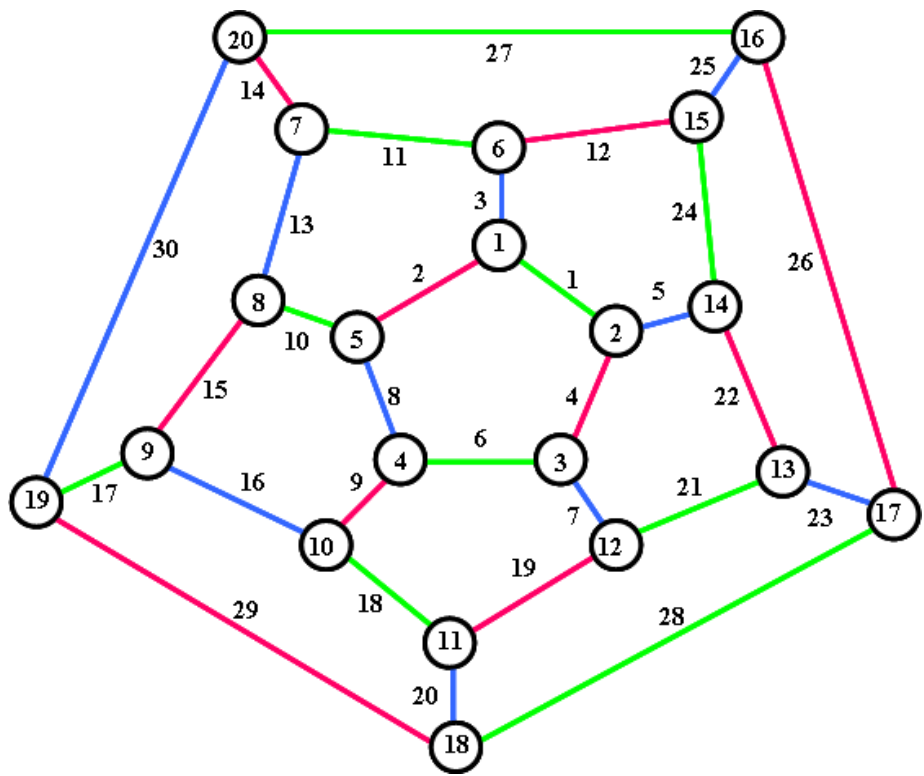

Рис. 4.14. Раскраска рёбер плоского кубического графа.

Данное вращение вершин индуцирует следующую систему ориентированных циклов, где длина каждого цикла кратна трем. Причем по некоторым ребрам цикл проходит дважды, однако, в разных направлениях (см. рис. 4. 13):

$c_1 = \langle v_1,v_5,v_4,v_3,v_2,v_{14},v_{15},v_6,v_1,v_2,v_3,v_{12},v_{13},v_{14},v_2 \rangle$;
$c_2 = \langle v_5,v_1,v_6,v_7,v_{20},v_{19},v_9,v_8,v_7,v_6,v_{15},v_{16},v_{20},v_7,v_8 \rangle$;
$c_3 = \langle v_{15},v_{14},v_{13},v_{17},v_{18},v_{19},v_{20},v_{16},v_{17},v_{13},v_{12},v_{11},v_{18},v_{17},v_{16} \rangle$;
$c_4 = \langle v_{11},v_{12},v_3,v_4,v_{10},v_9,v_{19},v_{18},v_{11},v_{10},v_4,v_5,v_8,v_9,v_{10} \rangle$.

Данную систему циклов можно записать в виде множества рёбер, сохраняя последовательность цветов R – B – G (см. рис. 4.14):

$c_1 = \langle e_2,e_8,e_6,e_4,e_5,e_{24},e_{12},e_3,e_1,e_4,e_7,e_{21},e_{22},e_5,e_1 \rangle$;
$c_2 = \langle e_2,e_3,e_{11},e_{14},e_{30},e_{17},e_{15},e_{13},e_{11},e_{12},e_{25},e_{27},e_{14},e_{13},e_{10} \rangle$;
$c_3 = \langle e_{22},e_{23},e_{28},e_{29},e_{30},e_{27},e_{26},e_{23},e_{21},e_{19},e_{20},e_{28},e_{26},e_{25},e_{24} \rangle$;
$c_4 = \langle e_{19},e_7,e_6,e_9,e_{16},e_{17},e_{29},e_{20},e_{18},e_9,e_8,e_{10},e_{15},e_{17},e_{16} \rangle$.

Таким образом, вращение $\hbar$ описывает рисунок плоского кубического графа и индуцирует изометрические циклы. А вращение, представленное на рис. 4.12, индуцирует цветные циклы с длиной кратной трём, которые раскрашиваются тремя цветами в последовательности R – B – G. Такое вращение будем называть цветным вращением плоского кубического графа.

### 4.3. Уравнение Хивуда

Для множества базовых треугольных граней $M_\alpha$ плоского биквадратного графа $B_2$, примыкающих к не базовой грани, можно записать систему уравнений Хивуда [1]:

$$\begin{cases} w_i^2 = 1, \\ \sum_{i \in M_\alpha} w_i \equiv 0 (\mod 3) \end{cases} \quad (4.4)$$

где $i = 1,2,\ldots,n$. Здесь n – количество вершин в соответствующем максимально плоском графе равное количеству не базовых граней биквадратного графа.

В свою очередь, для множества граней $M_\alpha$ примыкающих к вершине $i$ в исходном плоском графе $G'$, состоящем из треугольных граней, можно также записать систему уравнений Хивуда (4.4) для плоской триангуляции [1,7].

Если перейти к рассмотрению плоского биквадратного графа или максимально плоского графа $G'$, то треугольные грани перейдут во вращающиеся ориентированные треугольные грани. При этом решение системы Хивуда определяет раскраску, так как вращению граней по часовой стрелке можно поставить в соответствие «+1» в системе Хивуда, а вращению граней против часовой стрелки можно поставить в соответствие «-1».

**Пример 4.2**. Рассмотрим плоский биквадратный граф $B_2$ представленный на рис. 4.8 – 4.9. Будем отождествлять грань с изометрическим циклом или ободом. Тогда можно записать систему циклов описывающих базовые треугольные грани:

$c_1 = \{e_1, e_5, e_6\}$; $c_2 = \{e_2, e_3, e_4\}$; $c_3 = \{e_{16}, e_{17}, e_{18}\}$; $c_4 = \{e_{22}, e_{23}, e_{24}\}$;
$c_5 = \{e_7, e_{59}, e_{60}\}$; $c_6 = \{e_{10}, e_{35}, e_{36}\}$; $c_7 = \{e_{12}, e_{13}, e_{50}\}$; $c_8 = \{e_{15}, e_{19}, e_{29}\}$;
$c_9 = \{e_{21}, e_{25}, e_{26}\}$; $c_{10} = \{e_8, e_{39}, e_{58}\}$; $c_{11} = \{e_{11}, e_{34}, e_{51}\}$; $c_{12} = \{e_{14}, e_{30}, e_{49}\}$;
$c_{13} = \{e_{20}, e_{27}, e_{28}\}$; $c_{14} = \{e_{40}, e_{41}, e_{57}\}$; $c_{15} = \{e_9, e_{37}, e_{38}\}$; $c_{16} = \{e_{43}, e_{44}, e_{45}\}$;
$c_{17} = \{e_{32}, e_{33}, e_{52}\}$; $c_{18} = \{e_{31}, e_{48}, e_{53}\}$; $c_{19} = \{e_{47}, e_{54}, e_{55}\}$; $c_{20} = \{e_{42}, e_{46}, e_{56}\}$.

Система циклов описывающих не базовые грани:

$c_{21} = \{e_4, e_5, e_{11}, e_{12}, e_{35}\}$; $c_{22} = \{e_3, e_{13}, e_{14}, e_{15}, e_{18}\}$; $c_{23} = \{e_{16}, e_{19}, e_{20}, e_{21}, e_{22}\}$;
$c_{24} = \{e_{24}, e_{25}, e_{40}, e_{58}, e_{59}\}$; $c_{25} = \{e_6, e_7, e_8, e_9, e_{36}\}$; $c_{26} = \{e_{10}, e_{33}, e_{34}, e_{37}, e_{44}\}$;
$c_{27} = \{e_{31}, e_{49}, e_{50}, e_{51}, e_{52}\}$; $c_{28} = \{e_{28}, e_{29}, e_{30}, e_{48}, e_{54}\}$; $c_{29} = \{e_{26}, e_{27}, e_{41}, e_{55}, e_{56}\}$;
$c_{30} = \{e_{38}, e_{39}, e_{42}, e_{43}, e_{57}\}$; $c_{31} = \{e_{32}, e_{45}, e_{46}, e_{47}, e_{53}\}$; $c_0 = \{e_1, e_2, e_{17}, e_{23}, e_{60}\}$.

Определим вращение базовых треугольных граней:

$c_1 = +1$; $c_2 = +1$; $c_3 = -1$; $c_4 = +1$; $c_5 = +1$; $c_6 = -1$; $c_7 = +1$; $c_8 = +1$;
$c_9 = +1$; $c_{10} = +1$; $c_{11} = +1$; $c_{12} = +1$; $c_{13} = +1$; $c_{14} = -1$; $c_{15} = +1$; $c_{16} = +1$;
$c_{17} = +1$; $c_{18} = -1$; $c_{19} = +1$; $c_{20} = +1$.

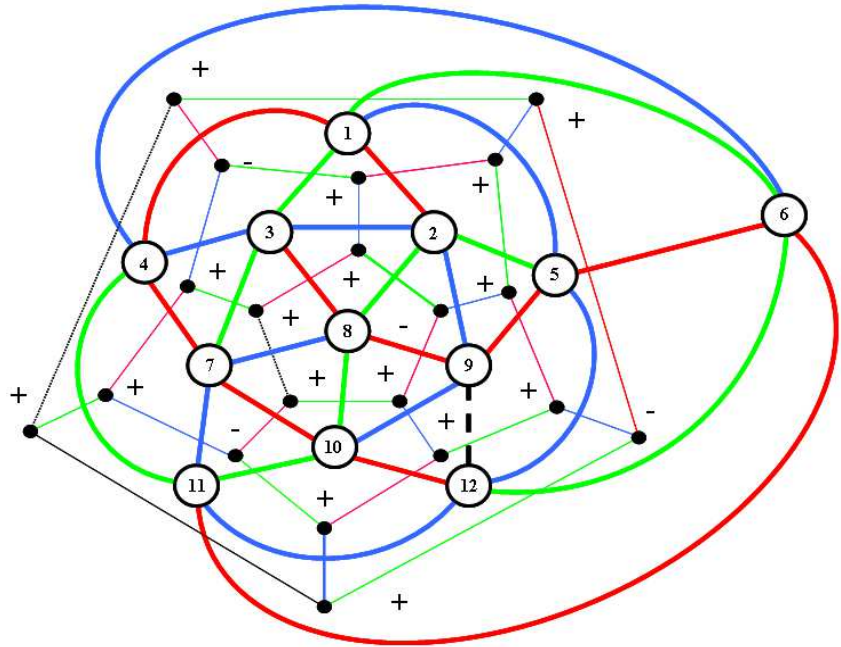

Рис. 4.15. Максимально плоский граф $G'$ и его двойственный кубический граф H.

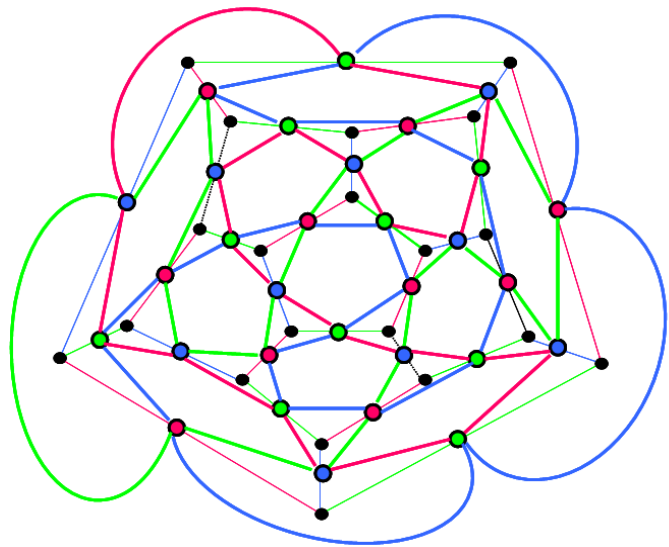

Рис. 4.16. Кубический граф H и его реберный биквадратный граф $B_2$.

Уравнения Хивуда имеют вид:

$w_{21} = c_1 + c_2 + c_6 + c_7 + c_{11} = (+1) + (+1) + (-1) + (+1) + (+1) = 0 \pmod 3$;
$w_{22} = c_2 + c_3 + c_7 + c_8 + c_{12} = (+1) + (-1) + (+1) + (+1) + (+1) = 0 \pmod 3$;
$w_{23} = c_3 + c_4 + c_8 + c_9 + c_{13} = (-1) + (+1) + (+1) + (+1) + (+1) = 0 \pmod 3$;
$w_{24} = c_4 + c_5 + c_9 + c_{10} + c_{14} = (+1) + (+1) + (+1) + (+1) + (-1) = 0 \pmod 3$;
$w_{25} = c_1 + c_5 + c_6 + c_{10} + c_{15} = (+1) + (+1) + (-1) + (+1) + (+1) = 0 \pmod 3$;
$w_{26} = c_6 + c_{11} + c_{15} + c_{16} + c_{17} = (-1) + (+1) + (+1) + (+1) + (+1) = 0 \pmod 3$;
$w_{27} = c_7 + c_{11} + c_{12} + c_{17} + c_{18} = (+1) + (+1) + (+1) + (+1) + (-1) = 0 \pmod 3$;
$w_{28} = c_8 + c_{12} + c_{13} + c_{18} + c_{19} = (+1) + (+1) + (+1) + (-1) + (+1) = 0 \pmod 3$;
$w_{29} = c_9 + c_{13} + c_{14} + c_{19} + c_{20} = (+1) + (+1) + (-1) + (+1) + (+1) = 0 \pmod 3$;
$w_{30} = c_{10} + c_{14} + c_{15} + c_{16} + c_{20} = (+1) + (-1) + (+1) + (+1) + (+1) = 0 \pmod 3$;

$w_{31} = c_{16} + c_{17} + c_{18} + c_{19} + c_{20} = (+1) + (+1) + (-1) + (+1) + (+1) = 0 \pmod 3$;

$w_0 = c_1 + c_2 + c_3 + c_4 + c_5 = (+1) + (+1) + (-1) + (+1) + (+1) = 0 \pmod 3$.

**Пример 4.3.** Рассмотрим плоский граф G, представленный на рис. 4.15. Раскрасим ребра в треугольных гранях в три цвета. Множеству треугольных граней $C_e$ поставим в соответствие изометрические циклы:

$c_1 = \{x_1, x_4, x_6\}$; $c_2 = \{x_1, x_3, x_4\}$; $c_3 = \{x_1, x_2, x_5\}$; $c_4 = \{x_1, x_5, x_6\}$;
$c_5 = \{x_3, x_4, x_7\}$; $c_6 = \{x_3, x_7, x_8\}$; $c_7 = \{x_2, x_3, x_8\}$; $c_8 = \{x_2, x_8, x_9\}$;
$c_9 = \{x_2, x_5, x_9\}$; $c_{10} = \{x_4, x_7, x_{11}\}$; $c_{11} = \{x_7, x_8, x_{10}\}$; $c_{12} = \{x_8, x_9, x_{10}\}$;
$c_{13} = \{x_5, x_9, x_{12}\}$; $c_{14} = \{x_5, x_6, x_{12}\}$; $c_{15} = \{x_7, x_{10}, x_{11}\}$; $c_{16} = \{x_{10}, x_{11}, x_{12}\}$;
$c_{17} = \{x_6, x_{11}, x_{12}\}$; $c_{18} = \{x_1, x_2, x_3\}$; $c_{19} = \{x_9, x_{10}, x_{12}\}$; $c_0 = \{x_4, x_6, x_{11}\}$.

Из раскраски ребер R – B – G определим вращение граней:

$c_1 = +1$; $c_2 = -1$; $c_3 = +1$; $c_4 = +1$; $c_5 = +1$; $c_6 = +1$; $c_7 = +1$; $c_8 = -1$;
$c_9 = +1$; $c_{10} = +1$; $c_{11} = +1$; $c_{12} = +1$; $c_{13} = +1$; $c_{14} = -1$; $c_{15} = -1$; $c_{16} = +1$;
$c_{17} = +1$; $c_{18} = +1$; $c_{19} = +1$; $c_{20} = +1$.

Тогда можно записать уравнения Хивуда в виде:

$x_1 = c_1 + c_2 + c_3 + c_4 + c_0 = (+1) + (-1) + (+1) + (+1) + (+1) = 0 \pmod 3$;
$x_2 = c_3 + c_7 + c_8 + c_9 + c_{18} = (+1) + (+1) + (-1) + (+1) + (+1) = 0 \pmod 3$;
$x_3 = c_2 + c_5 + c_6 + c_7 + c_{18} = (-1) + (+1) + (+1) + (+1) + (+1) = 0 \pmod 3$;
$x_4 = c_1 + c_2 + c_5 + c_{10} + c_0 = (+1) + (-1) + (+1) + (+1) + (+1) = 0 \pmod 3$;
$x_5 = c_3 + c_4 + c_9 + c_{13} + c_{14} = (+1) + (+1) + (+1) + (+1) + (-1) = 0 \pmod 3$;
$x_6 = c_1 + c_4 + c_{14} + c_{17} + c_0 = (+1) + (+1) + (-1) + (+1) + (+1) = 0 \pmod 3$;
$x_7 = c_5 + c_6 + c_{10} + c_{11} + c_{15} = (+1) + (+1) + (+1) + (+1) + (-1) = 0 \pmod 3$;
$x_8 = c_6 + c_7 + c_8 + c_{11} + c_{12} = (+1) + (+1) + (-1) + (+1) + (+1) = 0 \pmod 3$;
$x_9 = c_8 + c_9 + c_{12} + c_{13} + c_{19} = (-1) + (+1) + (+1) + (+1) + (+1) = 0 \pmod 3$;
$x_{10} = c_{11} + c_{12} + c_{15} + c_{16} + c_{19} = (+1) + (+1) + (-1) + (+1) + (+1) = 0 \pmod 3$;
$x_{11} = c_{10} + c_{15} + c_{16} + c_{17} + c_0 = (+1) + (-1) + (+1) + (+1) + (+1) = 0 \pmod 3$;
$x_{12} = c_{13} + c_{14} + c_{16} + c_{17} + c_{19} = (+1) + (-1) + (+1) + (+1) + (+1) = 0 \pmod 3$.

Для максимально плоских графов справедлива следующая теорема.

**Теорема 4.3.** [1,2]. Пусть M – множество треугольных граней сферической триангуляции. Тогда для любой правильной раскраски максимально плоского графа:

$$\sum_{\mu \in M} c_\mu \equiv 0 (\mod 4) \qquad (4.5)$$

И действительно:

$c_1 + c_2 + c_3 + c_4 + c_5 + c_6 + c_7 + c_8 + c_9 + c_{10} + c_{11} + c_{12} + c_{13} + c_{14} + c_{15} + c_{16} + c_{17} +$
$+ c_{18} + c_{19} + c_0 = (+1) + (-1) + (+1) + (+1) + (+1) + (+1) + (+1) + (-1) + (+1) +$
$+ (+1) + (+1) + (+1) + (+1) + (-1) + (-1) + (+1) + (+1) + (+1) + (+1) +$
$+ (+1) = 0 \pmod 4$.

Что касается плоских кубических графов Н дуальных максимально плоскому графу, то вращение треугольных граней максимально плоского графа соответствует цветному вращению вершин графа Н с сохранением всех свойств правильной раскраски графа. Здесь уравнения Хивуда описывают цветное вращение вершин относительно каждого изометрического цикла и обода плоского кубического графа Н.

Для плоских биквадратных графов – реберных к плоским кубическим графам цветному вращению вершин соответствует вращение базовых треугольных граней.

### 4.4. Соответствие между графами

Допустим, плоский кубический граф Н дуален к максимально плоскому графу (см. рис. 4.15) и рёберному (см. рис. 4.16) графу $B_2$. Тогда между максимально плоским графом $G'$ и графом Н существуют следующие соответствия:

- количество рёбер максимально плоского графа $G'$ соответствует количеству рёбер плоского кубического графа Н;
- количество треугольных граней максимально плоского графа $G'$ соответствует количеству вершин плоского кубического графа Н;
- количество вершин максимально плоского графа $G'$ соответствует количеству граней плоского кубического графа Н;
- количество вершин биквадратного плоского графа $B_2$ соответствует количеству рёбер плоского кубического графа Н;
- количество базовых треугольных граней биквадратного плоского графа $B_2$ соответствует количеству вершин плоского кубического графа Н;
- количество не базовых треугольных граней биквадратного плоского графа $B_2$ соответствует количеству граней плоского кубического графа Н;
- общее количество граней биквадратного плоского графа $B_2$ соответствует количеству вершин плоского кубического графа Н плюс количество его граней;
- количество рёбер биквадратного плоского графа $B_2$ соответствует утроенному количеству вершин плоского кубического графа Н или удвоенному количеству рёбер плоского кубического графа Н.

### Выводы

Введение в рассмотрение раскраску биквадратного графа позволило получить расширенную картину раскраски плоских графов. Это позволило связать в единое целое раскраску 2-фактора плоского кубического графа и раскраску рёбер эйлерова цикла в плоских биквадратных графах. Рассмотрены вопросы перехода от плоского кубического графа к реберному плоскому биквадратному графу $B_2$. Доказана теорема о необходимости и достаточности раскраски плоского биквадратного графа, на основе построения специального эйлерова цикла с длиной дисков кратным трем. Рассмотрены вопросы построения цветного вращения вершин в плоском кубическом графе и связь с уравнением Хивуда.

Теория вращения вершин графа позволяет строить математические модели для топологического рисунка графа, не производя никаких геометрических представлений на плоскости. Математическая модель топологического рисунка графа является хорошей иллюстрацией вопросов раскраски плоских графов.

**Заключение**

Рассмотрены вопросы перехода от максимально плоского графа к плоскому кубическому графу H. На основании теоремы Тэйта-Волынского [5], доказана теорема о существовании цветного диска проходящего по сцепленным ребрам плоского кубического графа. Из доказанной теоремы как следствие вытекает утверждение о раскраске плоского графа четырьмя красками. Введена операция ротации (вращения) цветного диска. Операция ротации цветного диска позволяет построить алгоритм рекурсивной раскраски ребер плоского кубического графа. Это в свою очередь, позволяет производить раскраску вершин плоского графа G. Представлен метод и алгоритм раскраски плоского кубического графа и определена его вычислительная сложность. Приведены примеры раскраски негамильтоновых плоских кубических графов, что является своеобразной проверкой метода и алгоритма раскраски.

Показана тесная связь раскраски ребер и раскраски граней плоского кубического графа. Раскраска ребер плоского кубического графа порождает (индуцирует) раскраску его граней. А раскраска граней порождает (индуцирует) раскраску ребер плоского кубического графа.

В работе рассматриваются свойства раскраски биквадратного плоского графа являющегося реберным графом для плоского кубического. Показана общая картина раскраски плоских графов на основании раскраски ребер плоского кубического графа. Теория вращений введенная Г.Рингелем [12] позволяет не только выявлять закономерности плоской раскраски графов, но и создавать топологические рисунки графов, не проводя геометрических построений на плоскости.

**Библиографические ссылки**

1. Донец Г.А. Алгебраический подход к проблеме раскраски плоских графов / Г.А. Донец, Н.З. Шор - Київ: Наук.думка.- 1982. - 144с.

2. Донец Г.А. Алгоритмы раскраски плоских графов / Г.А. Донец // Теорія оптимальних рішень, номер 5, 2006. – С.136-144.

3. Гроссман И. Группы и их графы / И. Гроссман, В. Магнус – М.: Мир, 1971. – 247 с.

4. Гэри М. Вычислительные машины и труднорешаемые задачи: Пер. с англ. - / М. Гэри, Д.М. Джонсон - М.: Мир, 1982. - 416с.

5. Зыков А.А. Теория конечных графов. / А.А. Зыков - Новосибирск: ГРФМЛ.- 1963.- 542 с.

6. Зыков А.А. Основы теории графов. / А.А. Зыков - М.: Наука, ГРФМЛ, 1987.- 384с.

7. Емеличев В.А. Лекции по теории графов / В.А. Емеличев, О.И. Мельников, В.И. Сарванов, Р.И. Тышкевич - М.: Наука. ГРФМЛ, 1990. - 384с.

8. Липский В. Комбинаторика для программистов: Пер. с польск. / В. Липский - М.: Мир, 1988.-213 с.

9. Курапов С.В., Давидовский М.В., Толок А.В.. Визуальный алгоритм раскраски плоских графов / Научная визуализация, 2018, том 10, номер 3, страницы 1 - 33, DOI: 10.26583/sv.10.3.01

10. Комбинаторное условие для плоских графов / С. Мак-Лейн // В кн.: Кибернетический сборник. Новая серия, 1970.-вып. 7.- С.68-77.

11. Рейнгольд Э. Комбинаторные алгоритмы, теория и практика / Э. Рейнгольд, Ю. Нивергельт, Н. Дер - М.: Мир.- 1980. - 480с.

12. Рингель Г. Теорема о раскраске карт / Г. Рингель - М.: Мир, 1977. - 126с.

13. Родионов В.В. Методы четырехцветной раскраски вершин плоских графов / В.В. Родионов - М.: КомКнига, 2005. – 48с.

14. Самохин А.В. Проблема четырех красок: Неоконченная история доказательства / А.В. Самохин // Соровский образовательный журнал, том 6, № 7, 2000. – С.91-96.

15. Свами М. Графы, сети и алгоритмы: Пер. с англ. / М. Свами, К. Тхуласираман - М.: Мир, 1984. - 455с.

16. Фляйшнер Г. М. Эйлеровы графы и смежные вопросы. / Г. М.Фляйшнер: Мир, 2002. – 335 с.

17. Харари Ф. Теория графов. - пер. с англ. Козырева В.П. / под ред. Гаврилова В.Г. / Ф. Харари - М.: Мир.- 1973, 300с.

18. Appel K. The solution of the four-color-map problem / K. Appel, W. Haken // Scientific American, October 1977. P. 108–121.

19. Appel K., Haken W. Every Planar Map / K. Appel, W. Haken Is Four Colorable. Contemporary

Mathematics. Providence (R.I.): Amer. Math Soc., 1989. Vol. 98. 308 p.
20. T. Kavitha et. al. Cycle bases in graphs characterization, algorithms, complexity, and applications. Computer Science Review, **3** (2009), 199-243.